\documentclass{amsart}
\usepackage{graphicx}
\usepackage{color}
\usepackage{amssymb}
\usepackage{esint}
\usepackage{tikz-cd}
\usepackage{array}
\usepackage{bbm}
\usetikzlibrary{matrix,arrows.meta}
\usepackage{bbding}

\newcommand\R{\mathbb R}

\newcommand\C{\mathbb C}
\newcommand\E{\mathbb E}
\newcommand\SSS{\mathbb S}
\newcommand\V{\mathbb V}

\newcommand\loc{{\text{\upshape loc}}}
\newcommand\negquad{\!\!\!\!}

\newcommand\Pot{{\mathcal{P}}}
\newcommand\Potc{{\mathcal{P}^c}}
\newcommand\Eul{{\mathcal{E}}}
\newcommand\Eulc{{\mathcal{E}^c}}
\newcommand\Lagr{{\mathcal{L}}}
\newcommand\Lagrn[1]{{{\mathcal{L}}^#1}}

\newcommand\PcLn[1]{{_{{\Potc\!\Lagr}^#1}}}
\newcommand\PL{{_{\Pot\!\Lagr}}}
\newcommand\PLn[1]{{_{\Pot\!\Lagr^#1}}}
\newcommand\EcL{{_{\Eulc\!\Lagr}}}
\newcommand\EcLn[1]{{_{{\Eulc\!\Lagr}^#1}}}
\newcommand\EL{{_{\Eul\!\Lagr}}}
\newcommand\ELn[1]{{_{{\Eul\!\Lagr}^#1}}}

\newcommand\LnPc[1]{{_{\Lagr^#1 \Potc}}}
\newcommand\LP{{_{\Lagr \Pot}}}
\newcommand\LnP[1]{{_{\Lagr^#1 \Pot}}}
\newcommand\LEc{{_{\Lagr \Eulc}}}
\newcommand\LE{{_{\Lagr \Eul}}}
\newcommand\LnE[1]{{_{\Lagr^#1 \Eul}}}
\newcommand\LnEc[1]{{_{\Lagr^#1 \Eulc}}}
\newcommand\PE{{_{\Pot \Eul}}}
\newcommand\PcEc{{_{\Potc\!\Eulc}}}
\newcommand\EP{{_{\Eul \Pot}}}
\newcommand\EcPc{{_{\Eulc\!\Potc}}}
\newcommand\vvect{\mathbf{v}}
\newcommand\rvect{\mathbf{r}}
\newcommand\svect{\mathbf{s}}
\newcommand\nuvect{\boldsymbol{\nu}}
\newcommand\phivect{\boldsymbol{\phi}}
\newcommand\fvect{\mathbf{f}}
\newcommand\gvect{\mathbf{g}}
\newcommand\hvect{\mathbf{h}}
\newcommand\kvect{\mathbf{k}}
\newcommand\xvect{\mathbf{x}}
\newcommand\yvect{\mathbf{y}}
\newcommand\wvect{\mathbf{w}}
\newcommand\uvect{\mathbf{u}}
\newcommand\avect{\mathbf{a}}
\newcommand\bvect{\mathbf{b}}
\newcommand\cvect{\mathbf{c}}
\newcommand\dvect{\mathbf{d}}
\newcommand\mvect{\mathbf{m}}
\newcommand\nvect{\mathbf{n}}
\newcommand\pvect{\mathbf{p}}
\newcommand\qvect{\mathbf{q}}

\DeclareMathOperator\Div{div}
\DeclareMathOperator\curl{curl}
\DeclareMathOperator\DivGamma{div_{_\Gamma}}
\DeclareMathOperator\Tr{Tr}

\numberwithin{equation}{section}
\newtheorem{thm}{Theorem}
\numberwithin{thm}{section}
\newtheorem*{thm*}{Theorem}
\newtheorem{cor}{Corollary}
\numberwithin{cor}{section}
\newtheorem{lem}{Lemma}
\numberwithin{lem}{section}
\newtheorem{prop}{Proposition}
\numberwithin{prop}{section}
\theoremstyle{definition}
\newtheorem{definition}{Definition}
\numberwithin{definition}{section}
\theoremstyle{remark}
\newtheorem{rem}{Remark}
\numberwithin{rem}{section}

\begin{document}
\title[Acoustic waves interacting with non--locally...]
{Acoustic waves interacting with non--locally reacting surfaces in a Lagrangian framework}
\author{Enzo Vitillaro}
\address[E.~Vitillaro]
       {Dipartimento di Matematica e Informatica, Universit\`a di Perugia\\
       Via Vanvitelli,1 06123 Perugia ITALY}
\email{enzo.vitillaro@unipg.it}




\thanks{
This work has been funded by the European Union - NextGenerationEU within the framework of PNRR  Mission 4 - Component 2 - Investment 1.1 under the Italian Ministry of University and Research (MUR) programme "PRIN 2022" - 2022BCFHN2 - Advanced theoretical aspects in PDEs and their applications - CUP: J53D23003700006
}

\begin{abstract} The paper deals with a family of evolution problems arising in the physical modeling of small amplitude acoustic phenomena occurring in a fluid, bounded by a surface of extended reaction. They are all derived  in a Lagrangian framework.

We study well--posedness of these problems, their mutual relations, and their relations with other evolution problems modeling  the same physical phenomena. They are those introduced in an Eulerian framework and those which deal with the (standard in Theoretical Acoustics) velocity potential. The latter reduce to the well--known wave equation with acoustic boundary conditions.

Finally, we prove that all problems are asymptotically stable provided the system is linearly damped.
\end{abstract}

\maketitle

\section{Introduction and main results}\label{Section 1} \subsection{Presentation of the problem and literature overview}\label{Section 1.1}

We deal with a family of evolution problems arising in the physical modeling of small amplitude acoustic phenomena occurring in a fluid, bounded by a surface of extended reaction. Since all these problems are derived in a Lagrangian framework, we shall call them  {\em Lagrangian models}.

As we are going to see, solutions of all these problems are also solutions of the boundary value problem
$$(\Lagr)\qquad
\begin{cases}
\rho_0\rvect_{tt}-B\nabla \Div \rvect=0\qquad &\text{in
$\R\times\Omega$,}\\
\curl\rvect_t=0\qquad &\text{in
$\R\times\Omega$,}\\
\mu v_{tt}- \DivGamma (\sigma \nabla_\Gamma v)+\delta v_t+\kappa v-B\Div \rvect =0\qquad
&\text{on
$\R\times \Gamma_1$,}\\
\rvect_t\cdot{\boldsymbol{\nu}} =0 \quad \text{on $\R\times \Gamma_0$,}\qquad\rvect_t\cdot{\boldsymbol{\nu}} =-v_t\qquad
&\text{on
$\R\times \Gamma_1$,}
\end{cases}
$$
where $\Omega$ is a bounded and simply connected domain  of $\R^3$, with boundary $\Gamma=\partial\Omega$ of class
\begin{footnote}{see \cite[Def.~1.2.1.1,~Chap.~1,~p.5]{grisvard}}\end{footnote}
$C^{2,1}$, and $\Gamma=\Gamma_0\cup\Gamma_1$, $\overline{\Gamma_0}\cap \overline{\Gamma_1}=\emptyset$. We assume $\Gamma_1$ to be nonempty and connected.
Moreover,  $\rvect=\rvect(t,x)$,  $v=v(t,y)$, $t\in\R$, $x\in\Omega$, $y\in\Gamma_1$. Finally,  $v$  and $\rvect$
respectively are $\C$ and $\C^3$ -- valued.
\begin{footnote}{all unknowns and function spaces in the paper will be complex, as common in Acoustics, even if in the physical derivation of the problems one uses real unknowns. Being all problems linear with real valued coefficients, the reduction to the real case is anyway trivial.}\end{footnote}

In problem $(\Lagr)$, and in all problems in the sequel, we respectively denote
the gradient, divergence and rotor operators by $\nabla$, $\Div$ and $\curl$.  Furthermore,  $\DivGamma$ and $\nabla_\Gamma$ stand for the
Riemannian versions of the first two operators on $\Gamma$. Moreover, we denote the outward normal to $\Omega$ by
${\boldsymbol{\nu}}$.

Furthermore,  $B$ and $\rho_0$ are fixed positive constants, and
$\mu,\sigma,\delta$ and  $\kappa$ are given real  functions on $\Gamma_1$
such that $\mu,\sigma,\delta,\kappa\in W^{1,\infty}(\Gamma_1)$, $\min_{\Gamma_1}\mu=\mu_0>0$, $\min_{\Gamma_1}\sigma=\sigma_0>0$, $\kappa\ge 0$ \text on $\Gamma_1$.
These properties will be assumed throughout the paper without further reference.

To the author's knowledge, problem ($\Lagr$) does not possesses  a mathematical literature. However, it is a variant of another boundary--value problem belonging to the family we are introducing. It  is the problem
$$(\Lagrn 0)\qquad
\begin{cases}
\rho_0\rvect_{tt}-B\nabla \Div \rvect=0\qquad &\text{in
$\R\times\Omega$,}\\
\curl\rvect=0\qquad &\text{in
$\R\times\Omega$,}\\
\mu v_{tt}- \DivGamma (\sigma \nabla_\Gamma v)+\delta v_t+\kappa v-B\Div \rvect =0\qquad
&\text{on
$\R\times \Gamma_1$,}\\
\rvect\cdot{\boldsymbol{\nu}} =0 \quad \text{on $\R\times \Gamma_0$,}\qquad\rvect\cdot{\boldsymbol{\nu}} =-v\qquad
&\text{on
$\R\times \Gamma_1$,}
\end{cases}
$$
which was studied by the author in the recent paper \cite{tremodelli}.
Its physical derivation  was given  in \cite[Chapter~7]{mugnvit}, although it was
written there in the alternative form
$$(\widetilde{\Lagr}^0)\qquad
\begin{cases}
p_\Lagr+B\Div \rvect =0 \qquad &\text{in
$\R\times\Omega$,}\\
\rho_0\rvect_{tt}=-\nabla p_\Lagr\qquad &\text{in
$\R\times\Omega$,}\\
\curl\rvect=0\qquad &\text{in
$\R\times\Omega$,}\\
\mu v_{tt}- \DivGamma (\sigma \nabla_\Gamma v)+\delta v_t+\kappa v+p_\Lagr =0\qquad
&\text{on
$\R\times \Gamma_1$,}\\
\rvect\cdot{\boldsymbol{\nu}} =0 \quad \text{on $\R\times \Gamma_0$,}\qquad\rvect\cdot{\boldsymbol{\nu}} =-v\qquad
&\text{on
$\R\times \Gamma_1$,}
\end{cases}
$$
which  trivially reduces to $(\Lagrn 0)$ when eliminating $p_\Lagr$.

We refer to \cite[Chapter~7]{mugnvit} for details on the physical model.
Here we simply recall that $\rvect$ stands for the displacement of fluid particles, $v$ stands for the normal displacement inside $\Omega$ of the moving part $\Gamma_1$ of the boundary. Moreover, in $(\widetilde{\Lagr}^0)$, $p_\Lagr$ stands for the acoustic excess pressure. We also recall that  equation $(\widetilde{\Lagr}^0)_1$ is derived from the Hooke's Law, equations $(\widetilde{\Lagr}^0)_2$ and
$(\widetilde{\Lagr}^0)_4$ are derived from the Newton Second Law, while equations $(\widetilde{\Lagr}^0)_3$ and
$(\widetilde{\Lagr}^0)_5$ are respectively derived from $(\Lagr)_2$ and
$(\Lagr)_4$ by integrating in time. Although such an integration has physical motivations, as explained in \cite[Chapter 7]{mugnvit}, clearly $(\Lagr)$  looks more adherent to the original physical model than problem $(\Lagrn 0)$.

The family of problems we are considering includes, together with the already introduced problems $(\Lagr)$ and $(\Lagrn 0)$, four of their variants. Although  $(\Lagr)$ and $(\Lagrn 0)$ are the most important problems, we are going to see that these further  four problems are  intermediate between them. They  help understanding the relation between $(\Lagr)$ and $(\Lagrn 0)$, which could otherwise remain hidden.

They are the following boundary value problems:
$$
(\Lagrn 1)\qquad
\begin{cases}
\rho_0\rvect_{tt}-B\nabla \Div \rvect=0\qquad &\text{in
$\R\times\Omega$,}\\
\curl\rvect_t=0\qquad &\text{in
$\R\times\Omega$,}\\
\mu v_{tt}- \DivGamma (\sigma \nabla_\Gamma v)+\delta v_t+\kappa v-B\Div \rvect =0\qquad
&\text{on
$\R\times \Gamma_1$,}\\
\rvect\cdot{\boldsymbol{\nu}} =0 \quad \text{on $\R\times \Gamma_0$,}\qquad\rvect\cdot{\boldsymbol{\nu}} =-v\qquad
&\text{on
$\R\times \Gamma_1$;}
\end{cases}
$$
$$
(\Lagrn 2)\qquad
\begin{cases}
\rho_0\rvect_{tt}-B\nabla \Div \rvect=0\qquad &\text{in
$\R\times\Omega$,}\\
\curl\rvect=0\qquad &\text{in
$\R\times\Omega$,}\\
\mu v_{tt}- \DivGamma (\sigma \nabla_\Gamma v)+\delta v_t+\kappa v-B\Div \rvect =0\qquad
&\text{on
$\R\times \Gamma_1$,}\\
\rvect_t\cdot{\boldsymbol{\nu}} =0 \quad \text{on $\R\times \Gamma_0$,}\qquad\rvect_t\cdot{\boldsymbol{\nu}} =-v\qquad
&\text{on
$\R\times \Gamma_1$,}\\
\int_\Gamma\rvect\cdot\nu+\int_{\Gamma_1}v=0 & \text{in $\R$;}
\end{cases}
$$
$$
(\Lagrn 3)\qquad
\begin{cases}
\rho_0\rvect_{tt}-B\nabla \Div \rvect=0\qquad &\text{in
$\R\times\Omega$,}\\
\curl\rvect_t=0\qquad &\text{in
$\R\times\Omega$,}\\
\mu v_{tt}- \DivGamma (\sigma \nabla_\Gamma v)+\delta v_t+\kappa v-B\Div \rvect =0\qquad
&\text{on
$\R\times \Gamma_1$,}\\
\rvect_t\cdot{\boldsymbol{\nu}} =0 \quad \text{on $\R\times \Gamma_0$,}\qquad\rvect_t\cdot{\boldsymbol{\nu}} =-v\qquad
&\text{on
$\R\times \Gamma_1$,}\\
\int_\Gamma\rvect\cdot\nu+\int_{\Gamma_1}v=0 & \text{in $\R$;}
\end{cases}
$$
$$(\Lagrn 4)\qquad
\begin{cases}
\rho_0\rvect_{tt}-B\nabla \Div \rvect=0\qquad &\text{in
$\R\times\Omega$,}\\
\curl\rvect=0\qquad &\text{in
$\R\times\Omega$,}\\
\mu v_{tt}- \DivGamma (\sigma \nabla_\Gamma v)+\delta v_t+\kappa v-B\Div \rvect =0\qquad
&\text{on
$\R\times \Gamma_1$,}\\
\rvect_t\cdot{\boldsymbol{\nu}} =0 \quad \text{on $\R\times \Gamma_0$,}\qquad\rvect_t\cdot{\boldsymbol{\nu}} =-v\qquad
&\text{on
$\R\times \Gamma_1$.}
\end{cases}
$$
Also problems $(\Lagrn 1)$ -- $(\Lagrn 4)$ do not possess, to the author's knowledge, a mathematical literature.

Summarizing, we are then going to study the following family of Lagrangian models:
\begin{equation}\label{family}
\mathfrak L=\{(\Lagr), (\Lagrn 0), (\Lagrn 1),(\Lagrn 2),(\Lagrn 3),(\Lagrn 4)\}.
\end{equation}
We point out that all problems in $\mathfrak L$  can be easily obtained as follows:
one adds the boundary conditions given by the following Table to equations $(\Lagr)_1$ and $(\Lagr)_3$.
\begin{center}
\setlength{\arrayrulewidth}{0.6pt}
\begin{tabular}{|>{\scriptsize}c|>{\scriptsize}c|>{\scriptsize}c|>{\scriptsize}c|}
\multicolumn{4}{>{\normalsize}c}{\sc Table. }\\
\hline
 &
 $\begin{aligned}
 & \rvect\cdot \boldsymbol{\nu}=0\quad\text{on $\R\times\Gamma_0$}\\
 & \rvect\cdot \boldsymbol{\nu}=-v\quad\text{on $\R\times\Gamma_1$}
 \end{aligned}$
 &
 $\begin{aligned}
 &\\
 &\rvect_t\cdot \boldsymbol{\nu}=0\quad\text{on $\R\times\Gamma_0$}\\
 & \rvect_t\cdot \boldsymbol{\nu}=-v_t\quad\text{on $\R\times\Gamma_1$}\\
 &\int_\Gamma\rvect\cdot\boldsymbol{\nu}+\int_{\Gamma_1}v=0\\
 &
 \end{aligned}$
  &
  $\begin{aligned}
 &\rvect_t\cdot \boldsymbol{\nu}=0\quad\text{on $\R\times\Gamma_0$}\\
 & \rvect_t\cdot \boldsymbol{\nu}=-v_t\quad\text{on $\R\times\Gamma_1$}
 \end{aligned}$
 \\
\hline
$\curl\rvect=0$ in $\R\times\Omega$
 &$(\Lagrn 0)$&$(\Lagrn 2)$ & $(\Lagrn 4)$\\
\hline
$\curl\rvect_t=0$ in $\R\times\Omega$
 &$(\Lagrn 1)$&$(\Lagrn 3)$ &$(\Lagr)$\\
\hline
\end{tabular}
\end{center}
 The Lagrangian models are strongly connected with other problems, which  were studied in \cite{tremodelli} and physically derived in \cite[Chapter~7]{mugnvit}. The first one is  \emph{the Eulerian model}
$$(\Eul)\qquad
\begin{cases}
p_t+B\Div \vvect =0 \qquad &\text{in
$\R\times\Omega$,}\\
\rho_0\vvect_t=-\nabla p\qquad &\text{in
$\R\times\Omega$,}\\
\curl\vvect=0\qquad &\text{in
$\R\times\Omega$,}\\
\mu v_{tt}- \DivGamma (\sigma \nabla_\Gamma v)+\delta v_t+\kappa v+p =0\qquad
&\text{on
$\R\times \Gamma_1$,}\\
\vvect\cdot{\boldsymbol{\nu}} =0 \quad \text{on $\R\times \Gamma_0$,}\qquad\vvect\cdot{\boldsymbol{\nu}} =-v_t\qquad
&\text{on
$\R\times \Gamma_1$,}
\end{cases}
$$
where $p=p(t,x)$, $\vvect=\vvect(t,x)$, $v=v(t,y)$, $t\in\R$, $x\in\Omega$,
$p$ and $v$  take complex  values, $\vvect$ takes values in $\C^3$. Moreover, $v$ has the same meaning it had  in the Lagrangian models, while $p$ stands for the acoustic excess pressure and $\vvect$ for the velocity field, in an Eulerian framework.
The reader should remark the difference between the notations $\vvect$ and $v$.

We shall also deal with the \emph{ constrained Eulerian model}, in short  ($\Eulc$). It is constituted by adding the  integral condition
\begin{equation}\label{1.1}
  \int_\Omega p=B\int_{\Gamma_1}v\qquad\text{in  $\R$,}
\end{equation}
to  ($\mathcal{E})$. To the author's knowledge the Eulerian models $(\Eul)$ and $(\Eulc)$ possess no mathematical literature, but for the paper \cite{tremodelli}.

The Lagrangian and Eulerian models presented above are strongly connected with the
\emph{potential model}. It is the boundary value problem
$$(\Pot)\qquad
\begin{cases} \rho_0u_{tt}-B\Delta u=0 \qquad &\text{in
$\R\times\Omega$,}\\
\mu v_{tt}- \DivGamma (\sigma \nabla_\Gamma v)+\delta v_t+\kappa v+\rho u_t =0\qquad
&\text{on
$\R\times \Gamma_1$,}\\
\partial_{\boldsymbol{\nu}} u=0 \quad \text{on $\R\times \Gamma_0$,}\qquad v_t =\partial_{\boldsymbol{\nu}} u\qquad
&\text{on
$\R\times \Gamma_1$,}
\end{cases}
$$
where $u=u(t,x)$,  $t\in\R$, $x\in\Omega$, $u$ takes complex values and $v$ is as in previous models. In $(\Pot)$,  the symbol $\Delta$ stands for the Laplacian operator with respect to the space variable $x$.

We shall also consider a constrained version of problem $(\Pot)$,   called the
\emph{constrained potential model}, in short  ($\mathcal{P}^c$). It is obtained  by adding  the  integral condition
\begin{equation}\label{1.2}
  \rho_0\int_\Omega u_t=B\int_{\Gamma_1}v\qquad\text{in  $\R$,}
\end{equation}
to  ($\mathcal{P})$.

Unlike  the Lagrangian and Eulerian models, problem ($\mathcal{P}$) possesses a wide mathematical literature, in which the wave equation ($\mathcal{P}$)$_1$ is equivalently written as $u_{tt}-c^2\Delta u=0$, with $c^2=B/\rho_0$. It is known as wave equation with acoustic boundary conditions.

This type of  boundary conditions have been introduced by Beale and Rosencrans, for bounded or external domains, in \cite{beale, beale2,BealeRosencrans} when $\sigma\equiv 0$ and $\Gamma_0=\emptyset$. In this case the boundary $\Gamma=\Gamma_1$ is called, using the terminology of  \cite[pp.~256]{morseingard}, {\em locally reacting}, since each point of it reacts like an harmonic oscillator.

The same model, in the case $\delta=0$, has been proposed in \cite{WeitzKeller} when $\Omega$ is a strip in $\R^2$ or $\R^3$, and in \cite{Peters}
when $\Omega$ is the half-space in $\R^3$, to describe acoustic wave propagation in an ice-covered ocean. We refer to~\cite{Belinsky} for a historical overview of these and related problems in Mathematical Physics.

After their introduction acoustic boundary conditions for locally reacting surfaces have been the subject of a huge literature. See for example \cite{Alcantara, Alcantara2025, Bautista, CFL2004,CavFrota, FL2006, FG2000, GGG,  graber, Guo2023, Hao, JGPhD, JGSB2012, KJR2016, KT2008, LLX2018,  Maatoug2017,  mugnolo, Shomberg, Vicente2023}.
When one dismisses the simplifying assumption that neighboring point do not interact,  such surfaces, again using the terminology of \cite[p.266]{morseingard}), are called of {\em extended reaction}.
We shall call those which react like a membrane {\em non-locally reacting}, since other types of reactions can be considered.

The simplest case in which $\sigma\equiv\sigma_0$ was briefly considered in \cite[\S 6]{beale}. In this case
the operator $\DivGamma(\sigma \nabla_\Gamma)$ reduces, up to $\sigma_0$, to the Laplace--Beltrami operator $\Delta_\Gamma$.
 Subsequently, it was studied in \cite{becklin2019global,FMV2011, FMV2014, VF2013bis,VF2016, VF2017}.
In all of them the authors assume that $\Gamma_0\not=\emptyset$ and that the homogeneous Neumann boundary condition on it is replaced by the (mathematically more attracting) homogeneous Dirichlet boundary condition. See also \cite{AMP,ArElst} for somehow related problems.

Problem $(\Pot)$ was widely studied, under much more general assumptions, in the recent booklet \cite{mugnvit}. Well--posedness of problems $(\Eul)$, $(\Eulc)$, $(\Pot)$,  $(\Potc)$ and $(\Lagrn 0)$, and the relations among all these problems,  were studied in the recent paper  \cite{tremodelli} by the author. The present paper constitutes, in some sense, a prosecution of this study.

Disregarding problem $(\Lagrn 0)$, which will be considered in the sequel, we now recall, from \cite{tremodelli}, some results concerning the Eulerian and potential models.
We then introduce some related notation.
Denoting for simplicity $H^k(\Omega)^3:=H^k(\Omega;\C^3)$, for $k=1,2$, we introduce
its closed subspace
\begin{equation}\label{1.3}
 H^k_{\curl 0}(\Omega)=\{\vvect\in H^k(\Omega)^3: \curl \vvect =0\},
\end{equation}
where, at least when $k=0$, $\curl\vvect$ is taken in the sense of distributions.
We also introduce, for $k=1,2$,  the Fr\'echet  spaces
\begin{equation}\label{1.3bis}
\left\{\begin{aligned}
&X_\Pot^k=\bigcap_{i=0}^k C^i(\R;H^{k-i}(\Omega)\times H^{k-i}(\Gamma_1)),\\
&X_\Potc^k=\{(u,v)\in X_\Pot^k: \text{\eqref{1.2} holds}\},\\
&X_\Eul^k\!=\!\bigcap_{i=0}^{k-1}\! C^i(\R;H^{k-1-i}(\Omega))\!\times\! \bigcap_{i=0}^{k-1} \!C^i(\R;H^{k-1-i}_{\curl 0}(\Omega))\!\times\! \bigcap_{i=0}^k \!C^i(\R;H^{k-i}(\Gamma_1)),\\
&X_\Eulc^k=\{(p,\vvect,v)\in X_\Eul^k: \text{\eqref{1.1} holds}\}.
\end{aligned}\right.
\end{equation}
In \cite{tremodelli}  weak and strong solutions of  problems $(\Eul)$, $(\Eulc)$, $(\Pot)$,  $(\Potc)$ were found. Referring to Definitions~\ref{Definition 4.1} and \ref{Definition 4.2} below for details,  weak solutions of the problems listed above respectively are $(p,\vvect,v)\in X^1_\Eul$, $(p,\vvect,v)\in X^1_\Eulc$, $(u,v)\in X^1_\Pot$, and $(u,v)\in X^1_\Potc$ satisfying these problems in suitable distributional forms.

Strong solutions of $(\Eul)$, $(\Eulc)$, $(\Pot)$,  $(\Potc)$ are weak solutions of them  respectively belonging to the spaces
$X^2_\Eul$, $X^2_\Eulc$, $X^2_\Pot$, $X^2_\Potc$. They also satisfy the respective problems  a.e. in $\R\times\Omega$  and on $\R\times\Gamma$.

In the sequel we shall denote the Fr\'echet  spaces constituted by the weak solutions of these problems respectively by
$\mathcal{S}^1_\Eul $, $\mathcal{S}^1_\Eulc$, $\mathcal{S}^1_\Pot $ and $\mathcal{S}^1_\Potc$. Those constituted by their
strong solutions are the spaces
\begin{equation}\label{1.3bisbis}
\mathcal{S}^2_\Eul=\mathcal{S}^1_\Eul\cap X_\Eul^2, \quad
\mathcal{S}^2_\Eulc=\mathcal{S}^1_\Eulc\cap X_\Eulc^2, \quad
\mathcal{S}^2_\Pot=\mathcal{S}^1_\Pot\cap X_\Pot^2, \quad
\mathcal{S}^2_\Potc=\mathcal{S}^1_\Potc\cap X_\Potc^2.
\end{equation}
Moreover, we shall denote by $\C_{\R\times\Omega}$ the space of complex--valued constant functions in $\R\times\Omega$,  and
we set $\C_{X_\Pot}:=\C_{\R\times\Omega}\times\{0\}$. Since trivially $\C_{X_\Pot}\subseteq\mathcal{S}^2_\Potc\subseteq \mathcal{S}^2_\Pot\cap \mathcal{S}^1_\Potc\subseteq\mathcal{S}^1_\Pot$, we can then introduce, for $k=1,2$, the quotient Fr\'echet spaces (see \cite[Chapter I, p.31]{Schaefer}), $\dot{\mathcal{S}}^k_\Pot:=\mathcal{S}^k_\Pot/\C_{X_\Pot}$ and $\dot{\mathcal{S}}^k_\Potc:=\mathcal{S}^k_\Potc/\C_{X_\Pot}$.

Denoting by $\mathcal{L}(X;Y)$ the space of continuous linear operators between two Fr\'echet spaces $X$ and $Y$, and $\mathcal{L}(X)=\mathcal{L}(X;X)$, we can now briefly recall the relations occurring among the problems above.
\begin{thm*}[{\cite[Theorems~1.7~and~5.11]{tremodelli}}]
There is a surjective operator $\Psi^1_\PE\in\mathcal{L}(\mathcal{S}^1_\Pot;\mathcal{S}^1_\Eul)$,
defined by
\begin{equation}\label{1.3ter}
\Psi^1_\PE(u,v)=(p,\vvect,v):=(\rho_0 u_t,-\nabla u,v)\qquad\text{for all $(u,v)\in \mathcal{S}^1_\Pot$}.
\end{equation}
It  restricts to  the surjective operators
\begin{equation}\label{1.3quater}
\Psi^1_\PcEc\in\mathcal{L}(\mathcal{S}^1_\Potc;\mathcal{S}^1_\Eulc), \quad
\Psi^2_\PE\in\mathcal{L}(\mathcal{S}^1_\Pot;\mathcal{S}^1_\Eul), \quad\text{and}\quad
\Psi^2_\PcEc\in\mathcal{L}(\mathcal{S}^2_\Potc;\mathcal{S}^2_\Eulc).
\end{equation}
Since $\text{Ker } \Psi^k_\PE = \text{Ker } \Psi^k_\PcEc =\C_{X_\Pot}$ for
$k=1,2$,   the operator $\Psi_\PE^1$ subordinates a bijective isomorphism $\dot{\Psi}^1_\PE\in\mathcal{L}(\dot{\mathcal{S}}^1_\Pot;\mathcal{S}^1_\Eul)$, which restricts to the bijective isomorphisms
\begin{equation}\label{1.3six}
\dot{\Psi}^1_\PcEc\in\mathcal{L}(\dot{\mathcal{S}}^1_\Potc;\mathcal{S}^1_\Eulc),\quad
\dot{\Psi}^2_\PE\in\mathcal{L}(\dot{\mathcal{S}}^2_\Pot;\mathcal{S}^2_\Eul),\quad
\dot{\Psi}^2_\PcEc\in\mathcal{L}(\dot{\mathcal{S}}^2_\Potc;\mathcal{S}^2_\Eulc).
\end{equation}
\end{thm*}
Denoting $\dot{\Psi}^k_\EP=(\dot{\Psi}^k_\PE)^{-1}$ and $\dot{\Psi}^k_\EcPc=(\dot{\Psi}^k_\PcEc)^{-1}$, the following diagram
graphically shows the relations occurring among those problems:
\tikzcdset{row sep/normal=4em}
\tikzcdset{column sep/normal=8em}
\begin{equation}\label{1.4}
\begin{tikzcd}
[column sep=large]
\dot{\mathcal{S}}^k_\Pot\arrow[r,shift left,"\dot{\Psi}^k_\PE"]\arrow[r,leftarrow,shift right,"\dot{\Psi}^k_\EP"']&\mathcal{S}^k_\Eul
\end{tikzcd}
\qquad\qquad
\begin{tikzcd}
\dot{\mathcal{S}}^k_\Potc\arrow[r,shift left,"\dot{\Psi}^k_\PcEc"]\arrow[r,leftarrow,shift right,"\dot{\Psi}^k_\EcPc"']
&\mathcal{S}^k_\Eulc.
\end{tikzcd}
\end{equation}
In this sense problems $(\Eul)$ and  $(\Eulc)$ are  respectively equivalent to  $(\Pot)$ and $(\Potc)$. In the quoted paper we also proved that $(\Eulc)$ and $(\Potc)$ are equivalent to $(\Lagrn 0)$, while $(\Eul)$ and $(\Pot)$  have no Lagrangian counterpart.

This lack thereof essentially motivates the present study, which has the following aims:
\renewcommand{\labelenumi}{{\arabic{enumi})}}
\begin{enumerate}
\item showing the nontrivial structural relations among the Lagrangian models, which consitute the main tool to achieve all other goals;
\item extending the well--posedness theory in \cite{tremodelli} to all Lagrangian models;
\item showing that problems $(\Lagrn n)$, $n=0,1,2,3$, are equivalent to $(\Eulc)$ and $(\Potc)$, while
$(\Lagrn 4)$ and $(\Lagr)$ are equivalent to $(\Eul)$ and $(\Pot)$;
\item  proving that all Lagrangian models are asymptotically stable, when they are dissipative, that is when $\delta\ge 0$ and $\delta\not\equiv 0$.
\end{enumerate}

\subsection{Main results I: well--posedness.}\label{Section 1.2}
In the sequel we shall consider the initial--value problem associated to problem $(\Lagr)$,  i.e. the initial and boundary value problem
$$(\Lagr_0 )\qquad
\begin{cases}
\rho_0\rvect_{tt}-B\nabla \Div \rvect=0\qquad &\text{in
$\R\times\Omega$,}\\
\curl\rvect_t=0\qquad &\text{in
$\R\times\Omega$,}\\
\mu v_{tt}- \DivGamma (\sigma \nabla_\Gamma v)+\delta v_t+\kappa v-B\Div \rvect =0\qquad
&\text{on
$\R\times \Gamma_1$,}\\
\rvect_t\cdot{\boldsymbol{\nu}} =0 \quad\text{on $\R\times \Gamma_0$,}\qquad
\rvect_t\cdot{\boldsymbol{\nu}} =-v_t\qquad
&\text{on
$\R\times \Gamma_1$,}\\
\rvect(0,x)=\rvect_0(x),\quad \rvect_t(0,x)=\rvect_1(x) &
 \text{in $\Omega$,}\\
v(0,x)=v_0(x),\quad v_t(0,x)=v_1(x) &
 \text{on $\Gamma_1$.}
\end{cases}
$$
We shall also deal with the initial value problem associated with $(\Lagrn n)$, for $n=0,1,2,3,4$. It is  obtained by adding the initial conditions $(\Lagr_0 )_5$--$(\Lagr_0 )_6$ to $(\Lagrn n)$, and  will be denoted by  $(\Lagr_0^n)$.

As we are going to show, the only difference among problems in $\mathfrak L$ consists in the configuration spaces associated to them,  i.e. in the functional spaces in which the couple $(\rvect(t),v(t))$ lives for all $t\in\R$.
We are now going to introduce these configuration spaces, denoting by $H^1_Q$ the Hilbert space associated with a problem $(Q)$ in
the family. We want to make the reader aware that the same notational convention, i.e. a subscript "$Q$", will be used for all $Q$-- related spaces which also have an Eulerian or a potential counterpart. When such a counterpart is missing we shall just use, as subscript, the number associated to the problem.

We hence introduce the Hilbert space $H^1_\Lagr:=H^1(\Omega)^3\times H^1(\Gamma_1)$, the functional $L\in (H^1_\Lagr)'$ given by
\begin{equation}\label{1.5}
L(\svect,z)=\int_\Gamma\svect\cdot\boldsymbol{\nu}+\int_{\Gamma_1}z=\int_\Omega \Div\svect +\int_{\Gamma_1}z\quad\text{for all $(\svect,z)\in H^1_\Lagr$,}
\end{equation}
where the Divergence Theorem was used, and the further Hilbert spaces
\begin{equation}\label{1.6}
\begin{cases}
H^1_{\Lagrn 4}=&H^1_{\curl 0}\times H^1(\Gamma_1),\\
H^1_{\Lagrn 3}=&\{(\svect,z)\in H^1_\Lagr: L(\svect,z)=0\}, \\
H^1_{\Lagrn 2}=&\{(\svect,z)\in H^1_{\Lagrn 4}: L(\svect,z)=0\}, \\
H^1_{\Lagrn 1}=&\{(\svect,z)\in H^1_\Lagr: \svect\cdot\nuvect=0\,\,\text{on $\Gamma_0$},\,\,  \svect\cdot\nuvect=-z\quad\text{on $\Gamma_1$}\}, \\
H^1_{\Lagrn 0}=&\{(\svect,z)\in H^1_{\Lagrn 4}: \svect\cdot\nuvect=0\,\text{on $\Gamma_0$},\,  \svect\cdot\nuvect=-z\quad\text{on $\Gamma_1$}\}.
\end{cases}
\end{equation}
They  are closed subspaces of  $H^1_\Lagr$ and will be equipped with the inherited norm. By using the Divergence Theorem, the following inclusions
\begin{equation}\label{1.7}
\begin{alignedat}5
&H^1_{\Lagrn 0} && \subseteq && H^1_{\Lagrn 2} && \subseteq && H^1_{\Lagrn 4}
\\
&\phantom{\,\,}\rotatebox[origin=c]{270}{$\subseteq$}&& && \phantom{\,\,}\rotatebox[origin=c]{270}{$\subseteq$} && && \phantom{\,\,}\rotatebox[origin=c]{270}{$\subseteq$}
\\
&H^1_{\Lagrn 1} && \subseteq && H^1_{\Lagrn 3} && \subseteq && H^1_{\Lagr}
\end{alignedat}
\end{equation}
hold. They constitute the fundamental hierarchy among most of the spaces associated to the problems in $\mathfrak{L}$ and among solutions of them.

The time derivatives of solutions $(\rvect,v)$ of the problems in $\mathfrak{L}$ belong to a common Hilbert space, that is to
$H^0:=H^0_{\curl 0}(\Omega)\times L^2(\Gamma_1)$. Consequently, the phase spaces associated with problems in $\mathfrak L$  are
\begin{equation}\label{1.8}
\mathcal{H}^1_\Lagr=H^1_\Lagr\times H^0,\qquad\text{and}\quad \mathcal{H}^1_\Lagrn n=H^1_\Lagrn n\times H^0,\qquad\text{for $n=0,1,2,3,4$.}
\end{equation}
To simplify the notation, in the sequel we shall identify
$$\mathcal{H}^1_\Lagr=H^1(\Omega)^3\times H^1(\Gamma_1)\times H^1_{\curl 0}(\Omega)\times L^2(\Gamma_1).$$
According with this identification, we shall also identify,  for $n=0,1,2,3,4$, $\mathcal{H}^1_\Lagrn n$ with the corresponding subspace of the four--components space above.

In the paper we shall also deal with highly regular solutions. We then introduce  the Hilbert spaces
\begin{equation}\label{1.10}
H^2_\Lagr=H^2(\Omega)^3\times H^2(\Gamma_1),\quad\text{and}\quad H^2_\Lagrn n=H^2_\Lagr\cap H^1_{\Lagrn n},\quad\text{for $n=0,1,2,3,4$.}
\end{equation}
They are all endowed with the standard norm of the product  $H^2(\Omega)^3\times H^2(\Gamma_1)$. The corresponding phase spaces are
\begin{equation}\label{1.11}
\mathcal{H}^2_\Lagr=H^2_\Lagr\times H^1_\Lagrn 0,\qquad\text{and}\quad \mathcal{H}^2_\Lagrn n=H^1_\Lagrn n\times H^1_\Lagrn 0,\qquad\text{for $n=0,1,2,3,4$.}
\end{equation}
 In Section~\ref{Section 2.3} we shall  define weak and strong solutions of $(\Lagr)$ and of $(\Lagrn n  )$, $n=0,1,2,3,4$.
Weak solutions of them respectively  belong to the Fr\'echet  spaces
\begin{equation}\label{1.9}
X^1_\Lagr :=\{(\rvect,v)\in C(\R;H^1_\Lagr): (\rvect,v)'\in C(\R;H^0)\},\,
X^1_\Lagrn n:=X^1_\Lagr\cap C(\R:H^1_\Lagrn n).
\end{equation}
 They are solutions in a suitable distributional sense. Their formal  definition, see Definition~\ref{Definition 2.1} below, is an essential outcome of the paper. Initial conditions in problems $(\Lagr_0)$, $(\Lagr^n_0)$ will be meant in the continuity sense given by the space $X^1_\Lagr$.

 Strong solutions of $(\Lagr)$ and of $(\Lagrn n  )$, $n=0,1,2,3,4$, are  weak solutions of the corresponding problem which also respectively  belong to the Fr\'echet  spaces
\begin{equation}\label{1.12}
\begin{aligned}
&X^2_\Lagr :=\{(\rvect,v)\in C(\R;H^2_\Lagr): (\rvect,v)'\in C(\R;H^1_{\Lagrn 0})\cap C^1(\R;H^0)\},\\
&X^2_\Lagrn n:=X^2_\Lagr\cap C(\R:H^2_\Lagrn n).
\end{aligned}
\end{equation}
All of them are equipped with the standard topology of $X^2_\Lagr$. Strong solutions also satisfy the respective problem a.e. in $\R\times\Omega$ and on $\R\times\Gamma$.

We can give our first main result.
\begin{thm}[\bf Well--posedness]\label{Theorem 1.1}
 For all data  $U_0=(\rvect_0,v_0,\rvect_1,v_1)\in \mathcal{H}_\Lagr^1$,
 problem $(\Lagr_0)$ has a unique weak solution $(\rvect,v)\in X_\Lagr^1$,
 continuously depending on  $U_0$ in the topologies of the respective spaces.

Furthermore, setting  the energy function $\mathcal{E}=\mathcal{E}(\rvect,v)\in C(\R)$
 by
\begin{footnote}{here and in the sequel $|\cdot|_\Gamma$ denotes the norm  associated to the  metric naturally induced by the Euclidean metric on the tangent bundle, see \$~\ref{Section 2.2}  below.}\end{footnote}
\begin{equation}\label{1.13}
\mathcal{E}(\rvect,v)= \tfrac 12  \int_\Omega \!\!\!\rho_0|\rvect_t|^2+ B |\Div\rvect|^2
+\tfrac 12 \int_{\Gamma_1}\!\!\!\sigma |\nabla_\Gamma v|_\Gamma^2+\mu |v_t|^2
+\kappa|v|^2,
\end{equation}
$(\rvect,v)$ satisfies the energy identity
\begin{equation}\label{1.14}
\mathcal{E}(\rvect,v)(t)-\mathcal{E}(\rvect,v)(s) =-\int_s^t\!\!\int_{\Gamma_1}\!\!\delta|v_t|^2\quad\text{for all  $s,t\in\R$.}
\end{equation}
Moreover, for $n=0,1,2,3,4$, $(\rvect,v)$ is the unique weak solution of problem $(\Lagrn n)$ if and only if $(\rvect,v)\in X_\Lagrn n^1$ and, in turn, if and only if
$U_0\in \mathcal{H}_\Lagrn n^1$. Consequently, for all $U_0\in \mathcal{H}_\Lagrn n^1$, problem $(\Lagrn n)$ has a unique weak solution continuously depending on   $U_0$ and satisfying the energy identity \eqref{1.14}.

Furthermore, $(\rvect,v)\in X_\Lagr^2$ if and only if  $U_0\in \mathcal{H}_\Lagr^2$, such data being dense in $\mathcal{H}_\Lagr^1$.
In this case $\rvect$ and $v$ satisfy equations $(\Lagr)_1$ -- $(\Lagr)_2$
a.e. in $\R\times\Omega$ and equations $(\Lagr)_3$--$(\Lagr)_4$ a.e. on $\R\times\Gamma$. Moreover, $(\rvect,v)$
continuously depends  on $U_0$ in the topologies of $\mathcal{H}_\Lagr^2$ and $X_\Lagr^2$.

Finally, for $n=0,1,2,3,4$, we can respectively replace $X^2_\Lagr$, $\mathcal{H}^2_\Lagr$ and  $\mathcal{H}^1_\Lagr$ with
$X^2_\Lagrn n$, $\mathcal{H}^2_\Lagrn n$ and  $\mathcal{H}^1_\Lagrn n$ in the last statement..
\end{thm}
\begin{rem}\label{Remark 1.1} As far as problem $(\Lagr_0^0)$ is concerned, Theorem~\ref{Theorem 1.1} is nothing but \cite[Theorem~1.3]{tremodelli}. This part of the statement of Theorem~\ref{Theorem 1.1} was repeated only for the sake of completeness.
Theorem~\ref{Theorem 1.1} extends the quoted result to all Lagrangian models.

For the sake of simplicity, in the present paper we shall not study further regularity of solutions, as done in \cite{tremodelli}.
However, the interested reader can extend \cite[Theorem~4.18]{tremodelli}, after formulating the related compatibility conditions.
\end{rem}
\begin{rem}\label{Remark 1.2} Theorem~\ref{Theorem 1.1} shows that the hierarchy \eqref{1.7} also applies to the whole Lagrangian models, since solution of a problem are solutions of another one provided the phase space of the former is included in the phase space of the latter. Hence, in particular, problems $(\Lagrn 0)$ and $(\Lagr)$  respectively are the least and the most general ones, and all problems in $\mathfrak{L}$ can be considered as restrictions of $(\Lagr)$.
\end{rem}
The proof of Theorem~\ref{Theorem 1.1} is completely different from the proofs of \cite[Theorem~1.3]{tremodelli} and \cite[Theorem~1.2.1]{mugnvit}. Indeed, while in the former the relation between problems $(\Lagrn 0)$ and $(\Potc)$ is used, using it for other Lagrangian models looks problematic. Indeed  this relation is more involved,  as we shall see.

Furthermore, while in the latter a semigroup approach is used, repeating a similar argument for other Lagrangian models looks problematic as well. Indeed, even the quasi--monotonicity of the operator that arises from an abstract formulation of the problems $(\Lagrn 4)$ and $(\Lagr)$ is not evident. Checking the range condition for this operator looks nontrivial as well.

Actually the proof of Theorem~\ref{Theorem 1.1} given in Section~\ref{Section 4} relies on combining \cite[Theorem~1.3]{tremodelli} with the structural decomposition result that we shall present in the sequel.
\subsection{Main results II: structural decomposition.} \label{Section 1.3}
Our second main result deals with equilibria and stationary solutions of the Lagrangian models. These two subtly different concepts are usually identified in the mathematical literature. By a \emph{weak} (\emph{strong}) \emph{equilibrium} of one of the Lagrangian models we mean an element of its configuration space such that the function constantly having this value in  $\R$ is as weak (strong) solution of the problem. Such a solution is then called a \emph{weak} (\emph{strong}) \emph{stationary solution} of the problem. See also Definition~\ref{Definition 2.2} below.

In the sequel we shall denote, for $k=1,2$,
\begin{equation}\label{2.10}
\begin{aligned}
\E^k=&\{(\svect,z)\in H^k_\Lagr: \text{$(\svect,z)$ is a weak equilibrium of $(\Lagr)$}\},
\\
\SSS^k=&\{(\rvect,v)\in X^k_\Lagr: \text{$(\rvect,v)$ is a weak stationary solution of $(\Lagr)$}\}.
\end{aligned}
\end{equation}
Since strong solutions of $(\Lagr)$  are exactly the weak solutions of it belonging to $X^2_\Lagr$, the spaces $\E^1$ and $\SSS^1$ above
respectively are the spaces of weak equilibria and stationary solutions, while the spaces  $\E^2$ and $\SSS^2$ are their strong counterparts.

We shall also denote by $i\in\mathcal{L}(\E^k;\SSS^k)$, $k=1,2$, the trivial bijective isomorphism given  by
\begin{equation}\label{1.15}
i(\svect,z)(t)=(\svect,z)\qquad\text{for all $t\in\R$ and $(\svect,z)\in \E^k$.}
\end{equation}
Moreover, we shall denote, for $k=1,2$ and $n=0,1,2,3,4$,
\begin{equation}\label{2.11}
\begin{aligned}
\E^k_n=&\{(\svect,z)\in H^k_\Lagrn n : \text{$(\svect,z)$ is a weak equilibrium of $(\Lagrn n)$}\},
\\
\SSS^k_n=&\{(\rvect,v)\in X^k_\Lagrn n: \text{$(\rvect,v)$ is a weak stationary solution of $(\Lagrn  n)$}\}.
\end{aligned}
\end{equation}
The same comments made above for  $\E^k$ and $\SSS^k$ apply to $\E^k_n$ and $\SSS^k_n$. Indeed, by Theorem~\ref{Theorem 1.1},
for $k=1,2$ and $n=0,1,2,3,4$, we have
\begin{equation}\label{2.12}
\SSS^k_n=\SSS^k\cap X^k_\Lagrn n =  i\E^k_n \quad\text{and}\quad \E^k_n=\E^k\cap H^k_\Lagrn n.
\end{equation}
To state our  main result we have to preliminarily introduce a particular strong equilibrium for problems $(\Lagrn 4)$ and $(\Lagr)$. Respectively denoting by $|\Omega|$, $\mathcal{H}^2(\Gamma_1)$, and $\mathbbm{1}_{\Gamma_1}$ the volume of $\Omega$,  the area of $\Gamma_1$, and the function taking constant value $1$ on $\Gamma_1$, this equilibrium takes two different forms, depending on the alternative $\kappa\equiv 0$ or $\kappa\not\equiv 0$:
\renewcommand{\labelenumi}{{\roman{enumi})}}
\begin{enumerate}
\item when $\kappa\equiv 0$, the couple $(0,\mathbbm{1}_{\Gamma_1})$
is a strong equilibrium of $(\Lagrn 4)$ and $(\Lagr)$;
\item when $\kappa\not\equiv 0$, we shall prove that the equation
\begin{equation}\label{1.17}
  -\Div_\Gamma(\sigma \nabla_\Gamma z^*)+\kappa z^*+B=0\qquad\text{on $\Gamma_1$}
\end{equation}
and  the problem
\begin{equation}\label{1.18}
\begin{cases}
\Div \svect^*=-1\quad&\text{in $\Omega$,}\\
\curl \svect^*=0\quad&\text{in $\Omega$,}\\
\svect^*\cdot\nuvect=0\quad&\text{on $\Gamma_0$,}\\
\svect^*\cdot\nuvect=-|\Omega|/\mathcal{H}^2(\Gamma_1)\quad&\text{on $\Gamma_1$,}
\end{cases}
\end{equation}
respectively have a unique solution $z^*\in H^2(\Gamma_1)$ and $\svect^*\in H^2_{\curl 0}(\Omega)$. Then $(\svect^*,z^*)$ is a strong equilibrium  of $(\Lagrn 4)$ and $(\Lagr)$.
\end{enumerate}
To treat  the two cases  above at once, we  set the strong equilibrium
\begin{equation}\label{1.19}
(\svect^\bullet ,z^\bullet)=
\begin{cases}
(0,\mathbbm{1}_{\Gamma_1}),\qquad&\text{when $\kappa\equiv 0$,}
\\
(\svect^*,z^*),\qquad&\text{when $\kappa\not\equiv 0$,}
\end{cases}
\end{equation}
 of $(\Lagrn 4)$ and $(\Lagr)$, and we respectively denote by $\E^\bullet$ and $\SSS^\bullet$ the one--dimensional spaces of strong equilibria and stationary
solutions which are generated by $(\svect^\bullet,z^\bullet)$ and  $i(\svect^\bullet,z^\bullet)$, that is
\begin{equation}\label{1.20}
\E^\bullet=\C (\svect^\bullet,z^\bullet),\qquad\text{and}\quad \SSS^\bullet=i\E^\bullet.
\end{equation}
In Lemma~\ref{Lemma 3.3} we shall recognize that $L(\svect^\bullet,z^\bullet)\not=0$, and as a consequence:
\begin{itemize}
\item[-]$(\svect^\bullet,z^\bullet)$ is not an equilibrium of the problem $(\Lagrn n )$ when $n=0,1,2,3$. Indeed, by \eqref{1.6} and \eqref{1.7}, one has $(\svect^\bullet,z^\bullet)\not\in H^1_\Lagrn n$;
\item[-] we can set the functional $\ell\in (H^1_\Lagr)'$   defined by
\begin{equation}\label{1.21}
\ell(\svect,z)=\frac{L(\svect,z)}{L(\svect^\bullet,z^\bullet)}=
\begin{cases}
\dfrac{\int_\Gamma \svect\cdot\nuvect+\int_{\Gamma_1}z}{\mathcal{H}^2(\Gamma_1)},\quad &\text{when $\kappa\equiv 0$,}\\
\dfrac{\int_\Gamma \svect\cdot\nuvect+\int_{\Gamma_1}z}{-|\Omega|+\int_{\Gamma_1}z^*},\quad &\text{when $\kappa\not\equiv 0$.}
\end{cases}
\end{equation}
\end{itemize}
The next result characterizes the equilibria  of  the Lagrangian models.
\begin{thm}[\bf Equilibria]\label{Theorem 1.2}
For $n=0,1,2,3$, the weak and strong equilibria of problem $(\Lagrn n)$  are given by the formula
\begin{equation}\label{1.22}
 \E^k_0=\{0\},\quad \E^k_n=\V^k_n\times \{0\},\quad\text{for $k=1,2$,}
\end{equation}
where
\begin{equation}\label{1.23}
\begin{cases}
\V^k_1:=&\{\svect\in H^k(\Omega)^3: \Div\svect=0\quad\text{and}\quad \svect\cdot\nuvect=0\quad\text{on $\Gamma$}\},\\
\V^k_2:=&\{\svect\in H^k_{\curl 0}(\Omega): \Div\svect=0\},\\
\V^k_3:=&\{\svect\in H^k(\Omega)^3: \Div\svect=0\}.
\end{cases}
\end{equation}
The weak and strong equilibria of problem $(\Lagrn 4)$ and $(\Lagr)$  are given by
\begin{equation}\label{1.24}
  \E^k_4=\E^k_2\oplus\E^\bullet \quad\text{and}\quad \E^k=\E^k\oplus\E^\bullet.
\end{equation}
\end{thm}
\begin{rem}\label{Remark 1.3}Theorem~\ref{Theorem 1.2} shows that, while the problem $(\Lagrn 0)$ (studied in \cite{tremodelli})  possesses only the trivial equilibrium, all other Lagrangian models possess many nontrivial equilibria, hence enjoying a richer structure.

More in detail, since one trivially has $\{\curl\svect: \svect\in \mathcal{D}(\Omega)^3\}\subseteq \V^2_1\subseteq \V^1_1$, both $\V^1_1$ and $\V^2_1$ are infinite dimensional, as well as $\E^1_1$ and $\E^2_1$.

Moreover, since it is well--known (see \cite[Chapter~IX]{dautraylionsvol3}) that, for $k=1,2$, one has
$\V^k_2=\{\nabla\phi: \phi\in H^{k+1}(\Omega)\,\,\text{and}\,\,\Delta\phi=0\}$, the same remark applies to the spaces
$\V^k_2$ and $\E^k_2$.

By \eqref{1.22}--\eqref{1.24}, one has $\E^k_2\subseteq \E^k_n\subset \E^K$ for $n=3,4,5$, so (in conclusion)  the remark applies to the spaces $\E^k$ and $\E^k_n$ for $n=1,2,3,4$.
\end{rem}

\begin{rem}\label{Remark 1.4}
One trivially has $\E^2\subseteq \E^1$ and $\E^2_n\subseteq \E^1_n$ for $n=0,1,2,3,4$, since strong equilibria are also weak ones. Actually, all these inclusions are proper, but for the case $\E^2_0=\E^1_0=\{0\}$.
Indeed (see Remark~\ref{Remark 3.1} below) one has $\E^2_n\not=\E^1_n$ for $n=1,2,3$. So, by \eqref{1.24}, this conclusion extends to $\E^2_4\not= \E^1_4$ and $\E^2\not=\E^1$.
This fact explains why distinguish between strong and weak equilibria is mandatory.
\end{rem}
\begin{rem}\label{Remark 1.5}
Theorem~\ref{Theorem 1.2} shows that $\mathfrak{L}$ can be conveniently divided into the following two subfamilies:
\begin{equation}\label{1.24bis}
\mathfrak{A}:=\{(\Lagrn 0), (\Lagrn 1), (\Lagrn 2), (\Lagrn 3)\},\qquad\text{and}\quad
\mathfrak{B}:= \{(\Lagrn 4), (\Lagr)\}.
\end{equation}
Equilibria of problems in $\mathfrak{A}$ have vanishing boundary deformation $z$, while elements of $\E^\bullet$ can have $z$ of arbitrarily large amplitude. Hence elements of $\E^\bullet$, and then most equilibria of problems in $\mathfrak{B}$, are sharply contrasting the assumptions made when building the physical model. According with them, both $\svect$ and $z$ have to be small. Furthermore, equilibria for problems in $\mathfrak{A}$ have divergence free displacement $\svect$ and vanishing energy $\mathcal{E}$, while one can take  equilibria in $\E^\bullet$  with arbitrarily large energy and  displacement.
  \end{rem}
  To state our next main result we introduce the Fr\'{e}chet spaces constituted by  weak and strong solutions of the Lagrangian models, that is for $k=1,2$ and $n=0,1,2,3,4$,
  \begin{equation}\label{1.25}
  \begin{alignedat}2
  &\mathcal{S}_\Lagr^k=&&\{(\rvect,v)\in X^k_\Lagr:\,\text{$(\rvect,v)$ is a weak solution of $(\Lagr)$}\}
  \\
  &\mathcal{S}_\Lagrn n^k=&&\{(\rvect,v)\in X^k_\Lagrn n:\,\,\text{$(\rvect,v)$ is a weak solution of $(\Lagrn n)$}\}=
  \mathcal{S}_\Lagr^k\cap X^l_\Lagrn n.
  \end{alignedat}
  \end{equation}
  Since strong solutions of $(\Lagr)$ are exactly weak solutions belonging to $X^2_\Lagr$, solutions in the spaces defined by \eqref{1.25}
  are strong if and only if  $k=2$.

Our third main result shows that differences among the Lagrangian models  only consist in their different equilibria.
\begin{thm}[\bf Structural decomposition]\label{Theorem 1.3}
For $k=1,2$ and $n=0,1,2,3,4$, we have
\begin{equation}\label{1.26}
  H^k_\Lagr=H^k_\Lagrn 0\oplus\E^k, \quad\text{and}\quad H^k_\Lagrn n =H^k_\Lagrn 0\oplus\E^k_n.
\end{equation}
The projectors $\Pi_{H^k_\Lagr;H^k_\Lagrn 0}\in\mathcal{L}(H^k_\Lagr;H^k_\Lagrn 0)$ and
$\Pi_{H^k_\Lagr;\E^k}\in\mathcal{L}(H^k_\Lagr;\E^k)$, associated with the first decomposition, are given, for any $(\svect,z)\in H^k_\Lagr$, by
\begin{footnote}{these projectors, which will be relevant in  Theorem~\ref{Theorem 1.6} below, have a more explicit form when distinguishing between the cases $\kappa\equiv 0$ and $\kappa\not\equiv 0$,  given in Appendix~\ref{appendice A}.}\end{footnote}
\begin{equation}\label{1.27}
 \Pi_{H^k_\Lagr;\E^k}(\svect,z)=(\fvect,\ell(\svect,z)z^\bullet),\quad\text{and}\quad
 \Pi_{H^k_\Lagr;H^k_\Lagrn 0}(\svect,z)=(\gvect,z-\ell(\svect,z)z^\bullet),
\end{equation}
where $\fvect, \gvect\in H^k(\Omega)^3$ respectively are the unique solutions of the  problems
$$
\boxed{(F)
\begin{cases}
\Div \fvect=\ell(\svect,z)\Div \svect^\bullet&\text{in $\Omega$,}\\
\curl \fvect=\curl\svect&\text{in $\Omega$,}\\
\fvect\cdot\nuvect=\svect\cdot\nuvect&\text{on $\Gamma_0$,}\\
\fvect\cdot\nuvect=\svect\cdot\nuvect+z-\ell(\svect,z)z^\bullet\negquad\!\!&\text{on $\Gamma_1$,}
\end{cases}
}
\,\,
\boxed{(G)
\begin{cases}
\Div \gvect=\Div\svect-\ell(\svect,z)\Div \svect^\bullet\!\!\negquad&\text{in $\Omega$,}\\
\curl \gvect=0&\text{in $\Omega$,}\\
\gvect\cdot\nuvect=0&\text{on $\Gamma_0$,}\\
\gvect\cdot\nuvect=-z+\ell(\svect,z)z^\bullet&\text{on $\Gamma_1$.}
\end{cases}
}
$$
Consequently,  for all data $(\rvect_0,v_0,\rvect_1,v_1)\in H^k_\Lagr$, the solution $(\rvect,v)$ of problem $(\Lagr_0)$, given by Theorem~\ref{Theorem 1.1}, admits the  decomposition
\begin{equation}\label{1.28}
 (\rvect(t),v(t))=\Pi_{H^k_\Lagr;\E^k}(\rvect_0,v_0)+(\rvect^0(t),v^0(t)),\qquad\text{for all $t\in\R$,}
\end{equation}
where $(\rvect^0,v^0)\in X^k_\Lagrn 0$ is the solution of problem $(\Lagr_0^0)$ (also given by Theorem~\ref{Theorem 1.1}) corresponding to initial data $(\rvect^0_0,v^0_0,\rvect_1,v_1)\in H^k_\Lagrn 0$, and $(\rvect_0^0,v_0^0):=\Pi_{H^k_\Lagr;H^k_\Lagrn 0}(\rvect_0,v_0)$.

Moreover, for any $n=0,1,2,3,4$, the projectors $\Pi_{H^k_\Lagrn n;H^k_\Lagrn 0}\in\mathcal{L}(H^k_\Lagrn n;H^k_\Lagrn 0)$ and
$\Pi_{H^k_\Lagrn n;\E^k_n}\in\mathcal{L}(H^k_\Lagrn n;\E^k_n)$, associated with the second decomposition in \eqref{1.26}, respectively are the restrictions of those in \eqref{1.27}.

Furthermore, the decomposition \eqref{1.28} continues to hold when respectively replacing $(\Lagr_0 )$, $H^k_\Lagr$, $\mathcal{H}^k_\Lagr$ and $\E^k$ with $(\Lagr^n_0 )$, $H^k_\Lagrn n$, $\mathcal{H}^k_\Lagrn n$ and $\E^k_n$.

Finally, the further decompositions hold true:
\begin{equation}\label{1.29}
\mathcal{S}^k_\Lagr=\mathcal{S}^k_\Lagrn 0\oplus \SSS^k\quad\text{and}\quad
\mathcal{S}^k_\Lagrn n=\mathcal{S}^k_\Lagrn 0\oplus \SSS^k_n\quad\text{for $n=0,1,2,3,4$.}
\end{equation}
\end{thm}
The proof of Theorem~\ref{Theorem 1.3} essentially reduces to proving the decompositions \eqref{1.26} and explicitly giving the projectors associated to them. It is worth pointing out that this goal is best achieved by considering the whole family $\mathfrak{L}$, starting from problem $(\Lagrn 0)$  to  finally achieve the result for problem $(\Lagr)$.

\subsection{Main results II: relations.}  In the next results we shall use, for $k=1,2$, the following quotient Fr\'{e}chet spaces:
\begin{equation}\label{1.30}
 \dot{\mathcal{S}}^k_\Lagrn n:=\mathcal{S}^n_\Lagrn n /\SSS_n^k,\,\,\text{for $n=0,1,2,3$},\quad  \dot{\mathcal{S}}^k_\Lagrn 4=\mathcal{S}^k_\Lagrn 4/\SSS^k_2,\quad\text{and}\quad
 \dot{\mathcal{S}}^k_\Lagr=\mathcal{S}^k_\Lagr /\SSS^k_3.
\end{equation}
\begin{rem}\label{Remark 1.6}
Since the quotient spaces introduced above are crucial in describing the relations involving the Lagrangian models, it is worth making some remarks:
\begin{itemize}
\item[-] since $\SSS^k_0=\{0\}$ for $k=0,1$, the spaces $\dot{\mathcal{S}}^k_{\Lagrn 0}$ coincide with the spaces $\mathcal{S}^k_{\Lagrn 0}$; they have been introduced  only for the sake of notation's consistency;
\item[-] since $\SSS_n^2\not=\SSS_n^1$ for $n=1,2,3$, the space $\dot{\mathcal{S}}^2_\Lagr$ is \emph{not included} in the space
$\dot{\mathcal{S}}^1_\Lagr$ and  the same remarks applies to the spaces related to problem $(\Lagrn n)$;
\item[-] the spaces $\mathcal{S}^k_\Lagr$ and $\mathcal{S}^k_\Lagrn n$ obey, for any $k=1,2$, the hierarchy established by \eqref{1.7}, but previous remark shows that the quotient spaces above \emph{do not obey it};
\item[-] the quotient spaces above are isomorphic to two different spaces , depending on the belonging of their problem to $\mathfrak{A}$ or to $\mathfrak{B}$ in \eqref{1.24bis}; indeed, since by \eqref{1.24} we also have $\SSS_4^k=\SSS_2^k\oplus \SSS^\bullet$ and
    $\SSS^k=\SSS_3^k\oplus \SSS^\bullet$, for $k=1,2$, we get
    \begin{equation}\label{1.30bis}
      \dot{\mathcal{S}}^k_\Lagrn n\simeq \mathcal{S}^k_\Lagrn 0\quad\text{when $n=0,1,2,3$, and}\quad
      \dot{\mathcal{S}}^k_\Lagrn 4\simeq \dot{\mathcal{S}}^k_\Lagr\simeq \mathcal{S}^k_\Lagrn 0\oplus \SSS^\bullet.
    \end{equation}
    The space $\mathcal{S}^k_\Lagrn 0\oplus \SSS^\bullet$ is not the  solutions space of any Lagrangian model, in contrast with $\mathcal{S}^k_\Lagrn 0$, so making the subfamily $\mathfrak{B}$ interesting.
\end{itemize}
\end{rem}
We can now give our two next main results, highlighting the relations among the Lagrangian, Eulerian and potential models.

\begin{thm}[\bf Relations between Lagrangian and  Eulerian models]\label{Theorem 1.4}
For any weak (strong) solution  $(\rvect,v)$ of  $(\Lagr)$,  the triple
$(p,\vvect,v)$ identified by the equations
\begin{equation}\label{1.31}
p=-B\Div \rvect, \quad \text{and}\quad \vvect=\rvect_t,\quad\text{in $\R$,}
\end{equation}
is a weak (strong) solution of  $(\Eul)$. Moreover, $(p,\vvect,v)$ also solves $(\Eulc)$ provided $(\rvect,v)$ solves  $(\Lagrn n)$ for any
$n=0,1,2,3$.

Furthermore, fixing  $k=1,2$ and $n=0,1,2,3$,  the map $(\rvect,v)\mapsto(p,\vvect,v)$ defines the surjective operators $\Psi^k_\LE \in \mathcal{L}(\mathcal{S}^k_\Lagr;\mathcal{S}^k_\Eul)$, $\Psi^k_\LnE 4\in \mathcal{L}(\mathcal{S}^k_\Lagrn 4;\mathcal{S}^k_\Eul)$, and $\Psi^k_\LnEc n\in \mathcal{L}(\mathcal{S}^k_\Lagrn n;\mathcal{S}^k_\Eulc)$, all of them being restrictions of $\Psi^1_\LE$.

Moreover, we have $\text{Ker } \Psi^k_\LE=\SSS^k_3$, $\text{Ker } \Psi^k_\LnE 4=\SSS^k_2$, and $\text{Ker } \Psi^k_\LnEc n=\SSS^k_n$.
Consequently, the operators introduced above respectively subordinate the bijective isomorphisms $\dot{\Psi}^k_\LE \in \mathcal{L}(\dot{\mathcal{S}}^k_\Lagr;\mathcal{S}^k_\Eul)$, $\dot{\Psi}^k_\LnE 4\in \mathcal{L}(\dot{\mathcal{S}}^k_\Lagrn 4;\mathcal{S}^k_\Eul)$, and $\dot{\Psi}^k_\LnEc n\in \mathcal{L}(\dot{\mathcal{S}}^k_\Lagrn n;\mathcal{S}^k_\Eulc)$,
\begin{footnote}{in general these operators are not restrictions of $\dot{\Psi}^1_\LE$, since their domains are not subspaces of $\dot{\mathcal S}^1_\Lagr$. See Remark~\ref{Remark 1.6}.}\end{footnote}  given by
 \begin{equation}\label{1.32}
\left\{
\begin{alignedat}4
  &\dot{\Psi}^k_\LE  &&[(\rvect,v)+\SSS^k_3]=&&\Psi^k_\LE  (\rvect,v)\,\, &&\text{for all $(\rvect,v)\in\mathcal{S}_\Lagr^k$,}\\
  &\dot{\Psi}^k_\LnE 4 &&[(\rvect,v)+\SSS^k_2]=&&\Psi^k_\LnE 4 (\rvect,v)\,\, &&\text{for all $(\rvect,v)\in\mathcal{S}_\Lagrn 4^k$,}\\
  &\dot{\Psi}^k_\LnEc n &&[(\rvect,v)+\SSS^k_n]=\quad&&\Psi^k_\LnEc n (\rvect,v)\quad &&\text{for all $(\rvect,v)\in\mathcal{S}_\Lagrn n^k$.}
\end{alignedat}
\right.
\end{equation}
Their inverses  are respectively the operators $\dot{\Psi}^k_\EL\in \mathcal{L}(\mathcal{S}^k_\Eul;\dot{\mathcal{S}}^k_\Lagr)$,
$\dot{\Psi}^k_\ELn 4\in \mathcal{L}(\mathcal{S}^k_\Eul;\dot{\mathcal{S}}^k_\Lagrn 4)$, and  $\dot{\Psi}^k_\EcLn n\in \mathcal{L}(\mathcal{S}^k_\Eulc;\dot{\mathcal{S}}^k_\Lagrn n)$, trivially given by
\begin{equation}\label{1.33}
\left\{
\begin{alignedat}3
&\dot{\Psi}^k_\EL(p,\vvect,v)=&&\{(\rvect,v)\in \mathcal{S}^k_\Lagr : \text{\eqref{1.31} hold}\}/\SSS^k_3&&\quad\text{for all $(p,\vvect,v)\in\mathcal{S}_\Eul^k$,}\\
&\dot{\Psi}^k_\ELn 4(p,\vvect,v)=&&\{(\rvect,v)\in \mathcal{S}^k_\Lagrn 4 : \text{\eqref{1.31} hold}\}/\SSS^k_2&&\quad\text{for all $(p,\vvect,v)\in\mathcal{S}_\Eul^k$,}\\
&\dot{\Psi}^k_\EcLn n(p,\vvect,v)=&&\{(\rvect,v)\in \mathcal{S}^k_\Lagrn n : \text{\eqref{1.31} hold}\}/\SSS^k_n &&\quad\text{for all $(p,\vvect,v)\in\mathcal{S}_\Eulc^k$.}
\end{alignedat}
\right.
\end{equation}
The operator $\dot{\Psi}^k_\EcLn n$ can be further characterized as follows.  For any $(p,\vvect,v)\in \mathcal{S}^k_\Eulc$ and $t\in\R$, the problem
\begin{equation}\label{1.34}
  \begin{cases}
  -B\Div \dot{\rvect}(t)=p(t)\quad &\text{in $\Omega$,}\\
  \quad \curl\dot{\rvect}(t)=0\quad &\text{in $\Omega$,}\\
  \quad \dot{\rvect}(t)\cdot\boldsymbol{\nu}=0\quad &\text{on $\Gamma_0$,}\\
 \quad \dot{\rvect}(t)\cdot\boldsymbol{\nu}=-v(t)\quad &\text{on $\Gamma_1$,}
  \end{cases}
\end{equation}
has a unique solution $\dot{\rvect}(t)\in H^k_{\curl 0}(\Omega)$. The couple $(\dot{\rvect},v)$ coincides with the unique weak solution $(\rvect,v)$ of  $(\Lagrn 0)$ satisfying equations \eqref{1.31}. We have $\dot{\Psi}^k_\EcLn 0(p,\vvect,v)=(\dot{\rvect},v)$ and,
more generally,
$\dot{\Psi}^k_\EcLn n(p,\vvect,v)=(\dot{\rvect},v)+\SSS^k_n$ for $n=0,1,2,3$.
\end{thm}

\begin{thm}[\bf Relations between Lagrangian and  potential models]\label{Theorem 1.5}
For any weak (strong) solution  $(\rvect,v)$ of $(\Lagr)$,  we set $u$, up to a space--time constant, by
\begin{equation}\label{1.37}
u(t)=u(0)-\tfrac B{\rho_0}\int_0^t \Div \rvect(\tau)\,d\tau,\quad\text{and}\quad -\nabla u(0)=\rvect_t(0).
\end{equation}
Then the couple $(u,v)$ is a weak (strong) solution of $(\Pot)$,
satisfying the equations
\begin{equation}\label{1.38}
-B\Div \rvect =\rho_0 u_t,\quad \rvect_t=-\nabla u\quad\text{in $\R$,}
\end{equation}
which also solves $(\Potc)$ whenever $(\rvect,v)$ solves $(\Lagrn n)$ for some $n=0,1,2,3$.

Moreover, fixing $k=1,2$ and $n=0,1,2,3$, the map $(\rvect,v)\mapsto(u,v)+\C_{X_\Pot}$ defines the surjective operators $\Psi^k_\LP \in \mathcal{L}(\mathcal{S}^k_\Lagr;\dot{\mathcal{S}}^k_\Pot)$,  $\Psi^k_\LnP 4\in \mathcal{L}(\mathcal{S}^k_\Lagrn 4;\dot{\mathcal{S}}^k_\Pot)$, and $\Psi^k_\LnPc n\in \mathcal{L}(\mathcal{S}^k_\Lagrn n;\dot{\mathcal{S}}^k_\Potc)$,
all of them being restrictions of $\Psi^1_\LP$.

Furthermore, we have $\text{Ker } \Psi^k_\LP =\SSS^k_3$, $\text{Ker } \Psi^k_\LnP 4=\SSS^k_2$, and $\text{Ker } \Psi^k_\LnPc n=\SSS^k_n$. Consequently, the operators introduced above subordinate the bijective isomorphisms $\dot{\Psi}^k_\LP \in \mathcal{L}(\dot{\mathcal{S}}^k_\Lagr;\dot{\mathcal{S}}^k_\Pot)$, $\dot{\Psi}^k_\LnP 4\in \mathcal{L}(\dot{\mathcal{S}}^k_\Lagrn 4;\dot{\mathcal{S}}^k_\Pot)$, and  $\dot{\Psi}^k_\LnPc n\in \mathcal{L}(\dot{\mathcal{S}}^k_\Lagrn n;\dot{\mathcal{S}}^k_\Potc)$,
\begin{footnote}{in general these operators are not restrictions of $\dot{\Psi}^1_\LP$, as in the previous Theorem. }\end{footnote}
 given by
\begin{equation}\label{1.39}
\left\{
\begin{alignedat}4
  &\dot{\Psi}^k_\LP  &&[(\rvect,v)+\SSS^k_3]=&&\Psi^k_\LP  (\rvect,v)\,\, &&\text{for all $(\rvect,v)\in\mathcal{S}_\Lagr^k$,}\\
  &\dot{\Psi}^k_\LnP 4 &&[(\rvect,v)+\SSS^k_2]=&&\Psi^k_\LnP 4 (\rvect,v)\,\, &&\text{for all $(\rvect,v)\in\mathcal{S}_\Lagrn 4^k$,}\\
  &\dot{\Psi}^k_\LnPc n &&[(\rvect,v)+\SSS^k_n]=&&\Psi^k_\LnPc n (\rvect,v)\quad &&\text{for all $(\rvect,v)\in\mathcal{S}_\Lagrn n^k$.}
\end{alignedat}
\right.
\end{equation}
Their inverses  are respectively the operators $\dot{\Psi}^k_\PL\in \mathcal{L}(\dot{\mathcal{S}}^k_\Pot;\dot{\mathcal{S}}^k_\Lagr)$,
$\dot{\Psi}^k_\PLn 4\in \mathcal{L}(\dot{\mathcal{S}}^k_\Pot;\dot{\mathcal{S}}^k_\Lagrn 4)$, and $\dot{\Psi}^k_\PcLn n\in \mathcal{L}(\dot{\mathcal{S}}^k_\Potc;\dot{\mathcal{S}}^k_\Lagrn n)$, trivially  given by
\begin{equation}\label{1.40}
\left\{
\begin{alignedat}3
&\dot{\Psi}^k_\EL &&[(u,v)+\C_{X_\Pot}]=&&\quad\{(\rvect,v)\in \mathcal{S}^k_\Lagrn 4 : \text{\eqref{1.38} hold}\}/\SSS^k_3,\\
&\dot{\Psi}^k_\PLn 4 &&[(u,v)+\C_{X_\Pot}]=&&\quad\{(\rvect,v)\in \mathcal{S}^k_\Lagrn 4 : \text{\eqref{1.38} hold}\}/\SSS^k_2,\\
&\dot{\Psi}^k_\PcLn n &&[(u,v)+\C_{X_\Pot}]=&&\quad\{(\rvect,v)\in \mathcal{S}^k_\Lagrn n : \text{\eqref{1.38} hold}\}/\SSS^k_n,
\end{alignedat}
\right.
\end{equation}
for all $(u,v)+\C_{X_\Pot}$ belonging to the respective domain.

The operator $\dot{\Psi}^k_\PcLn n$ can be further characterized as follows. For any $(u,v)\in \dot{\mathcal{S}}^k_\Potc$ and $t\in\R$,  the problem
\begin{equation}\label{1.41}
  \begin{cases}
  -B\Div \dot{\rvect}(t)=\rho_0 u_t(t),\quad &\text{in $\Omega$,}\\
  \quad \curl\dot{\rvect}(t)=0\quad &\text{in $\Omega$,}\\
  \quad \dot{\rvect}(t)\cdot\boldsymbol{\nu}=0\quad &\text{on $\Gamma_0$,}\\
 \quad \dot{\rvect}(t)\cdot\boldsymbol{\nu}=-v(t)\quad &\text{on $\Gamma_1$,}
  \end{cases}
\end{equation}
has a unique solution $\dot{\rvect}(t)\in H^k_{\curl 0}(\Omega)$.
The couple $(\dot{\rvect},v)$ coincides with the unique weak solution $(\rvect,v)$ of  $(\Lagrn 0)$ satisfying equations \eqref{1.38}. Hence  $\dot{\Psi}^k_\PcLn 0[(u,v)+\C_{X_\Pot}]=(\dot{\rvect},v)$ and, more generally,
$\dot{\Psi}^k_\PcLn n[(u,v)+\C_{X_\Pot}]=(\dot{\rvect},v)+\SSS^k_n$ for $n=0,1,2,3$.
\end{thm}
\begin{rem}\label{Remark 1.7}The inverse of the isomorphism related to problem $(\Lagrn n)$, for $n=0,1,2,3$, is more explicit than those related to problems  $(\Lagrn 4)$ and $(\Lagr)$. Actually, it would be possible to give a more explicit form of the latter, still using $(\dot{\rvect},v)$, but this explicit form is involved and  does not seem of interest, hence it is omitted.
\end{rem}
The isomorphisms given in \eqref{1.4} and those given in
Theorems~\ref{Theorem 1.4}--\ref{Theorem 1.5} are best illustrated by the  following commutative diagrams. In them $k=1,2$ and $n=0,1,2,3$. In the first diagram problem $(\Lagr)$ can be systematically replaced by problem $(\Lagrn 4)$.

\tikzcdset{row sep/normal=6em}
\tikzcdset{column sep/normal=10em}
\begin{tikzcd}
\dot{\mathcal{S}}^k_\Pot\arrow[r,shift left,"\dot{\Psi}^k_\PE"]\arrow[r,leftarrow,shift right,"\dot{\Psi}^k_\EP"']
\arrow[d,shift left,"\dot{\Psi}^k_\LP"]\arrow[d,leftarrow,shift right,"\dot{\Psi}^k_\PL"']&\mathcal{S}^k_\Eul\\
\dot{\mathcal{S}}^k_\Lagr\arrow[ru,shift left,"\Psi^k_\LP"]\arrow[ru,leftarrow,shift right,"\Psi^k_\PL"']&
\end{tikzcd}
\begin{tikzcd}
\dot{\mathcal{S}}^k_\Potc\arrow[r,shift left,"\dot{\Psi}^k_\PcEc"]\arrow[r,leftarrow,shift right,"\dot{\Psi}^k_\EcPc"']
\arrow[d,shift left,"\dot{\Psi}^k_\PcLn n"]\arrow[d,leftarrow,shift right,"\dot{\Psi}^k_\LnPc n"']&\mathcal{S}^k_\Eulc\\
\dot{\mathcal{S}}^k_\Lagrn n\arrow[ru,shift left,"\dot{\Psi}^k_\LnPc n"]\arrow[ru,leftarrow,shift right,"\dot{\Psi}^k_\PcLn n"']&
\end{tikzcd}

\subsection{Main results III: asymptotic stability.}
Our final main result asserts the asymptotic stability of all Lagrangian models when some dissipation is present.
\begin{thm}[\bf Asymptotic stability]\label{Theorem 1.6}
Also assume that $\delta\ge 0$ on $\Gamma_1$ and $\delta\not\equiv 0$. For all data $(\rvect_0,v_0,\rvect_1,v_1)\in\mathcal{H}^1_\Lagr$, let $(\rvect,v)$ denote the unique solution of problem $(\Lagr_0)$ given by Theorem~\ref{Theorem 1.1}. Then  we have
\begin{equation}\label{1.43}
 (\rvect(t), v(t), \rvect_t(t), v_t(t))\to (\rvect_\infty, v_\infty,0,0)\quad\text{in $\mathcal{H}^1_\Lagr$,\, as $t\to\infty$,}
\end{equation}
where $(\rvect_\infty,v_\infty):=\Pi_{H^1_\Lagr;\E^1}(\rvect_0,v_0)$, the projector $\Pi_{H^1_\Lagr;\E^1}$ being given by \eqref{1.27}.

The asymptotic behavior given by \eqref{1.43} continues to hold when respectively replacing $(\Lagr_0)$, $\mathcal{H}^1_\Lagr$, $H^1_\Lagr$ and $\E^1$ with  $(\Lagr^n_0)$, $\mathcal{H}^1_\Lagrn n$, $H^1_\Lagrn n$ and $\E^1_n$ for $n=0,1,2,3,4$.

In particular, when $(\rvect_0,v_0,\rvect_1,v_1)\in\mathcal{H}^1_\Lagrn 0$, so $(\rvect,v)$ solves $(\Lagr_0^0)$, we have
$$(\rvect(t), v(t), \rvect_t(t), v_t(t))\to 0\quad\text{in $\mathcal{H}^1_\Lagr$, \,as $t\to\infty$.} $$
\end{thm}
The conclusion of Theorem~\ref{Theorem 1.6} is new also for problem  $(\Lagr_0^0)$, already considered in \cite{tremodelli}. Its proof essentially consists in combining Theorem~\ref{Theorem 1.3} with \cite[Theorem~1.3.1]{mugnvit} and  \cite[Theorem~1.5]{tremodelli}.

\noindent{\bf Organization of the paper.} The paper is organized as follows. In Section ~\ref{Section 2} we set the notation and give some  preliminaries. Section~\ref{Section 3} deals with our structural decomposition result, so including the proof of Theorem~\ref{Theorem 1.2} and formula \eqref{1.26}. The Section~\ref{Section 4} is devoted to prove all remaining main results.
\section{Notation and preliminaries}\label{Section 2}
\subsection{Notation}
Borrowing a convention in Physics, for vectors in $\C^3$ and vector--valued functions we shall use boldface. For $\mathbf{x}=(x_1,x_2,x_3)$, $\mathbf{y}=(y_1,y_2,y_3)\in \C^3$ we shall denote $\mathbf{x}\cdot \mathbf{y}=\sum_{i=1}^3 x_iy_i$ and  by $\overline{\mathbf{x}}$  the vector conjugated with $\mathbf{x}$.

We shall use the standard notation  for functions spaces on $\Omega$, referring to \cite{adamsfournier}. As already done in formula \eqref{1.7}, where  $H^k(\Omega)^3=H^k(\Omega;\C^3)$, to simplify the notation we shall systematically identify the $\C^3$--valued versions of all spaces above with the Cartesian cubes of the corresponding scalar spaces. Moreover $\|\cdot\|_p$, $1\le p\le \infty$, will denote the norm in $L^p(\Omega)$ and in $L^p(\Omega)^3$, since no confusion will arise.

Moreover, for any Fr\'{e}chet space $X$,  we shall denote by $X'$ its dual, by $\langle\cdot,\cdot\rangle_X$ the duality product.  When $X$ is a Banach space we shall use the standard notation for Bochner--Lebesgue and Sobolev spaces of $X$--valued functions.
\subsection{Function spaces and operators on $\Gamma$}\label{Section 2.2}The assumption made on $\Omega$, $\Gamma_0$, and $\Gamma_1$, assures that $\Gamma$ inherits from $\R^3$ the structure of a Riemannian surface of class $C^2$. Hence,  in the sequel, we shall use some notation of geometric nature. It is quite common in the smooth case, see \cite{taylor}, and  can be easily extended to the $C^2$ case. See for example \cite{mugnvit} or \cite{Dresda1,Dresda2}.

Moreover, since $\overline{\Gamma_0}\cap \overline{\Gamma_1}=\emptyset$, both $\Gamma_0$ and $\Gamma_1$ are relatively open on $\Gamma$. Hence all geometrical concepts apply to them as well. To avoid making repetition, in the sequel we shall denote by $\Gamma'$ any relatively open subset of $\Gamma$.

We shall denote by $(\cdot,\cdot)_\Gamma$ the Riemannian metric inherited from $\R^3$ and uniquely extended to an Hermitian metric on the complexified  tangent bundle $T(\Gamma')$, and also the associated bundles metric on the complexified cotangent bundle $T^*(\Gamma')$. By $|\cdot|_\Gamma^2= (\cdot,\cdot)_\Gamma$  we shall denote the associated bundle norms.

The standard surface elements $\omega$ associated to $(\cdot,\cdot)_\Gamma$ is  the density of the Lebesgue surface measure on $\Gamma$, coinciding with the restriction to $\Gamma$ of the Hausdorff measure $\mathcal{H}^2$, i.e. $\omega=d\mathcal{H}^2$. In the sequel $\Gamma$ will be equipped, without further comments,  with this measure and  the corresponding  notions of a.e. equivalence, integrals and Lebesgue spaces $L^p(\Gamma')$, $1\le p\le \infty$. For simplicity, we shall denote $\|\cdot\|_{p,\Gamma'}=\|\cdot\|_{L^p(\Gamma')}$. Moreover, the notation $d\mathcal{H}^2$ will be dropped from boundary integrals, and a.e. equivalence on $\R\times\Gamma$ will be referred to the Hausdorff measure $\mathcal{H}^3$ in $\R^4$.

Sobolev spaces on $\Gamma'$ are treated in many textbooks in the smooth case, see for example
\cite{hebey,lionsmagenes1}. The $C^2$ case  is treated in  \cite{grisvard} and, when $\Gamma$ is possibly non--compact, in \cite{mugnvit}. Here we shall refer, for simplicity, to \cite{grisvard}, and we shall use the standard notation.
Moreover, mainly to simplify the notation, since $\overline{\Gamma_0}\cap \overline{\Gamma_1}=\emptyset$, by identifying the elements of $W^{s,q}(\Gamma_i)$, $i=0,1$,  with their trivial extensions to $\Gamma$, we have the decomposition
\begin{equation}\label{2.01}
W^{s,q}(\Gamma)=W^{s,q}(\Gamma_0)\oplus W^{s,q}(\Gamma_1)\qquad\text{for  $-2\le s\le 2$ and $1\le q<\infty$.}
\end{equation}

We refer to \cite{mugnvit} for details on the  Riemannian gradient operator $\nabla_\Gamma$ and on the
Riemannian divergence operator $\Div_\Gamma$. Here we just recall that one gets the operator
\begin{equation}\label{2.02}
 \Div_\Gamma(\sigma\nabla_\Gamma)\in\mathcal{L}\left(H^2(\Gamma_1),L^2(\Gamma_1)\right),
\end{equation}
and that, $\Gamma_1$ being compact, one also gets the integration by parts formula on $\Gamma_1$
\begin{equation}\label{2.03}
- \int_{\Gamma_1} \Div_\Gamma(\sigma\nabla_\Gamma u)v=\int_{\Gamma_1}\sigma (\nabla_\Gamma u,\nabla_\Gamma \overline{v})_\Gamma
\end{equation}
for all  $u\in H^2(\Gamma_1)$, $v\in H^1(\Gamma_1)$.

Finally, in the sequel, we shall use the well--known Trace Theorem, i.e. the  existence of the trace operator
$\Tr\in \mathcal{L}\left(H^m(\Omega),H^{m-1/2}(\Gamma)\right)$ for $m=1,2$.
We shall denote, as usual, $\Tr u=u_{|\Gamma}$. By $u_{|\Gamma_0}$ and $u_{|\Gamma_1}$ we shall mean the restrictions of $u_{|\Gamma}$ to $\Gamma_0$ and $\Gamma_1$. When clear, we shall omit trace related subscripts.
\subsection{Solutions of problems $\boldsymbol{(\Lagr)}$, $\boldsymbol{(\Lagr_0)}$ and $\boldsymbol{(\Lagrn n)}$, $\boldsymbol{(\Lagr^n_0)}$. }\label{Section 2.3}
At first we make precise which types of solutions we shall consider in the sequel.
\begin{definition}\label{Definition 2.1}
We say that
\renewcommand{\labelenumi}{{\roman{enumi})}}
\begin{enumerate}
\item $(\rvect,v)\in X^2_\Lagr$ is a {\em strong solution} of $(\Lagr)$ provided $(\Lagr)_1$ holds a.e. in $\R\times\Omega$ and
$(\Lagr)_3$ --$(\Lagr)_4$ holds a.e. on $\R\times\Gamma$;
\begin{footnote}{by \eqref{1.6} and \eqref{1.12}, the equation $(\Lagr)_2$  is automatically verified.}\end{footnote}
\item $(\rvect,v)\in X^1_\Lagr$ is a {\em generalized solution} of $(\Lagr)$ provided it is the limit in $X^1_\Lagr$ of a sequence of strong solutions of it;
\item $(\rvect,v)\in X^1_\Lagr$ is a {\em weak solution} of $(\Lagr)$ provided the following distributional identities hold true:
\begin{alignat}2
\label{2.1}
  &\int_{\R\times\Omega} \curl \rvect\,\varphi_t=0 \quad&&\text{for all $\varphi\in C^\infty_c(\R\times\Omega)$,}\\
\label{2.2}
&\int_{\R\times\Gamma_0} \rvect\cdot\nuvect\, \psi_t=0 \quad&&\text{for all $\psi\in C^2_c(\R\times\Gamma_0)$,}\\
\label{2.3}
&\int_{\R\times\Gamma_1} (\rvect\cdot\nuvect+v)\, \psi_t=0 \quad&&\text{for all $\psi\in C^2_c(\R\times\Gamma_1)$,}
\end{alignat}
and
\begin{multline}\label{2.4}
 \int_{\R\times\Omega} \rho_0\rvect_t\cdot \phivect_t-B\Div \rvect\Div \phivect\\
 +\int_{\R\times\Gamma_1}\mu v_t\psi_t-\sigma(\nabla_\Gamma v,\nabla_\Gamma \overline{\psi})_\Gamma-(\delta v_t+\kappa)\psi
 =0
\end{multline}
for all $\phivect\in C^1_c(\R\times\R^3)^3$ such that $\phivect\cdot\boldsymbol{\nu}=0$ on $\R\times\Gamma_0$, where $\psi:=-\phivect\cdot\boldsymbol{\nu}$ on $\R\times\Gamma_1$;
\item $(\rvect,v)$ is a strong, generalized or weak solution \emph{of problem $(\Lagrn n)$, for $n=0,1,2,3,4$},  provided it is a solution of $(\Lagr)$ of the same type and it belongs to the space $X^1_\Lagrn n$;
\item $(\rvect,v)$ is a strong, generalized or weak solution \emph{of problem $(\Lagr_0)$ (or of problem $(\Lagr^n_0)$), for $n=0,1,2,3,4$}, provided is is a solution of $(\Lagr)$ (of problem $(\Lagrn n)$) and the initial conditions hold in the space $X^1_\Lagr$.
\end{enumerate}
\end{definition}
Since solutions of $(\Lagrn n)$ are also solutions of $(\Lagr)$, the following discussion  on solutions of the latter will concern the former as well.

Strong solutions are a.e. classical solutions, so deserving some attention. Unfortunately, living in the space $X^2_\Lagr$, strong solutions of problem $(\Lagrn 0)$ can exist only when $(\rvect_0,v_0,\rvect_1,v_1)\in \mathcal{H}_\Lagr^2$ (as stated in Theorem~\ref{Theorem 1.1}), i.e. data are regular enough and satisfy the compatibility conditions $\rvect_1\cdot\nuvect=0$ on $\Gamma_0$, $\rvect_1\cdot\nuvect=-v_1$ on $\Gamma_1$.

Hence we are leaded to consider solutions merely belonging to the energy space $X^1_\Lagr$, like the generalized ones. They keep important properties of strong solutions, like the energy identity \eqref{1.13}, but they satisfy $(\Lagr)$ in a quite indirect way. Characterizing them as weak solutions is then important, as it will be clear from the following  discussion on them.

Clearly equations \eqref{2.1} and \eqref{2.2}--\eqref{2.3} respectively are the distributional version of equations $(\Lagr)_2$ and
$(\Lagr)_4$. On the other hand, while for $(\rvect,v)\in X^1_\Lagr$ the equation $(\Lagr)_1$ assumes in $\mathcal{D}'(\R\times\Omega)^3$ the natural form
\begin{equation}\label{2.5}
 \int_{\R\times\Omega} \rho_0\rvect_t\cdot \phivect_t-B\Div \rvect\Div \phivect=0\qquad\text{for all $\phivect\in [\mathcal{D}(\R\times\Omega)]^3$,}
 \end{equation}
 equation $(\Lagr)_3$ can not be written in a distributional sense, unless the term $\Div\rvect$ in it, merely belonging to $C(\R;L^2(\Omega))$, has some trace sense on $\R\times\Gamma_1$.

 This type of difficulty was already found in \cite{tremodelli}, when dealing with problem $(\Lagrn 0)$. It was solved
 by  combining equations $(\Lagr)_1$ and $(\Lagr)_3$ in the single distributional identity \eqref{2.4}.

 The following result, which will be also useful in the sequel, shows that when $\Div\rvect$ has such a trace sense, then the distributional identity \eqref{2.4} is equivalent to the combination of \eqref{2.5} with the distributional version of the equation $(\Lagr)_3$, that is
 \begin{equation}\label{2.6}
 \int_{\R\times\Gamma_1}\mu v_t\psi_t-\sigma(\nabla_\Gamma v,\nabla_\Gamma \overline{\psi})_\Gamma-\delta v_t\psi
 -\kappa v\psi+B\Div\rvect\psi=0
 \end{equation}
 for all $\psi\in C^2_c(\R\times\Gamma_1)$.
 Hence the definition of weak solutions given above is the closest possible approximation of a notion of distributional solution of the whole problem $(\Lagr)$.
  \begin{prop}\label{Proposition 2.1} Let $(\rvect,v)\in X^1_\Lagr$ be such that $\Div\rvect\in L^1_\loc(\R; H^1(\Omega))$. Then
 $(\rvect,v)$ is a weak solution of $(\Lagr)$ if and only if it satisfies the distributional identities \eqref{2.1}--\eqref{2.3} and \eqref{2.5}--\eqref{2.6}.
  \end{prop}
 \begin{proof} Let $(\rvect,v)$ be as in the statement. We first claim that \eqref{2.5} holds if and only if
 \begin{equation}\label{2.7}
 \rvect_t\in W^{1,1}_\loc (\R;L^2(\Omega)^3)\quad\text{and}\quad \rho_0\rvect_{tt}=B\nabla\Div \rvect
 \quad\text{in $L^1_\loc (\R;L^2(\Omega)^3)$.}
 \end{equation}
 Indeed, if \eqref{2.7} holds, one gets \eqref{2.5}, simply  by multiplying the equation in \eqref{2.7} by $\phivect$ and integrating by parts in space and time. Conversely, if \eqref{2.5} holds, by integrating by parts in $\Omega$ one gets
 $$\int_{\R\times\Omega} \rho_0\rvect_t\cdot \,\phivect_t+B\nabla\Div \rvect\cdot\phivect=0\quad\text{for all $\phivect\in\mathcal{D}(\R\times\Omega)^3$.}$$
  Consequently, taking test functions in the separate form $\phivect(t,x)=\varphi(t)\phivect_0(x)$, with $\varphi\in\mathcal{D}(\R)$ and $\phivect_0\in \mathcal{D}(\Omega)^3$, we get
 \begin{equation}\label{2.8}
  \int_\Omega\left[\int_{-\infty}^\infty \rho_0\rvect_t\varphi'+B\nabla\Div \rvect \varphi\right]\cdot\phivect_0=0.
  \end{equation}
 Moreover, \eqref{2.8} trivially extends by density to $\phivect_0\in H^1_0(\Omega)^3$. Consequently, we get that $\rvect_t\in W^{1,1}_\loc(\R;H^{-1}(\Omega)^3)$ and that the equation in \eqref{2.7} holds true in $L^1_\loc(\R;H^{-1}(\Omega)^3)$.
  Since $\nabla\Div \rvect\in L^1_\loc(\R; L^2(\Omega)^3)$  and $L^2(\Omega)^3\simeq [L^2(\Omega)']^3\hookrightarrow H^{-1}(\Omega)^3$, we then get that $\rvect_t\in W^{1,1}_\loc (\R;L^2(\Omega)^3)$  and that the equation in \eqref{2.7} holds in the space $L^1_\loc(\R;L^2(\Omega)^3)$, so proving our claim.

 Now let $(\rvect,v)$ be a weak solution of $(\Lagr)$. Taking in \eqref{2.4} test functions $\phivect\in \mathcal{D}(\R\times\Omega)^3$ one obtains \eqref{2.5}. By the previous claim then \eqref{2.7} holds. Consequently, integrating by parts in space and time, we can rewrite \eqref{2.4} as \eqref{2.6}, but with test functions $\psi$ belonging to a different class. Indeed they are of the form $\psi=-\phivect\cdot\boldsymbol{\nu}$, where $\phivect\in C^1_c(\R\times\R^3)^3$ is such that $\phivect\cdot\boldsymbol{\nu}=0$ on $\R\times\Gamma_0$.

 Trivially, such $\psi$ restrict to $\psi\in C_c^1(\R\times\Gamma_1)$, and we actually claim that \eqref{2.6} holds true {\em for all $\psi\in C_c^1(\R\times\Gamma_1)$}, proving the direct implication in the assertion. Our claim follows since, for any $\psi\in C^1_c(\R\times\Gamma_1)$, its trivial extension $\widetilde{\psi}$ on the whole of $\R\times\Gamma$ belongs, as $\overline{\Gamma_0}\cap\overline{\Gamma_1}=\emptyset$, to $C^1_c(\R\times\Gamma)$. Then, defining $\phivect\in C^1_c(\R\times\Gamma)^3$ by $\phivect=-\widetilde{\psi}\boldsymbol{\nu}$, we have $\psi=-\phivect\cdot\boldsymbol{\nu}$ on $\R\times\Gamma_1$ and  $\phivect\cdot\boldsymbol{\nu}=0$ on $\R\times\Gamma_0$. By using the compactness of $\Gamma$, \cite[Definition 1.2.1.1, Chapter 1, p. 5]{grisvard}, local equations, cut--off arguments and partitions of the unity,
 it is then straightforward to extend $\phivect$ to $\phivect\in C^1_c(\R\times\R^3)^3$,  so proving our claim.

 To prove the reverse implication in the statement, we suppose that $(\rvect,v)$ satisfies \eqref{2.5} and \eqref{2.6}. We remark that, by density, \eqref{2.6} holds for all
 $\psi\in C_c^1(\R\times\Gamma_1)$. By using our claim, we can rewrite \eqref{2.5} as \eqref{2.7}. We multiply it by a test function $\phivect$ as in Definition~\ref{Definition 2.1}--iii), and we integrate in space and time the resulting equation. In this way we get
 $$\int_{\R\times\Omega} \rho_0\rvect_t\cdot\phivect_t-B\Div\rvect\Div\phivect=B\int_{\R\times\Gamma_1}\Div\rvect \psi$$
 which, when combined  with \eqref{2.6},  gives \eqref{2.4}, concluding the proof.
 \end{proof}
 Proposition~\ref{Proposition 2.1}  allows to point some trivial relations among the three types of solutions of $(\Lagr)$ introduced in Definition~\ref{Definition 2.1}.
\begin{lem}\label{Lemma 2.2}
Let $(\rvect,v)\in X^1_\Lagr$ be a solution of $(\Lagr)$ according to Definition~\ref{Definition 2.1}. Then strong $\Rightarrow$ generalized $\Rightarrow$ weak and, if $(\rvect,v)\in X^2_\Lagr$, weak $\Rightarrow$ strong.
\end{lem}
\begin{proof} Strong solutions are trivially also generalized and, by Proposition~\ref{Proposition 2.1}, also weak. Since \eqref{2.1}--\eqref{2.4} are stable with respect to the convergence in $X^1_\Lagr$, we also get that generalized $\Rightarrow$ weak.

To prove the  final conclusion let $(\rvect,v)\in X^2_\Lagr$ be a weak solution. By  Proposition~\ref{Proposition 2.1} it  satisfies  all equations in $(\Lagr)$ in a distributional sense. Being $\rvect$ and $v$ regular enough, we can integrate by parts and use \eqref{2.03}. In this way we show that all equations in $(\Lagr)$ also hold a.e.. Hence $(\rvect,v)$ is a strong solution.
\end{proof}
\begin{rem}\label{Remark 2.1}
Hence solutions $(\rvect,v)\in X^2_\Lagr$ are equivalently strong, generalized or weak.
When $(\rvect,v)\in X^1_\Lagr\setminus X^2_\Lagr$, Lemma~\ref{Lemma 2.2}  only says that generalized $\Rightarrow$ weak.

However, as we shall see in the sequel, also the reverse implication weak $\Rightarrow$ generalized holds true. Hence
\emph{all types of solutions in Definition~\ref{Definition 2.1} actually coincide}, strong solutions being defined only in the restricted regularity class $X^2_\Lagr$.
Consequently, in the sequel, we are going to deal only with weak solutions.
\end{rem}
The following result shows that problems $(\Lagr_0^n)$, $n=0,1,2,3,4$, are just restrictions of  problem $(\Lagr_0)$ to data $(\rvect_0,v_0)$ in $H^1_\Lagrn n$.
\begin{prop}\label{Proposition 2.2}
Let $(\rvect,v)\in X^1_\Lagr$ be a weak solution of $(\Lagr_0)$. Then
\renewcommand{\labelenumi}{{\roman{enumi})}}
\begin{enumerate}
\item $\curl \rvect(t)=\curl \rvect_0$ for all $t\in\R$;
\item $\rvect(t)\cdot\nuvect=\rvect_0\cdot\nuvect$ on $\Gamma_0$ for all $t\in\R$;
\item $\rvect(t)\cdot\nuvect+v(t)=\rvect_0\cdot\nuvect+v_0$ on $\Gamma_1$ for all $t\in\R$;
\item $L(\rvect(t),v(t))=L(\rvect_0,v_0)$ for all $t\in\R$.
\end{enumerate}
Consequently,   for any $n=0,1,2,3,4$,
$$(\rvect,v)\quad\text{is a weak solution of $(\Lagr_0^n)$} \Longleftrightarrow  (\rvect,v)\in X^1_\Lagrn n \Longleftrightarrow (\rvect_0,v_0)\in H^1_\Lagrn n.$$
\end{prop}
\begin{proof}To prove i) we point out that, by \eqref{2.1}, for any $\varphi\in \mathcal{D}(\R)$ and $w\in \mathcal{D}(\Omega)$, we have
\begin{equation}\label{2.9}
\int_{\R\times\Omega} \varphi'(t) \curl \rvect(t,x) w(x) dx\,dt=0.
\end{equation}
By density \eqref{2.9} hold for all $w\in L^2(\Omega)$. Then, switching to Bochner's integrals, we get
$$\int_{-\infty}^\infty \varphi' \curl\rvect =0\quad\text{in $L^2(\Omega)^3$, for all $\varphi\in \mathcal{D}(\R)$.}$$
Consequently, we have $\curl \rvect\in W^{1,1}_\loc(\R; L^2(\Omega)^3)$ and $(\curl \rvect)'=0$. Hence
$\curl \rvect\in C^1(\R;L^2(\Omega)^3)$ and i) holds.

By  essentially using the same arguments, replacing $L^2(\Omega)^3$ with $L^2(\Gamma_0)$,
we prove that \eqref{2.2} yields ii). Moreover, replacing it with $L^2(\Gamma_1)$, we prove that \eqref{2.3} implies iii).

To prove iv) we respectively integrate ii) and iii) on $\Gamma_0$ and  $\Gamma_1$. We thus get
\begin{equation}\label{2.9BISS}
\int_{\Gamma_1}\rvect(t)\cdot\nuvect+v(t)+\int_{\Gamma_0}\rvect(t)\cdot\nuvect
=\int_{\Gamma_1}\rvect_0\cdot\nuvect+v_0+\int_{\Gamma_0}\rvect_0\cdot\nuvect\quad\text{for all $t\in\R$.}
\end{equation}
 By \eqref{1.5}, \eqref{2.9BISS} is exactly iv).
 The final conclusion trivially follows by Definition~\ref{Definition 2.1}, i)--iv) above and \eqref{1.6}.
\end{proof}
A natural question arising concerns the  uniqueness of weak solutions of $(\Lagr_0)$.
\begin{prop}\label{Proposition 2.3}
Weak solutions of $(\Lagr_0)$, and hence also of  $(\Lagr_0^n)$, for $n=0,1,2,3,4$, are unique.
\end{prop}
\begin{proof} Uniqueness for weak solutions of $(\Lagr_0^0)$ is known, see \cite[Theorem~4.16]{tremodelli}.

Now let $(\xvect,w)$ and $(\yvect,u)$ be weak solutions of $(\Lagr_0)$, and set $(\rvect,v)=(\xvect,w)-(\yvect,u)$.
By linearity $(\rvect,v)$ is still a weak solution of $(\Lagr_0)$, corresponding to  initial data $\rvect_0=\rvect_1=0$ and $v_0=v_1=0$. Since $(0,0)\in H^1_\Lagrn 0$, by Proposition~\ref{Proposition 2.2} one gets that  $(\rvect,v)$ is a weak solution of $(\Lagr_0^0)$.  By the already recalled result it vanishes, concluding he proof.
\end{proof}
\subsection{Equilibria of problems $\boldsymbol{(\Lagr)}$ and $\boldsymbol{(\Lagrn n)}$. }\label{Section 2.4}
We start by making precise the notions of weak and strong equilibria, as well as stationary solutions, which were quickly introduced in Section~\ref{Section 1.3}.
\begin{definition}\label{Definition 2.2} By a \emph{weak (strong) stationary solution} of problem $(\Lagr)$ we mean a weak (strong) solution $(\rvect,v)$ of it which does not depends on time, i.e. there is $(\svect,z)$ in $H^1_\Lagr$ (in $H^2_\Lagr$) such that $\rvect(t)=\svect$ and $v(t)=z$ for all $t\in\R$.
We call such an $(\svect,z)$ a \emph{weak (strong) equilibrium} of $(\Lagr)$.

Moreover, if $(\svect,z)\in H^1_\Lagrn n$, we  call it a \emph{weak equilibrium} of $(\Lagrn n)$, while if
$(\svect,z)\in H^2_\Lagrn n$, we  call it a \emph{strong equilibrium} of $(\Lagrn n)$.
\end{definition}

We then make precise what we mean by a solution of an equation like \eqref{1.15}.
\begin{definition}\label{Definition 2.3}Let $p_0\in\C$. We say that $z\in H^1(\Gamma_1)$
\renewcommand{\labelenumi}{{\roman{enumi})}}
\begin{enumerate}
\item   is a \emph{strong solution} of the elliptic equation
\begin{equation}\label{2.13}
  -\Div_\Gamma(\sigma \nabla_\Gamma z^*)+\kappa z^*+p=0\qquad\text{on $\Gamma_1$,}
\end{equation}
provided $z\in H^2(\Gamma_1)$ and it satisfies \eqref{2.13} a.e. on $\Gamma_1$;
\item is a \emph{weak solution} of \eqref{2.13} provided
\begin{equation}\label{2.14}
  \int_{\Gamma_1}\sigma (\nabla_\Gamma z,\nabla_\Gamma \psi)_\Gamma +\int_{\Gamma_1}\kappa z\overline{\psi}+p_0\int_{\Gamma_1}\overline{\psi}=0\quad\text{for all $\psi\in H^1(\Gamma_1)$.}
\end{equation}
\end{enumerate}
\end{definition}
By combining \cite[Theorem~5.0.1 and Lemma~6.1.5]{mugnvit} and  some trivial arguments, one immediately gets the following result.
\begin{lem}\label{Lemma 2.3}For any $p_0\in\C$, $z\in H^1(\Gamma_1)$ is a strong solution of \eqref{2.14} if and only if it is a weak solution of it. Moreover,
\renewcommand{\labelenumi}{{\roman{enumi})}}
\begin{enumerate}
\item when $\kappa\equiv 0$ the equation \eqref{2.13} has weak solutions if and only if $p_0=0$.
In this case the weak solutions of \eqref{2.13} constitute the space $\C_{\Gamma_1}$ of the complex constant functions on $\Gamma_1$;
\item when $\kappa\not\equiv 0$ for any $p_0\in\C$ the equation \eqref{2.13} has the unique weak solution $\rho_0z^*/B$, where $z^*$ is the unique weak solution of \eqref{1.17} and it is real--valued.
\end{enumerate}
\end{lem}
The following result gives a preliminary characterization of strong and weak equilibria of $(\Lagr)$.
\begin{lem}\label{Lemma 2.4} Let $(\svect,z)\in H^1_\Lagr$. Then
\renewcommand{\labelenumi}{{\roman{enumi})}}
\begin{enumerate}
\item $(\svect,z)$ is a weak equilibrium of $(\Lagr)$ if and only if $z\in H^2(\Gamma_1)$ and there is $p_0\in\C$ such that
\begin{equation}\label{2.15}
\begin{cases}
  -B\Div \svect=p_0\quad &\text{in $\Omega$,}\\
   -\DivGamma(\sigma\nabla_\Gamma z)+\kappa z+p_0=0\quad &\text{on $\Gamma_1$;}
  \end{cases}
\end{equation}
\item $(\svect,z)$ is a strong equilibrium of $(\Lagr)$ if and only if $(\svect,z)\in H^2_\Lagr$ and there is $p_0\in\C$ such that \eqref{2.15} holds.
\end{enumerate}
\end{lem}
\begin{proof}At first, we set $(\rvect,v)=i(\svect,z)$, see \eqref{1.15}. We shall keep this notation all along the proof.
 To prove i), we first take $z\in H^2(\Gamma_1)$ and $p_0\in\C$ such that \eqref{2.15} holds, claiming that
$(\svect,z)$ is a weak equilibrium of $(\Lagr)$.  Trivially, we have $\Div\rvect\in L^1_\loc(\R;H^1(\Omega))$ and equations \eqref{2.1}--\eqref{2.3},  \eqref{2.5}--\eqref{2.6} hold. Hence, by Proposition~\ref{Proposition 2.1}, $(\rvect,v)$ is a stationary solution of $(\Lagr)$, so that $(\svect,z)$ is a weak equilibrium of $(\Lagr)$.

Conversely, we suppose that $(\svect,z)$ is a weak equilibrium of $(\Lagr)$, so that $(\rvect,v)$ is a weak stationary solution of it. By \eqref{2.4}, we get $\int_{\R\times\Omega} \Div \rvect \Div \phivect=0$ for all $\phivect\in\mathcal{D}(\R\times\Omega)^3$. Consequently, taking test functions $\phivect$ in the separate form $\phivect(t,x)=\varphi(t)\wvect(x)$,  with $\varphi\in\mathcal{D}(\R)$ and $\wvect\in \mathcal{D}(\Omega)^3$, we get
$$\int_{-\infty}^\infty \varphi(t)\left(\int_\Omega \Div \svect(x)\Div \wvect(x)\,dx\right)\,dt=0\quad\text{for all $\varphi\in\mathcal{D}(\R)$.}$$
Consequently, we have $\int_\Omega \Div \svect\Div \wvect=0$. Integrating it by parts we obtain $\Div \svect\in H^1(\Omega)$ and $\nabla\Div \svect=0$. Being $\Omega$ connected, then there is $p_0\in\C$ such that $-B\Div\svect=p_0$. By Proposition~\ref{Proposition 2.1}, we then get that $v$ satisfies \eqref{2.6}, that is
$$\int_{\R\times\Gamma_1}\sigma(\nabla_\Gamma v,\nabla_\Gamma \psi)_\Gamma+
 (\kappa v+\rho_0)\overline{\psi}=0\quad\text{for all $\psi\in C^2_c(\R\times\Gamma_1)$.}$$
Then, taking  $\psi$ in the separate form $\psi(t,y)=\psi_0(t)\psi_1(y)$, where $\psi_0\in C^2_c(\R)$ and
$\psi_1\in C^2(\Gamma_1)$, we get
$$\int_{-\infty}^\infty \overline{\psi_0}(t)\left(\int_{\Gamma_1}\sigma(y)(\nabla_\Gamma z(y),\nabla_\Gamma \psi_1(y))_\Gamma+
 (\kappa z(y)+\rho_0)\overline{\psi_1}(y)\,dy\right)\,dt=0.$$
Since $\psi_0$ is arbitrary we  deduce
$$\int_{\Gamma_1}\sigma(\nabla_\Gamma z,\nabla_\Gamma \psi_1)_\Gamma+
 (\kappa z+\rho_0)\overline{\psi_1}=0\quad\text{for all $\psi_1\in C^2(\Gamma_1)$.}$$
Since $C^2(\Gamma_1$ is dense in $H^1(\Gamma_1)$ we then get \eqref{2.14}, that is to say $z$ is a weak solution of \eqref{2.13}. Hence, by Lemma~\ref{Lemma 2.3}, $z\in H^2(\Gamma_1)$ and \eqref{2.15} holds.

The proof of part i) is then concluded.  To prove part ii) we simply point out that $(\rvect,v)\in X^2_\Lagr$ if and only if $(\svect,z)\in H^2_\Lagr$.
 \end{proof}
\subsection{An auxiliary result.}\label{Section 2.5} In the sequel we shall use the following preliminary result, which
is a variant of \cite[Lemma~4.2]{tremodelli}.
\begin{lem}\label{Lemma 2.5}
For all $w\in L^2(\Omega)$, $z_0\in H^{1/2}(\Gamma_0)$,  and $z_1\in H^{1/2}(\Gamma_1)$, satisfying the compatibility condition
\begin{equation}\label{2.16}
 \int_\Omega w+\int_{\Gamma_0}z_0+\int_{\Gamma_1}z_1=0,
\end{equation}
the problem
\begin{equation}\label{2.17}
  \begin{cases}
\Div \svect=w,\quad &\text{in $\Omega$,}\\
\curl\svect=0\quad &\text{in $\Omega$,}\\
\svect\cdot\boldsymbol{\nu}=-z_0\quad &\text{on $\Gamma_0$,}\\
\svect\cdot\boldsymbol{\nu}=-z_1\quad &\text{on $\Gamma_1$,}
  \end{cases}
\end{equation}
has a unique solution $\svect\in H^1_{\curl 0}(\Omega)$. Moreover, when $w\in H^1(\Omega)$,  $z_0\in H^{3/2}(\Gamma_0)$, and
$z_1\in H^{3/2}(\Gamma_1)$, one has $\svect\in H^2_{\curl 0}(\Omega)$.
\end{lem}
\begin{proof}
By \cite[Lemma 4.1]{tremodelli}, solving problem \eqref{2.17} in the space $H^1_{\curl 0}(\Omega)$ is equivalent to solving the inhomogeneous Neumann problem
\begin{equation}\label{2.18}
  \begin{cases}
  -\Delta \varphi= w,\quad &\text{in $\Omega$,}\\
  \quad \partial_{\boldsymbol{\nu}}\varphi=-z_0\quad &\text{on $\Gamma_0$,}\\
\quad \partial_{\boldsymbol{\nu}}\varphi=-z_1\quad &\text{on $\Gamma_1$,}
  \end{cases}
\end{equation}
with $\varphi\in H^2(\Omega)/\C$, and hence take $\svect=-\nabla\varphi$. Using the decomposition \eqref{2.01} and standard elliptic theory, see for example \cite[Chapter I, Theorem 1.10, p. 15]{GiraultRaviart}, we get that problem  \eqref{2.18}
has a unique solution $\varphi\in H^2(\Omega)/\C$. Moreover, when  $w\in H^1(\Omega)$,  $z_0\in H^{3/2}(\Gamma_0)$, and  $z_1\in H^{3/2}(\Gamma_1)$, we have $\varphi\in H^3(\Omega)/\C$, completing the proof.
\end{proof}
\section{Structural decomposition}\label{Section 3}
This section is devoted to prove Theorem~\ref{Theorem 1.2} and the structural decomposition formula \eqref{1.26} in Theorem~\ref{Theorem 1.3}. To avoid boring repetition, we shall adopt the following convention in this Section. When a Hilbert space $H$ is the direct sum of two among its closed subspaces, i.e. $H=V\oplus W$, we will always denote by $\Pi_{H;V}$ and  $\Pi_{H;W}$ the associated projectors. We notice that, by the Closed Graph Theorem, one has $\Pi_{H;V}\in\mathcal{L}(H;V)$ and $\Pi_{H;W}\in\mathcal{L}(H;W)$, see for example \cite[Chapter III, Theorem~5.20]{Kato}.

Since problems $(\Lagrn 1)$, $(\Lagrn 2)$, and $(\Lagrn 3)$ are the closest ones to problem $(\Lagrn 0)$, we shall preliminarily deal with them, including $(\Lagrn 0)$ for the sake of completeness.
\subsection{Problems $(\Lagrn n)$, $n=0,1,2,3$.}\label{Section 3.1}
At first we characterize weak and strong equilibria.
\begin{prop}\label{Proposition 3.1}
The weak and strong equilibria of problems $(\Lagrn n)$, $n=0,1,2,3$, are given by formulas \eqref{1.22}--\eqref{1.23}.

Moreover, for $k=1,2$, we have
\begin{equation}\label{3.3}
\V^k_3= \V^k_1\oplus \V^k_2,\quad\text{and, consequently,}\quad  \E^k_3= \E^k_1\oplus \E^k_2.
\end{equation}
Furthermore, for $k=1,2$ and $\svect\in \V^k_3$, we have
\begin{equation}\label{3.4}
 \Pi_{\E^k_3;\E^k_1}(\svect,0)=(\avect,0),\quad\text{and}\quad
 \Pi_{\E^k_3;\E^k_2}(\svect,0)=(\bvect,0),
\end{equation}
where $\avect, \bvect\in H^k(\Omega)^3$ respectively are the unique solutions of the  problems
$$
\boxed{(A)
\begin{cases}
\Div \avect=0&\text{in $\Omega$,}\\
\curl \avect=\curl\svect&\text{in $\Omega$,}\\
\avect\cdot\nuvect=0&\text{on $\Gamma$,}
\end{cases}
}
\qquad
\boxed{(B)
\begin{cases}
\Div \bvect=0&\text{in $\Omega$,}\\
\curl \bvect=0&\text{in $\Omega$,}\\
\bvect\cdot\nuvect=\svect\cdot\nuvect&\text{on $\Gamma$.}
\end{cases}
}
$$
\end{prop}
\begin{proof} We fix $k=1,2$. At first we claim that $\E^k_3=\V^k_3\times\{0\}$. To prove our claim we notice that the elements of
$\E^k_3$ are the couples $(\svect,z)\in H^k_\Lagr$ such that $z\in H^2(\Gamma_1)$ and there is $p_0\in\C$ for which the following system holds:
\begin{equation}\label{3.5}
\begin{cases}
  -B\Div \svect=p_0\quad &\text{in $\Omega$,}\\
   -\DivGamma(\sigma\nabla_\Gamma z)+\kappa z+p_0=0\quad &\text{on $\Gamma_1$,}\\
   \int_\Omega \Div\svect+\int_{\Gamma_1}z=0. &
  \end{cases}
\end{equation}
At first, we are going to prove that, when $p_0\not=0$, the system \eqref{3.5} has no solutions.

Indeed, supposing by contradiction that $p_0\not=0$ and $(\svect,z)$ solves it, by normalizing we can suppose that $p_0=1$. Hence, taking real parts, we can also suppose $z$ to be real--valued. Hence, taking $\psi=z$ in \eqref{2.14}, we get
\begin{equation}\label{3.6}
\int_{\Gamma_1}\sigma |\nabla_\Gamma z|^2_\Gamma +\int_{\Gamma_1}\kappa |z|^2+\int_{\Gamma_1}z=0.
\end{equation}
On the other hand, since $p_0=1$, by \eqref{3.5} we have $\int_{\Gamma_1}z=-\int_\Omega\Div\svect=|\Omega|/B$. Plugging it into \eqref{3.6}, since $\kappa,\sigma\ge 0$, we  get that $|\Omega|/B\le 0$, which is the required contradiction.

Consequently, we can rewrite \eqref{3.5} as
\begin{equation}\label{3.7}
\begin{cases}
  -B\Div \svect=0\quad &\text{in $\Omega$,}\\
   -\DivGamma(\sigma\nabla_\Gamma z)+\kappa z=0\quad &\text{on $\Gamma_1$,}\\
   \int_{\Gamma_1}z=0. &
  \end{cases}
\end{equation}
We now apply Lemma~\ref{Lemma 2.3}. Consequently, when $\kappa\equiv 0$, we get $z\in\C_{\Gamma_1}$. By \eqref{3.7}$_3$ we thus obtain $z=0$. When $\kappa\not\equiv 0$ we directly get $z=0$. Hence, in both cases, \eqref{3.7} reduces to $(\svect,z)\in \V^k_3\times\{0\}$, proving our claim.

Our second claim is that $\E^k_2=\V^k_2\times\{0\}$. We point out that, by \eqref{1.7}, we have
$H^k_{\Lagrn 2}\subseteq H^k_{\Lagrn 3}$. Hence, by \eqref{2.11} and \eqref{1.23}, we have
$\E^k_2=\E^k_3\cap H^k_\Lagrn 2=\E^k_3\cap H^1_\Lagrn 2$. Consequently, by \eqref{1.6}, we get
$\E^k_2=\{(\svect,0)\in \E^k_3: \curl \svect=0\}=\V^k_2\times\{0\}$, proving our second claim.

Since $H^k_{\Lagrn 1}\subseteq H^k_{\Lagrn 3}$, the same arguments used above yield that
$\E^k_1=\E^k_3\cap  H^1_\Lagrn 1$. Hence, by \eqref{1.6}, we get
$\E^k_1=\{(\svect,0)\in \E^k_3: \svect\cdot\nuvect=0\quad\text{on $\Gamma$}\}=\V^k_1\times\{0\}$.

To complete the proof of formulas \eqref{1.22}--\eqref{1.23}, we then have to prove that $\E^k_0=\{0\}$.
To achieve this goal we notice that, by \eqref{1.6}, we have $H^1_\Lagrn 0=H^1_\Lagrn 1\cap H^1_\Lagrn 2$. Hence, by \eqref{2.11} and our two first claims, we get $\E^k_0=\E^k_1\cap \E^k_2=(\V_1^k\cap\V_2^k)\times\{0\}$. Consequently, by \eqref{1.23}, we obtain
$$\E^k_0=\{(\svect,0)\in H^k_\Lagr: \Div\svect=0,\quad \curl\svect=0,\quad\text{and}\quad \svect\cdot\nuvect=0\quad\text{on $\Gamma$}\}.$$
By applying Lemma~\ref{Lemma 2.5} we then get $\E^k_0=\{0\}$, completing the proof of \eqref{1.22}--\eqref{1.23}.

We now turn to the decomposition \eqref{3.3}. At first, we remark that, since we just proved that $\E^k_1\cap \E^k_2=\{0\}$, proving \eqref{3.3} reduces to proving that $\E^k_3=\E^k_1+\E^k_2$. To achieve this result, we take any $(\svect,z)=(\svect,0)\in \E^k_3$ and consider problem $(B)$ in the statement.

Since $(\svect,0)\in H^1_\Lagrn 3$, by \eqref{1.6} we have $\int_\Gamma \svect\cdot\nuvect=0$. Hence, by Lemma~\ref{Lemma 2.5}, problem $(B)$ has a unique solution $\bvect\in H^k(\Omega)^3$. By \eqref{1.22}--\eqref{1.23}, we then get $(\bvect,0)\in \E^k_2$. We now set $\avect=\svect-\bvect\in H^k(\Omega)^3$, which trivially solves problem $(A)$ in the statement. Moreover,
by \eqref{1.22}--\eqref{1.23}, one has $(\avect,0)\in \E^k_1$. We thus proved that the decomposition \eqref{3.3} hold and that the associated projectors are given by  \eqref{3.4}.

To complete the proof we point out that, as a consequence of Lemma~\ref{Lemma 2.5}, beside those of problem $(B)$, also solutions of problem $(A)$ are unique. Indeed the quoted result can be applied to the difference between any couple of solutions of $(A)$.
\end{proof}
\begin{rem}\label{Remark 3.1}
In Remark~\ref{Remark 1.4} we noticed that, for each $n=1,2,3$, problem $(\Lagrn n)$  possesses weak equilibria which are not strong.
The proof of this assertion was postponed up to now.

We first claim that $\V^2_1\subsetneq\V^1_1$, so $\E^2_1\subsetneq\E^1_1$, the inclusions being trivial. To prove this fact, let us fix an open subset $\Omega'$ of $\Omega$, compactly contained in it, and any $\uvect\in H^2(\Omega')^3\setminus H^3(\Omega')^3$. We fix a cut--off function $\psi\in \mathcal{D}(\Omega)$ such that $\psi\equiv 1$ in $\Omega'$. Clearly one has $\svect:=\curl (\psi\uvect)\in \V^1_1$. Our claim will then follows, provided $\svect\not\in \V^2_1$, or equivalently $\svect\not\in H^2(\Omega)^3$. Supposing by contradiction $\svect\in H^2(\Omega)^3$, by \cite[Chapter  IX,~Corollary~7, p. 236]{dautraylionsvol3}, one would get  $\psi\uvect\in H^3(\Omega)^3$. Consequently, one would obtain $\uvect\in H^3(\Omega')^3$, so reaching the required contradiction.

We now turn to proving that $\V^2_2\subsetneq\V^1_2$, so $\E^2_2\subsetneq\E^1_2$, also in this case the inclusions being trivial.
Since, as recalled in Remark~\ref{Remark 1.3}, for $k=1,2$ we have
$\V^k_2=\{\nabla\phi: \phi\in H^{k+1}(\Omega)\,\,\text{and}\,\, \Delta\phi=0\}$, the assertion reduces to proving that there are harmonic functions in $H^{k+1}(\Omega)\setminus H^k(\Omega)$, for $k=1,2$. However, this fact is a well--known consequence of standard elliptic theory.

We finally point out that, by \eqref{3.3}, $\E^2_1\subsetneq\E^1_1$ (or $\E^2_2\subsetneq\E^1_2$) yields $\E^2_3\not=\E^1_3$.
\end{rem}
The following result establishes, for $k=1,2$ and $n=1,2,3$, the key decomposition of the space $H^k_{\Lagrn n}$ which allows   to analyse problem $(\Lagrn n)$.
\begin{prop}\label{Proposition 3.2}The following decomposition holds true:
\begin{equation}\label{3.8}
  H^k_\Lagrn n =H^k_\Lagrn 0\oplus\E^k_n,\quad\text{for $k=1,2$ and $n=1,2,3$.}
\end{equation}

Moreover, for all $(\svect,z)\in H^k_\Lagrn n$, we have
\begin{equation}\label{3.9}
 \Pi_{H^k_\Lagrn n ;H^k_\Lagrn 0}(\svect,z)=(\cvect,z),\quad\text{and}\quad
\Pi_{H^k_\Lagrn n;\E^k_n}(\svect,z)=
\begin{cases}
(\avect,0)\quad &\text{when $n=1$,}\\
(\dvect,0)\quad &\text{when $n=2,3$,}
\end{cases}
\end{equation}
where $\avect, \cvect, \dvect\in H^k(\Omega)^3$ respectively are the unique solutions of problem $(A)$ in Proposition~\ref{Proposition 3.1} and of the  following problems:
$$
\boxed{(C)
\begin{cases}
\Div \cvect=\Div \svect&\text{in $\Omega$,}\\
\curl \cvect=0&\text{in $\Omega$,}\\
\cvect\cdot\nuvect=0&\text{on $\Gamma_0$,}\\
\cvect\cdot\nuvect=-z&\text{on $\Gamma_1$,}
\end{cases}
}
\qquad
\boxed{(D)
\begin{cases}
\Div \dvect=0&\text{in $\Omega$,}\\
\curl \dvect=\curl\svect&\text{in $\Omega$,}\\
\dvect\cdot\nuvect=\svect\cdot\nuvect&\text{on $\Gamma_0$,}\\
\dvect\cdot\nuvect=\svect\cdot\nuvect+z&\text{on $\Gamma_1$.}
\end{cases}
}
$$
\end{prop}
\begin{proof}We fix $k=1,2$. By \eqref{2.11} and Proposition~\ref{Proposition 3.1}, for $n=1,2,3$, we have
$H^k_\Lagrn 0\cap\E^k_n=\E^k_0=\{0\}.$ Hence proving \eqref{3.8} reduces to proving that
$H^k_\Lagrn n =H^k_\Lagrn 0+\E^k_n$.

At first, we consider the case $n=1$. Consequently, let us take $(\svect,z)\in H^k_\Lagrn 1$, that is
$(\svect,z)\in H^k_\Lagr$ such that $\svect\cdot\nuvect=0$ on $\Gamma_0$ and $\svect\cdot\nuvect=-z$ on $\Gamma_1$.
By the Divergence Theorem we thus have $\int_\Omega \Div\svect+\int_{\Gamma_1}z=0$. Hence, by Lemma~\ref{Lemma 2.5}, problem $(C)$ has a unique solution $\cvect\in H^k(\Omega)^3$. By \eqref{1.6}, we have $(\cvect,z)\in  H^k_{\Lagrn 0}$. We then set
$\avect:=\svect-\cvect\in H^k(\Omega)^3$, which trivially is the (unique, by Proposition~\ref{Proposition 3.1}) solution of problem $(A)$. Since $(\avect,0)\in\E^k_1$, we proved \eqref{3.8}--\eqref{3.9} when $n=1$.

We now turn to the case $n=2$, and then we take  $(\svect,z)\in H^k_\Lagrn 2$. By \eqref{1.6}, we  have $(\svect,z)\in H^k_\Lagr$, with $\curl \svect=0$ and $\int_\Omega \Div \svect+\int_{\Gamma_1}z=0$. As in the case $n=1$, problem $(C)$ has a unique solution $\cvect\in H^k(\Omega)^3$, and $(\cvect,z)\in H^k_{\Lagrn 0}$. However, in this case, setting $\dvect=\svect-\cvect\in H^k(\Omega)^3$, we get that $\dvect$ solves problem $(D)$. Furthermore, since $\curl \svect=0$, one has $(\dvect,0)\in\E^k_2$.  By applying Lemma~\ref{Lemma 2.5} to the difference between any couple of solutions of problem $(D)$, we get they are unique.

When $n=3$ we simply repeat the arguments used when $n=2$. The only difference is that, in this case, we get $(\dvect,0)\in\E^k_3$.
\end{proof}
\begin{rem}\label{Remark 3.2}
One immediately sees that the projectors $\Pi_{H^k_\Lagrn 1 ;H^k_\Lagrn 0}$  and $\Pi_{H^k_\Lagrn 2;H^k_\Lagrn 0}$ are  restrictions of
$\Pi_{H^k_\Lagrn 3 ;H^k_\Lagrn 0}$. With a little effort one also sees that also the projectors
$\Pi_{H^k_\Lagrn 1 ;\E^k_1}$  and $\Pi_{H^k_\Lagrn 2;\E^k_2}$ are restrictions of
$\Pi_{H^k_\Lagrn 3 ;\E^k_3}$.
\end{rem}
\subsection{Problems $(\Lagr)$ and $(\Lagrn n)$, $n=0,1,2,3,4$.}\label{Section 3.2}
Beside the equilibria already found when analysing problems $(\Lagrn n)$ for $n=1,2,3$, problems
$(\Lagrn 4)$ and $(\Lagr)$ possess the additional one--dimensional space of strong equilibria  $\E^\bullet$, which was introduced, without proofs, in \S~\ref{Section 1.3}.

We are now going to explain how this space is found, in this way proving the assertions made  to setup
the strong equilibrium $(\svect^\bullet ,z^\bullet)$ in \eqref{1.19}.

We point out that, by Lemma~\ref{Lemma 2.4}, weak equilibria of $(\Lagrn 4)$ are exactly couples $(\svect,z)\in H^1_\Lagr$, with $z\in H^2(\Gamma_1)$, verifying the system
  \begin{equation}\label{3.10}
\begin{cases}
  -B\Div \svect=p_0\quad &\text{in $\Omega$,}\\
  \curl \svect=0\quad &\text{in $\Omega$,}\\
   -\DivGamma(\sigma\nabla_\Gamma z)+\kappa z+p_0=0\quad &\text{on $\Gamma_1$,}
  \end{cases}
\end{equation}
whereas the weak equilibria of $(\Lagr)$ have to verify the less demanding system \eqref{2.15}.
We then apply Lemma~\ref{Lemma 2.3} to \eqref{3.10}$_3$,  distinguishing between the cases $\kappa\equiv 0$ and
$\kappa\not\equiv 0$.
\renewcommand{\labelenumi}{{\roman{enumi})}}
\begin{enumerate}
\item When $\kappa\equiv 0$, equation \eqref{3.10}$_3$ has solutions only when $p_0=0$. In this case they are exactly  the elements of $\C_{\Gamma_1}$. Hence, choosing $\svect=0$ in \eqref{3.10}, we get the trivial strong equilibrium
    $(0,\mathbbm{1}_{\Gamma_1})$ of $(\Lagrn 4)$ and $(\Lagr)$.
\item  When $\kappa\not\equiv 0$, for any $p_0\in\C$,  equation \eqref{3.10}$_3$   has a unique solution $z\in H^2(\Gamma_1)$. Taking $p_0=0$, one recovers the space $\E^1_2$, when dealing with $(\Lagrn 4)$, or the space $\E^1_3$, when dealing with $(\Lagr)$.

    Hence we consider the case $p_0=B\not=0$, in which equation \eqref{3.10}$_3$ reduces to equation \eqref{1.17}. It  as the unique solution $z^*\in H^2(\Gamma_1)$ already introduced in \S~\ref{Section 1.3}.

    We then have to look for solutions $\svect$ of \eqref{3.10}$_1$--\eqref{3.10}$_2$.  We also ask that
    $\svect\cdot\nu=0$ on $\Gamma_0$ and  $\svect\cdot\nu=-|\Omega|/\mathcal{H}^2(\Gamma_1)$ on $\Gamma_1$. These, somehow arbitrary, boundary conditions were chooses to fulfill in the simplest possible way the compatibility condition \eqref{2.16} in Lemma~\ref{Lemma 2.5}. The resulting system is nothing but the problem \eqref{1.18}, which then has the unique solution $\svect^*\in H^2(\Omega)^3$ already introduced in \S~\ref{Section 1.3}.

    From what preceeds we got the strong equilibrium $(\svect^*,z^*)\in H^2_\Lagrn 4\subseteq H^2_\Lagr$.
    \end{enumerate}
We can then set $(\svect^\bullet ,z^\bullet)$ by  \eqref{1.19} and $\E^\bullet$ by \eqref{1.20}. Moreover, as
we announced in \S~\ref{Section 1.3}, by the following result  $\E^\bullet$ is not included in $H^1_\Lagrn n$ when $n=0,1,2,3$.
\begin{lem}\label{Lemma 3.3}
 We have $L(\svect^\bullet,z^\bullet)\not=0$, so that
 \begin{equation}\label{3.11}
   \E^\bullet\cap H^1_\Lagrn n=\E^\bullet\cap H^2_\Lagrn n=\{0\}\qquad\text{for $n=0,1,2,3$.}
 \end{equation}
 Hence we can set $\ell\in (H^1_\Lagr)'$ by formula \eqref{1.21}, and
  $\ell(\svect^\bullet,z^\bullet)=1$.

  Moreover, for all $(\svect,z)\in H^1_\Lagr$, we have
  \begin{equation}\label{3.1}
\ell(\svect,z)=
\begin{cases}
\dfrac{\int_\Gamma \svect\cdot\nuvect+\int_{\Gamma_1}z}{\mathcal{H}^2(\Gamma_1)}\quad &\text{when $\kappa\equiv 0$,}\\
\dfrac{\int_\Gamma \svect\cdot\nuvect+\int_{\Gamma_1}z}{-|\Omega|+\int_{\Gamma_1}z^*}\quad &\text{when $\kappa\not\equiv 0$.}
\end{cases}
\end{equation}
\end{lem}
\begin{proof}When $\kappa\equiv 0$, by \eqref{1.5} and \eqref{1.19}, we have
\begin{equation}\label{3.2}
L(\svect^\bullet,z^\bullet)=L(0,\mathbbm{1}_{\Gamma_1})=\mathcal{H}^2(\Gamma_1)>0.
\end{equation}
When
$\kappa\not\equiv 0$, by \eqref{1.5}, \eqref{1.18} and \eqref{1.19}, we get
\begin{equation}\label{3.12}
L(\svect^\bullet,z^\bullet)= L(\svect^*,z^*)=\int_\Omega \Div \svect^*+\int_{\Gamma_1}z^*=-|\Omega|+\int_{\Gamma_1}z^*.
\end{equation}
Taking $\psi=z^*$ in \eqref{2.14}, since $z^*$ is real--valued, we get
$$
 \int_{\Gamma_1}\sigma |\nabla_\Gamma z^*|^2_\Gamma+\int_{\Gamma_1}\kappa |z^*|^2+B\int_{\Gamma_1}z^*=0.
$$
Consequently, since $\sigma,\kappa\ge 0$, we have $\int_{\Gamma_1}z^*\le 0$. Then \eqref{3.12} yields $L(\svect^*,z^*)<0$.

We thus proved that $L(\svect^\bullet,z^\bullet)\not=0$. We can then set $\ell\in (H^1_\Lagr)'$ by formula \eqref{1.21}, so
  $\ell(\svect^\bullet,z^\bullet)=1$. Moreover, by \eqref{3.1} and \eqref{3.12}, also \eqref{3.1} follows.
Finally, by \eqref{1.6}, we get \eqref{3.11}.
\end{proof}
Our next result establishes a preliminary decomposition of the spaces $H^k_\Lagrn 4$ and $H^k_\Lagr$, for $k=1,2$, making explicit the projectors associated with it.
\begin{lem}\label{Lemma 3.4}
The following decompositions hold true:
\begin{equation}\label{3.13}
H^k_\Lagr=H^k_\Lagrn 3\oplus \E^\bullet,\quad\text{and}\quad H^k_\Lagrn 4=H^k_\Lagrn 2\oplus \E^\bullet,\quad\text{for $k=1,2$.}
\end{equation}
Moreover, for all  $(\svect,z)\in H^k_\Lagr$, we have
\begin{equation}\label{3.14}
 \Pi_{H^k_\Lagr;H^k_\Lagrn 3}(\svect,z)=(\svect,z)-\ell(\svect,z)(\svect^\bullet,z^\bullet),\,\text{and}\,\,\,\,
 \Pi_{H^k_\Lagr;\E^\bullet}(\svect,z)=\ell(\svect,z)(\svect^\bullet,z^\bullet).
\end{equation}
Furthermore, the projectors $\Pi_{H^k_\Lagrn 4;H^k_\Lagrn 2}$ and $\Pi_{H^k_\Lagrn 4;\E^\bullet}$  are the  restrictions to $H^k_\Lagrn 4$ of those in \eqref{3.14}.
\end{lem}
\begin{proof}We fix $k=1,2$.
Since, by \eqref{1.6}, $H^1_\Lagrn 3=\text{Ker }L$ and $H^2_\Lagrn 3=\text{Ker }L\cap H^2_\Lagr$, the operators in \eqref{3.14} are well--defined. They show that $H^k_\Lagr=H^k_\Lagrn 3+ \E^\bullet$, so by Lemma~\ref{Lemma 3.3} we get the decomposition $H^k_\Lagr=H^k_\Lagrn 3\oplus \E^\bullet$.

Since $\curl\svect^\bullet=0$, for all $(\svect,z)\in H^k_\Lagrn 4$ we have $\Pi_{H^k_\Lagr;H^k_\Lagrn 3}(\svect,z)\in H^k_\Lagrn 2$. Hence we also get the decomposition  $H^k_\Lagrn 4=H^k_\Lagrn 2\oplus \E^\bullet$, the projectors associated to it being the restrictions of the ones given in \eqref{3.14}.
\end{proof}

Lemma~\ref{Lemma 3.4} allows to characterize the equilibria of problems $(\Lagrn 4)$ and   $(\Lagr)$.
\begin{prop}\label{Proposition 3.5}
The weak and strong equilibria of problems $(\Lagrn 4)$ and   $(\Lagr)$ are given by \eqref{1.24}. Moreover, we also have the decompositions
\begin{equation}\label{3.15}
 \E^k=\E^k_1\oplus \E^k_2\oplus \E^\bullet=\E^k_1\oplus \E^k_4,\qquad\text{for $k=1,2$.}
\end{equation}
\end{prop}
\begin{proof}We fix $k=1,2$. At first, to prove  formula \eqref{1.24}, we claim that $\E^k_4=\E^k_2\oplus \E^\bullet$.

By \eqref{2.12} and Lemma~\ref{Lemma 3.3} we have $\E^k_2\cap\E^\bullet=H^k_\Lagrn 2\cap \E^\bullet=\{0\}$. Moreover, for all $(\svect,z)\in \E^k_4$, we have
$$\Pi_{H^k_\Lagrn 4; H^k_\Lagrn 2}(\svect,z)=(\svect,z)-\ell(\svect,z)(\svect^\bullet, z^\bullet)\in \E^k\cap H^k_\Lagrn 2=\E^k_2.$$
Since $(\svect^\bullet, z^\bullet)\in\E^\bullet$, we then get $\E^k_4=\E^k_2\oplus \E^\bullet$, proving our claim.

Exactly the same arguments show that $\E^k=\E^k_3\oplus \E^\bullet$, so proving  \eqref{1.24}.
 By combining these decompositions with \eqref{3.3}, we get \eqref{3.15}, concluding the proof.
\end{proof}
We can now give the
\begin{proof}[\bf Proof of Theorem~\ref{Theorem 1.2}]
One simply combines Propositions~\ref{Proposition 3.1} and \ref{Proposition 3.5}.
\end{proof}
Lemma~\ref{Lemma 3.4} and Proposition~\ref{Proposition 3.2} allow to prove the following intermediate decompositions of
$H^k_\Lagrn 4$ and $H^k_\Lagr$, for $k=1,2$, making explicit the projectors associated with them.
\begin{lem}\label{Lemma 3.6}
The following decompositions hold true:
\begin{equation}\label{3.16}
H^k_\Lagr=H^k_\Lagrn 0\oplus \E^k_3\oplus\E^\bullet,\quad\text{and}\quad H^k_\Lagrn 4=H^k_\Lagrn 0\oplus \E^k_2\oplus\E^\bullet,\quad\text{for $k=1,2$.}
\end{equation}
Moreover, the projectors associated with the first splitting in \eqref{3.16} are the operators
$\Pi_{H^k_\Lagr;\E^\bullet}$, still given by \eqref{3.14}, and the operators
$\Pi_{H^k_\Lagr;H^k_\Lagrn 0}$, $\Pi_{H^k_\Lagr;\E^k_3}$. They are given,  for all $(\svect,z)\in H^k_\Lagr$, by
\begin{equation}\label{3.17}
 \Pi_{H^k_\Lagr;H^k_\Lagrn 0}(\svect,z)=(\gvect, z-\ell(\svect,z)z^\bullet),\quad\text{and}\quad
 \Pi_{H^k_\Lagr;\E^k_3}(\svect,z)(\hvect,0),
\end{equation}
where $\gvect, \hvect\in H^k(\Omega)^3$ respectively are the unique solutions of the  problem $(G)$ in Theorem~\ref{Theorem 1.3} and of the following one:
$$
(H)
\begin{cases}
\Div \hvect=0&\text{in $\Omega$,}\\
\curl \hvect=\curl\svect&\text{in $\Omega$,}\\
\hvect\cdot\nuvect=\svect\cdot\nuvect&\text{on $\Gamma_0$,}\\
\hvect\cdot\nuvect=\svect\cdot\nuvect+z-\ell(\svect,z)(\svect^\bullet\cdot\nuvect+z^\bullet)\negquad\!\!&\text{on $\Gamma_1$.}
\end{cases}
$$
Moreover, the projectors $\Pi_{H^k_\Lagrn 4;\E^\bullet}$,
$\Pi_{H^k_\Lagrn 4;H^k_\Lagrn 0}$, and  $\Pi_{H^k_\Lagrn 4;\E^k_2}$ respectively are the restrictions of the ones associated with the first splitting in  \eqref{3.16}.
\end{lem}
\begin{proof}We fix $k=1,2$. The splitting $H^k_\Lagr=H^k_\Lagrn 0\oplus \E^k_3\oplus\E^\bullet$
simply follows by combining those given in Lemma~\ref{Lemma 3.4} and Proposition~\ref{Proposition 3.2}, that is $H^k_\Lagr=H^k_\Lagrn 3\oplus\E^\bullet$ and $H^k_\Lagrn 3=H^k_\Lagrn 0\oplus \E^k_3$. Moreover, $\Pi_{H^k_\Lagr;\E^\bullet}$ is still given by \eqref{3.14}, while the projectors on $\E^k_3$ and $H^k_\Lagrn 0$ are obtained by using the trivial formulas
$$\Pi_{H^k_\Lagr;\E^k_3}=\Pi_{H^k_\Lagrn 3;\E^k_3}\cdot \Pi_{H^k_\Lagr; H^k_\Lagrn 3},\quad\text{and}\quad \Pi_{H^k_\Lagr;H^k_\Lagrn 0}=\Pi_{H^k_\Lagrn 3;H^k_\Lagrn 0}\cdot \Pi_{H^k_\Lagr; H^k_\Lagrn 3}.$$
Consequently, by \eqref{3.14}, for all $(\svect,z)\in H^k_\Lagr$, we have
\begin{equation}\label{3.18}
 \begin{aligned}
 &\Pi_{H^k_\Lagr;\E^k_3}(\svect,z)=\Pi_{H^k_\Lagrn 3;\E^k_3}
 (\svect-\ell(\svect,z)\svect^\bullet, z-\ell(\svect,z)z^\bullet),\quad\text{and}\\
 &\Pi_{H^k_\Lagr;H^k_\Lagrn 0}(\svect,z)=\Pi_{H^k_\Lagr;H^k_\Lagrn 0}
 (\svect-\ell(\svect,z)\svect^\bullet, z-\ell(\svect,z)z^\bullet).
\end{aligned}
\end{equation}
To use formulas \eqref{3.9} in \eqref{3.18} we point out that, by \eqref{1.18}--\eqref{1.19}, we have $\curl\svect^\bullet=0$ in $\Omega$ and $\svect\cdot\nuvect=0$ on $\Gamma_0$, in both cases. In this way, starting from \eqref{3.18} and performing a trivial calculation, we get \eqref{3.17}.

We then apply Lemma~\ref{Lemma 2.5} to the difference of two arbitrary solutions
of each one between problems $(G)$ and $(H)$. In this way we get that solutions of these problems are unique.

As far as the splitting $H^k_\Lagrn 4=H^k_\Lagrn 0\oplus \E^k_2\oplus\E^\bullet$ is concerned, it simply follows by combining the decompositions $H^k_\Lagrn 4=H^k_\Lagrn 2\oplus\E^\bullet$ and
$H^k_\Lagrn 2=H^k_\Lagrn 0\oplus\E^k_2$, already given in Lemma~\ref{Lemma 3.4} and Proposition~\ref{Proposition 3.2}. Moreover, $\Pi_{H^k_\Lagrn 4; \E^\bullet}$ is still given by \eqref{3.14}.

Finally, since we trivially have
$$\Pi_{H^k_\Lagrn 4;\E^k_2}=\Pi_{H^k_\Lagrn 2;\E^k_2}\cdot \Pi_{H^k_\Lagrn 4; H^k_\Lagrn 2},\quad\text{and}\quad \Pi_{H^k_\Lagrn 4;H^k_\Lagrn 0}=\Pi_{H^k_\Lagrn 2;H^k_\Lagrn 0}\cdot \Pi_{H^k_\Lagrn 4; H^k_\Lagrn 2},$$
by Lemma~\ref{Lemma 3.4} and Proposition~\ref{Proposition 3.2}, these operators are obtained by restriction.
\end{proof}
The following result establishes the key  formula \eqref{1.26}.
\begin{prop}\label{Proposition 3.7}
The decompositions \eqref{1.26} hold true. Furthermore, the projectors associated to it are given as  stated in Theorem~\ref{Theorem 1.3}.
\end{prop}
\begin{proof}We fix $k=1,2$. The splitting $H^k_\Lagr=H^k_\Lagrn 0\oplus\E^k$ simply follows by combining the splitting $H^k_\Lagr=H^k_\Lagrn 0\oplus\E_3^k\oplus\E^\bullet$, given in
Lemma~\ref{Lemma 3.6}, with the one $\E^k=\E^k_3\oplus\E^\bullet$, given in Proposition~\ref{Proposition 3.5}.
By combining the decompositions $H^k_\Lagrn 4=H^k_\Lagrn 0\oplus\E_2^k\oplus\E^\bullet$ and $\E^k_4=\E^k_2\oplus\E^\bullet$, given in the same result, we get $H^k_\Lagrn 4=H^k_\Lagrn 0\oplus\E^k_4$. Moreover, by Proposition~\ref{Proposition 3.2}, we have
$H^k_\Lagrn n=H^k_\Lagrn 0\oplus\E^k_n$ for $n=0,1,2,3$, so completing the proof of \eqref{1.26}.

To recognize that the projectors associated with the splitting $H^k_\Lagr=H^k_\Lagrn 0\oplus\E^k$ are given by formula \eqref{1.27}, we remark that $\Pi_{H^k_\Lagr;H^k_\Lagrn 0}$ was already made explicit in \eqref{3.17}, while the form of $\Pi_{H^k_\Lagr;\E^k}$ can be simply obtained by using the trivial formula $\Pi_{H^k_\Lagr;\E^k}=\Pi_{H^k_\Lagr;\E^k_3}+\Pi_{H^k_\Lagr;\E^\bullet}$, which is a consequence of the decomposition $\E^k=\E^k_3\oplus\E^\bullet$. By using \eqref{3.17} and \eqref{3.14}, we thus get
\begin{equation}\label{3.21}
\Pi_{H^k_\Lagr;\E^k}(\svect,z)=(\hvect,0)+\ell(\svect,z)(\svect^\bullet,z^\bullet)=(\fvect,\ell(\svect,z)z^\bullet)\quad\text{for all $(\svect,z)\in H^k_\Lagr$,}
\end{equation}
where $\hvect\in H^k(\Omega)^3$ solves problem $(H)$ in Lemma~\ref{Lemma 3.6} and $\fvect:=\hvect+\ell(\svect,z)\svect^\bullet$. Since $\curl\svect^\bullet=0$ and $\svect^\bullet\cdot\nuvect=0$ on $\Gamma_0$, we  $\fvect$ solves the problem
$$\begin{cases}
\Div \fvect=\Div \hvect+\ell(\svect,z)\Div \svect^\bullet=\ell(\svect,z)\Div \svect&\text{in $\Omega$,}\\
\curl\fvect=\curl \hvect=\curl\svect&\text{in $\Omega$,}\\
\fvect\cdot\nuvect=\hvect\cdot\nuvect=\svect\cdot\nuvect&\text{on $\Gamma_0$,}\\
\fvect\cdot\nuvect=\hvect\cdot\nuvect+\ell(\svect,z)\svect^\bullet\cdot\nuvect=\svect\cdot\nuvect+z-\ell(\svect,z)z^\bullet\negquad\!\!&\text{on $\Gamma_1$,}
\end{cases}
$$
which is nothing but problem $(F)$ in Theorem~\ref{Theorem 1.3}. We thus got \eqref{1.27}.
Furthermore, by applying Lemma~\ref{Lemma 2.5} once more, we get that solutions of problem $(F)$ are unique, like the ones  of
problem $(G)$, this fact having already been  proved in Lemma~\ref{Lemma 3.6}.

As far as the projectors associated with the splitting $H^k_\Lagrn n=H^k_\Lagrn 0\oplus\E^k_n$, for $n=01,2,3,4$, are concerned, we notice that since $H^k_\Lagrn n\subset H^k_\Lagr$, the projector $\Pi_{H^k_\Lagrn n;H^k_\Lagrn 0}$ is restriction of $\Pi_{H^k_\Lagr;H^k_\Lagrn 0}$. Consequently, since $\Pi_{H^k_\Lagr;\E^k_n}=I-\Pi_{H^k_\Lagr;H^k_\Lagrn 0}$,  also
$\Pi_{H^k_\Lagr;\E^k_n}$ is restriction of $\Pi_{H^k_\Lagr;\E^k}$.
\end{proof}
Together with the main decomposition results established above, in the next section we shall also use a couple of their consequences.
The first one, which could be of independent interest, gives the ``atomic'' decompositions of the spaces $H^k_\Lagr$ and $H^k_\Lagrn n$, for $k=1,2$ and $n=1,2,3,4$, as well as the expressions of the  projectors associated with them.

\begin{cor}\label{Corollary 3.8} The following decompositions hold true for $k=1,2$:
\begin{equation}\label{3.22}
\begin{gathered}
H^k_\Lagrn 1=H^k_\Lagrn 0\oplus \E^k_1,\quad H^k_\Lagrn 2=H^k_\Lagrn 0\oplus \E^k_2, \quad H^k_\Lagrn 3=H^k_\Lagrn 0\oplus \E^k_1\oplus \E^k_2,\\
H^k_\Lagrn 4=H^k_\Lagrn 0\oplus \E^k_2\oplus \E^\bullet,\qquad
\text{and}\quad H^k_\Lagr=H^k_\Lagrn 0\oplus \E^k_1\oplus\E^k_2\oplus \E^\bullet.
\end{gathered}
\end{equation}
The four projectors associated with the splitting of $H^k_\Lagr$ in \eqref{3.22} are the operators
$\Pi_{H^k_\Lagr;H^k_\Lagrn 0}$, given by \eqref{1.27}, $\Pi_{H^k_\Lagr;\E^\bullet}$, given by \eqref{3.14}, together with the operators
$\Pi_{H^k_\Lagr;\E^k_1}$ and $\Pi_{H^k_\Lagr;\E^k_2}$. They are given, for all $(\svect,z)\in H^k_\Lagr$, by
\begin{equation}\label{3.23}
\Pi_{H^k_\Lagr;\E^k_1}(\svect,z)=(\avect,0),\qquad   \Pi_{H^k_\Lagr;\E^k_2}(\svect,z)=(\kvect,0),
\end{equation}
where $\avect,\kvect\in H^k(\Omega)^3$ respectively are the unique solutions of problem $(A)$ in Proposition~\ref{Proposition 3.1} and of the following problem
$$
(K)
\begin{cases}
\Div \kvect=0&\text{in $\Omega$,}\\
\curl \kvect=0&\text{in $\Omega$,}\\
\kvect\cdot\nuvect=\svect\cdot\nuvect&\text{on $\Gamma_0$,}\\
\kvect\cdot\nuvect=\svect\cdot\nuvect+z-\ell(\svect,z)(\svect^\bullet\cdot\nuvect+z^\bullet)\negquad\!\!&\text{on $\Gamma_1$.}
\end{cases}
$$
Moreover, any projector associated with one of the decompositions in \eqref{3.22}, onto one of the spaces  $H^k_\Lagrn 0$, $\E^k_1$,
$\E^k_2$, and $\E^\bullet$, is obtained as the restriction of the projector from  $H^k_\Lagr$ onto the same space.
\end{cor}
\begin{proof}At first, we fix $k=1,2$ and we point out that, by Proposition~\ref{Proposition 3.2}, we have $H^k_\Lagrn n=H^k_\Lagrn 0\oplus \E^k_n$ for $n=1,2,3$. Since, by Proposition~\ref{Proposition 3.1}, we have $\E^k_3=\E^k_1\oplus \E^k_2$, all decompositions in the first line of \eqref{3.22} have already been proved. Moreover, by \eqref{3.16}, we also obtain $H^k_\Lagrn 4=H^k_\Lagrn 0\oplus \E^k_2\oplus \E^\bullet$.  Furthermore, by combining \eqref{3.16} with \eqref{3.3}, we also get the decomposition
$H^k_\Lagr=H^k_\Lagrn 0\oplus \E^k_3\oplus \E^\bullet=H^k_\Lagrn 0\oplus \E^k_1\oplus\E^k_2\oplus \E^\bullet$, proving \eqref{3.22}.

Turning to the forms of the projectors listed in the statement, the first two of them are already known, while the last two are simply obtained from Lemma~\ref{Lemma 3.6} and Proposition~\ref{Proposition 3.1}, since for $n=1,2$ we have
$\Pi_{H^k_\Lagr;\E^k_n}=\Pi_{\E^k_3;\E^k_n}\cdot \Pi_{H^k_\Lagr; \E^k_3}$. Hence, by \eqref{3.17}, for all $n=1,2$ and $(\svect,z)\in H^k_\Lagr$, we have $\Pi_{H^k_\Lagr;\E^k_n}(\svect,z)=\Pi_{\E^k_3;\E^k_n}(\hvect,0)$, where $\hvect\in H^k(\Omega)^3$ solves problem $(H)$ in Lemma~\ref{Lemma 3.6}. Then, by \eqref{3.4}, we get $\Pi_{H^k_\Lagr;\E^k_1}(\svect,z)=(\avect,0)$, where
$\avect\in H^k(\Omega)^3$ solves problem $(A)$. In the same way we also get $\Pi_{H^k_\Lagr;\E^k_2}(\svect,z)=(\kvect,0)$, where
$\kvect\in H^k(\Omega)^3$ is the unique solution of  problem $(K)$.
The final assertion in the statement is trivial, since for $n=1,2,3,4$ we have $H^k_\Lagrn n\subset H^k_\Lagr$.
\end{proof}
\begin{rem}\label{Remark 3.3}Corollary~\ref{Corollary 3.8} shows that, for any fixed $k=1,2$, the spaces $H^k_\Lagr$ and
$H^k_\Lagrn n$, $n=0,1,2,3,4$, can all be  obtained as direct sums of the main base space $H^k_\Lagrn 0$ with some (if any) of the spaces $\E^k_1$, $\E^k_2$ and $\E^\bullet$, according to the following Table. These four spaces can then be considered as the ``atoms'' of all decompositions.
\begin{center}
\begin{tabular}{|c|c|c|c|c|}
\hline
  &$H^k_\Lagrn 0$  &$\E^k_1$&$\E^k_2$& $\E^\bullet$\\
\hline
$H^k_\Lagrn 0$&\Checkmark&\XSolidBrush &\XSolidBrush & \XSolidBrush\\
\hline
$H^k_\Lagrn 1$&\Checkmark&\Checkmark&\XSolidBrush & \XSolidBrush\\
\hline
$H^k_\Lagrn 2$&\Checkmark&\XSolidBrush &\Checkmark & \XSolidBrush\\
\hline
$H^k_\Lagrn 3$&\Checkmark&\Checkmark &\Checkmark & \XSolidBrush\\
\hline
$H^k_\Lagrn 4$&\Checkmark&\XSolidBrush &\Checkmark & \Checkmark\\
\hline
$H^k_\Lagr$&\Checkmark&\XSolidBrush &\XSolidBrush & \XSolidBrush\\
\hline
\end{tabular}
\end{center}
As it was shown in Theorem~\ref{Theorem 1.3}, the space $H^k_\Lagrn 0$ is the only atomic space where non--stationary solutions
of problem $(\Lagr)$ live. For this reason, giving a look to all possible direct sums of this type can be of some interest.

As it is clear from the Table above, there are only $2=2^3-6$ possible spaces of this type which are not the configuration space of some Lagrangian model. The first one is the space $H^k_\Lagrn 0\oplus \E^\bullet$, while the second one is the space
$H^k_\Lagrn 0\oplus \E^k_1\oplus\E^\bullet$. It looks clear that the latter is the ``rotational'' version of the former. The space $H^k_\Lagrn 0\oplus \E^\bullet$ does not look to be the configuration space of a reasonable evolution problem. On the other hand, elements of $\mathcal{S}^k_\Lagrn 0\oplus\mathcal{S}^\bullet$,  which importance was pointed out in Remark~\ref{Remark 1.6}, live in it.
\end{rem}
The second consequence of our main decomposition result will be useful in the next section. While Corollary~\ref{Corollary 3.8} shows how configuration spaces of the Lagrangian models are obtained by adding something to $H^k_\Lagrn 0$, the following result shows how adding something to them one gets $H^k_\Lagr$.
\begin{cor}\label{Corollary 3.9}
The following decompositions hold true for $k=1,2$:
\begin{equation}\label{3.24}
H^k_\Lagr =H^k_\Lagrn 4\oplus \E^k_1= H^k_\Lagrn 3\oplus \E^\bullet=H^k_\Lagrn 2\oplus \E^k_1\oplus \E^\bullet=H^k_\Lagrn 1\oplus \E^k_2\oplus \E^\bullet.
\end{equation}
Moreover, $\Pi_{H^2_\Lagr;H^2_\Lagrn n}$  is the restriction of  $\Pi_{H^1_\Lagr;H^1_\Lagrn n}$ to $H^2_\Lagr$ for $n=1,2,3,4$.
\end{cor}
\begin{proof}
The decompositions \eqref{3.24} immediately follow from \eqref{3.22}, as the Table in Remark~\ref{Remark 3.3} makes evident. Moreover, since $H^2_\Lagr\subset H^1_\Lagr$ and $H^2_\Lagrn n\subset H^1_\Lagrn n$ for $n=1,2,3,4$, the last assertion is trivial.
\end{proof}
\section{Proofs of our main results}\label{Section 4}
This section is devoted to prove our main results stated in Section~\ref{Section 1}, but for Theorem~\ref{Theorem 1.2} which was already proved in the previous Section.
\begin{proof}[\bf Proof of Theorem~\ref{Theorem 1.1}] We fix data $U_0=(\rvect_0,v_0,\rvect_1,v_1)\in \mathcal{H}^1_\Lagr=H^1_\Lagr\times H^0$. By Proposition~\ref{Proposition 3.7}, we can uniquely decompose $(\rvect_0,v_0)$ as the sum of
$(\rvect_0^0,v_0^0):=\Pi_{H^1_\Lagr;H^1_\Lagrn 0}(\rvect_0,v_0)\in H^1_\Lagrn 0$ and $(\rvect_0^1,v_0^1):=\Pi_{H^1_\Lagr;\E^1 }(\rvect_0,v_0)\in \E^1$. Hence, setting $U_0^0:=(\rvect_0^0,v_0^0,\rvect_1,v_1)$ and $U_0^1:=(\rvect_0^1,v_0^1,0,0)$, we clearly have $U_0=U_0^0+U_0^1$.

Trivially $i(\rvect_0^1,v_0^1)\in \SSS^1$ is a weak stationary solution of problem $\Lagr_0$ corresponding to initial data $U_0^1$.
We then consider problem $\Lagr_0$ in correspondence to initial data $U_0^0$.

Since $U_0^0\in H^1_\Lagrn 0\times H^0$, by Proposition~\ref{Proposition 2.2},  for these data problem $(\Lagrn 0)$ is equivalent to problem  $(\Lagr_0^0)$. Now, as far as problem
$(\Lagr_0^0)$ is concerned, Theorem~\ref{Theorem 1.1} reduces to \cite[Theorem~1.3]{tremodelli}. We can then apply it, also recalling (see \cite[p. 36]{tremodelli}) that for problem $(\Lagr_0^0)$ the concepts of generalized and weak solutions are equivalent.
In this way we obtain that $(\Lagr_0^0)$ has a unique generalized and weak solution $(\rvect^0,v^0)\in X^1_\Lagrn 0$, continuously depending on $U^0_0$. We also get that $(\rvect^0,v^0)\in X^2_\Lagrn 0$ if and only if $U^0_0\in\mathcal{H}^2_\Lagrn 0$, in this case
$(\rvect^0,v^0)$ being a strong solution.

Since $U_0=U_0^0+U_0^1$, the couple $(\rvect,v)\in X^1_\Lagr$ given by \eqref{1.28}
is a generalized and weak solution of problem $(\Lagr_0)$, continuously depending on $U_0$.
It is given by
\begin{equation}\label{4.1}
  \rvect(t)=\rvect^0(t)+\rvect_0^1,\qquad v(t)=v^0(t)+v^1_0,\qquad\text{for all $t\in\R$.}
\end{equation}
 Moreover,
Proposition~\ref{Proposition 2.3} asserts its uniqueness among weak solutions of $(\Lagr_0)$.

We now claim that $(\rvect,v)\in X^2_\Lagr$ if and only if $U_0\in \mathcal{H}^2_\Lagr=H^2_\Lagr\times H^1_\Lagrn 0$. To prove our claim we first suppose that $(\rvect,v)\in X^2_\Lagr$. By \eqref{1.12}, we get $(\rvect_0,v_0)=(\rvect(0),v(0))\in H^2_\Lagr$ and $(\rvect_1,v_1)=(\rvect_t(0), v_t(0))\in H^1_\Lagrn 0$, so $U_0\in \mathcal{H}^2_\Lagr$.

Conversely, let us suppose that $U_0\in \mathcal{H}^2_\Lagr$, so that $(\rvect_0,v_0)\in H^2_\Lagr$ and $(\rvect_1,v_1)\in H^1_\Lagrn 0$. By Proposition~\ref{Proposition 3.7}, the splitting $H^2_\Lagr=H^2_\Lagrn 0\oplus \E^2$ holds true, so $(\rvect_0^0,v_0^0)=\Pi_{H^2_\Lagr; H^2_\Lagrn 0}(\rvect_0,v_0)\in H^2_\Lagrn 0$ and $(\rvect_0^1,v_0^1)\in\E^2$. By applying \cite[Theorem~1.3]{tremodelli} once again, we get $(\rvect^0,v^0)\in X^2_\Lagr$. Then, by \eqref{4.1}, we have $(\rvect,v)\in X^2_\Lagr$, proving our claim.

We now claim that $\mathcal{H}^2_\Lagr=H^2_\Lagr\times H^1_\Lagrn 0$ is dense in  $\mathcal{H}^1_\Lagr=H^1_\Lagr\times H^0$.
As stated in \cite[Theorem~1.3]{tremodelli}, the space $\mathcal{H}^2_\Lagrn 0=H^2_\Lagrn 0\times H^1_\Lagrn 0$ is dense in
$\mathcal{H}^1_\Lagrn 0=H^1_\Lagrn 0\times H^0$. Hence $H^1_\Lagrn 0$ is dense in
$H^0$. Proving our claim then reduces to proving the density of $\mathcal{H}^2_\Lagr$ in $\mathcal{H}^1_\Lagr$.
Now $H^k_\Lagr=H^k(\Omega)^3\times H^k(\Gamma_1)$ for $k=1,2$. Moreover, by standard density results, $H^2(\Omega)^3$ and $H^2(\Gamma_1)$ are  dense, respectively in $H^1(\Omega)^3$ and in $H^1(\Gamma_1)$. Hence our claim follows.

Our next claim concerns the validity of the energy identity \eqref{1.14}. By our two previous claims, we can use a density argument and restrict to strong solutions $(\rvect,v)\in X^2_\Lagr$. We then take any such solution, multiply $(\Lagr_0)$ by $\overline{\rvect_t}$ and integrate by parts in $\Omega$. Then, for all $s,t\in\R$, we get
$$
\int_s^t\left[\int_\Omega\rho_0\rvect_{tt}\overline{\rvect_t}+B\Div\rvect\Div \overline{\rvect_t}-B\int_\Gamma\Div\rvect\,\,\overline{\rvect_t}\cdot\nuvect\right]=0.
$$
First  multiplying $(\Lagr_0)_3$ by $\overline{v_t}$,  then using \eqref{2.3}, for all $s,t\in\R$, we get
$$
\int_s^t\int_{\Gamma_1}\mu v_{tt}\overline{v_t}+\sigma(\nabla_\Gamma v,\nabla_\Gamma v_t)_\Gamma+\delta |v_t|^2+\kappa v \overline{v_t}-B\Div \rvect\,\overline{v_t}=0.
$$
Summing the last two equations and using $(\Lagr_0)_4$, we get
$$
\int_s^t\left[\int_\Omega\rho_0\rvect_{tt}\overline{\rvect_t}+B\Div\rvect\Div \overline{\rvect_t}+\int_{\Gamma_1}\mu v_{tt}\overline{v_t}+\sigma(\nabla_\Gamma v,\nabla_\Gamma v_t)_\Gamma+\delta |v_t|^2+\kappa v \overline{v_t}\right]=0.
$$
By taking the real part we thus obtain
$$
\tfrac 12\int_s^t\frac d{dt}\left(\int_\Omega\rho_0|\rvect_t|^2+B|\Div\rvect|^2+\int_{\Gamma_1}\mu |v_t|^2+\sigma|\nabla_\Gamma v|^2_\Gamma+\kappa |v|^2\right)+\int_s^t\int_{\Gamma_1}\delta |v_t|^2=0,
$$
which is nothing but \eqref{1.14}, so proving our claim.

We thus proved Theorem~\ref{Theorem 1.1}, as far as problem $(\Lagr^0_0)$ is concerned. We then turn to problem
$(\Lagr^n_0)$, for $n=1,2,3,4$.
By Definition~\ref{Definition 2.1} and Proposition~\ref{Proposition 2.2}, the solution $(\rvect,v)\in X^1_\Lagr$ of
$(\Lagr_0)$ is the unique solution of $(\Lagr^n_0)$ if and only if $(\rvect,v)\in X^1_\Lagrn n$, and in turn if and only if $U_0\in H^1_\Lagrn n$. Consequently, problem $(\Lagr^n_0)$ is well--posed  in $H^1_\Lagrn n$, in the Hadamard sense. Moreover, since $X^2_\Lagrn n=X^2_\Lagr\cap C(\R;H^1_\Lagrn n)$, by \eqref{1.10}--\eqref{1.12} the equivalence
$(\rvect,v)\in X^2_\Lagrn n\Longleftrightarrow U_0\in \mathcal{H}^2_\Lagrn n$ holds true.

To complete the proof, we claim that $\mathcal{H}^2_\Lagrn n=H^2_\Lagrn n\times H^1_\Lagrn 0$ is dense in
$\mathcal{H}^1_\Lagrn n=H^1_\Lagrn n\times H^0$.  Since, as already recognized, $H^1_\Lagrn 0$ is dense in $H^0$, our claim reduces to the density of  $H^2_\Lagrn n$  in $H^1_\Lagrn n$. To prove this fact we shall use the (already proved) density of $H^2_\Lagr$ in $H^1_\Lagr$, together with Corollary~\ref{Corollary 3.9}. Indeed, since $H^1_\Lagrn n\subset H^1_\Lagr$, for any $(\svect,z)\in H^1_\Lagrn n$ there is a sequence $((\svect_k,z_k))_k$ in $H^2_\Lagr$ such that
$(\svect_k,z_k)\to (\svect,z)$ in $H^1_\Lagr$. Since $\Pi_{H^1_\Lagr;H^1_\Lagrn n}$ is bounded, we  have
$$(\svect^n_k,z^n_k):=\Pi_{H^1_\Lagr;H^1_\Lagrn n}(\svect_k,z_k)\to \Pi_{H^1_\Lagr;H^1_\Lagrn n}(\svect,z)=(\svect,z)\quad\text{in $H^1_\Lagrn n$, as $k\to\infty$.}$$
Since $(\svect_k,z_k)\in H^2_\Lagr$, by Corollary~\ref{Corollary 3.9} we have $(\svect^n_k,z^n_k)\in H^2_\Lagrn n$, hence proving our claim and completing the proof.
\end{proof}
\begin{rem}\label{Remark 4.1}
As announced in Remark~\ref{Remark 2.1}, Theorem~\ref{Theorem 1.1} also shows the equivalence between weak and generalized solutions of $(\Lagr)$ and $(\Lagrn n)$, $n=0,1,2,3,4$. Indeed, we already pointed out that generalized solutions are also weak. By Theorem~\ref{Theorem 1.1}, any weak solution of $(\Lagr)$ coincides, by uniqueness, with the generalized solution of problem $(\Lagr_0)$ with initial data $U_0:=(\rvect(0),v(0),\rvect_t(0),v_t(0))$.
\end{rem}
\begin{proof}[\bf Proof of Theorem~\ref{Theorem 1.3}] At first, we point out that the decompositions \eqref{1.26} hold true, their associated projectors being given by \eqref{1.27}. Indeed, this result is nothing but Proposition~\ref{Proposition 3.7}. Moreover, the decomposition \eqref{1.28} of the solution coincides with the one given by \eqref{4.1}.

Furthermore, the properties of the projectors having $H^k_\Lagrn n$ as domain, for $k=1,2$ and $n=0,1,2,3,4$, were proved in Proposition~\ref{Proposition 3.7}. Next, by applying Proposition~\ref{Proposition 3.2}, one recognizes that the decomposition \eqref{1.28} continues to hold when respectively replacing  $(\Lagr_0)$, $H^k_\Lagr$, $\mathcal{H}^k_\Lagr$,
and $\E^k$ with $(\Lagr^n_0)$, $H^k_\Lagrn n$, $\mathcal{H}^k_\Lagrn n$, and $\E^k_n$. The decompositions \eqref{1.29} are  consequences of \eqref{1.28}.
\end{proof}
Before proving Theorems~\ref{Theorem 1.4} and \ref{Theorem 1.5}, we recall some notions from \cite{tremodelli}. At first, we introduce the Hilbert spaces $H^k_\Eul:=H^{k-1}(\Omega)\times H^{k-1}_{\curl 0}(\Omega)\times H^k(\Gamma_1)$ for $k=1,2$, the functional $L_\Eul\in (H^1_\Eul)'$ given by
$$L_\Eul(q,\wvect,z):=\int_\Omega q-B\int_{\Gamma_1}z\qquad\text{for all $(q,\wvect,z)\in H^1_\Eul$},$$
together with the further Hilbert spaces $H^k_\Eulc=H^k_\Eul\cap\text{Ker }L_\Eul$, $k=1,2$.

The spaces $X^k_\Eul$ and $X^k_\Eulc$, $k=1,2$, in \eqref{1.3bis} can  be rewritten as follows:
\begin{equation}\label{4.4}
\begin{aligned}
&X_\Eul^1=\{(p,\vvect,v)\in C(\R;H^1_\Eul): v\in C^1(\R;L^2(\Gamma_1))\},\\
&X_\Eul^2=\{(p,\vvect,v)\in C(\R;H^2_\Eul)\cap C^1(\R;H^1_\Eul): v\in C^2(\R;L^2(\Gamma_1))\},\\
&X_\Eulc^1=X_\Eul^1\cap C(\R;H^1_\Eulc),\quad X_\Eulc^2=X_\Eul^2\cap C(\R;H^2_\Eulc).
\end{aligned}
\end{equation}
As it was announced in Section~\ref{Section 1.1}, we are now going to recall the   definitions of solutions of problems $(\Eul)$, $(\Eulc)$, $(\Pot)$, and $(\Potc)$, given in \cite{tremodelli}.  We are also going to correct a missprint in formula $(5.21)$ of the quoted paper (a ``$B$'' was missing).
\begin{definition}\label{Definition 4.1}
We say that
\renewcommand{\labelenumi}{{\roman{enumi})}}
\begin{enumerate}
\item $(p,\vvect,v)\in X^1_\Eul$ is a {\em weak solution} of $(\Eul)$ provided the distributional identities
\begin{align}\label{4.5}
&\int_{\R\times\Omega} p\varphi_t+B\vvect\cdot\nabla\varphi+B\int_{\R\times\Gamma_1}v_t\varphi=0,\\
\label{4.6}
&\int_{\R\times\Omega}  \rho_0\vvect\cdot \phivect_t  +p\Div \phivect + \int_{\R\times\Gamma_1}\mu v_t\psi_t-\sigma(\nabla_\Gamma v,\nabla_\Gamma \overline{\psi})_\Gamma-(\delta v_t
 +\kappa v)\psi=0,
 \end{align}
 hold for all $\varphi\in C^2_c(\R\times\R^3)$ and $\phivect\in C^1_c(\R\times\R^3)^3$ such that $\phivect\cdot\boldsymbol{\nu}=0$ on $\R\times\Gamma_0$, where $\psi=-\phivect\cdot\boldsymbol{\nu}$ on $\R\times\Gamma_1$;
\item $(p,\vvect,v)\in X^2_\Eul$ is a {\em strong solution} of $(\Eul)$ provided $(\Eul)_1$ -- $(\Eul)_2$ hold a.e. in $\R\times\Omega$ and $(\Eul)_4$ -- $(\Eul)_5$ hold a.e. on $\R\times\Gamma_1$;
\item a weak solution $(p,\vvect,v)$ of $(\Eul)$ is a {\em weak solution} of $(\Eulc)$ provided $(p,\vvect,v)\in X^1_\Eulc$;
\item a strong solution $(p,\vvect,v)$ of $(\Eul)$ is a {\em strong solution} of $(\Eulc)$ provided $(p,\vvect,v)\in X^2_\Eulc$;
\item solutions in i)--iv) above are \emph{stationary} if they are constant in time, such a constant being called
an \emph{equilibrium}.
\end{enumerate}
\end{definition}
\begin{rem}\label{Remark 4.2} By \cite[Lemma~5.8]{tremodelli}, strong solutions in Definition~\ref{Definition 4.1}are exactly weak solutions belonging to $X^2_\Eul$.
\end{rem}
\begin{definition}\label{Definition 4.2}
We say that
\renewcommand{\labelenumi}{{\roman{enumi})}}
\begin{enumerate}
\item $(u,v)\in X^1_\Pot$ is a {\em weak solution} of $(\Pot)$ provided the distributional identities
\begin{align}\label{4.7}
&\int_{\R\times\Omega}  -\rho_0 u_t\varphi_t+B\nabla u\nabla \varphi-B\int_{\R\times\Gamma_1}v_t\varphi=0,\\
\label{4.8}
 &\int_{\R\times \Gamma_1}-\mu v_t\psi_t+\sigma(\nabla_\Gamma v,\nabla_\Gamma \overline{\psi})_\Gamma+(\delta v_t
 +\kappa v)\psi-\rho_0 u\psi_t=0,
 \end{align}
hold for all $\varphi\in C^\infty_c(\R\times\R^3)$ and $\psi\in C^2_c(\R\times\Gamma_1)$;
\item $(u,v)\in X^2_\Pot$ is a {\em strong solution} of $(\Pot)$ provided $(\Pot)_1$ holds a.e. in $\R\times\Omega$ and
$(\Pot)_2$--$(\Pot)_3$ hold a.e. on $\R\times\Gamma_1$;
\item a weak solution $(u,v)$ of $(\Pot)$ is a {\em weak solution} of $(\Potc)$ provided $(u,v)\in X^1_\Potc$;
\item a strong solution $(u,v)$ of $(\Pot)$ is a {\em strong solution} of $(\Potc)$ provided $(u,v)\in X^2_\Potc$.
\end{enumerate}
\end{definition}
\begin{rem}\label{Remark 4.3} By \cite[Lemma~3.4]{tremodelli}, strong solutions in Definition~\ref{Definition 4.2} are exactly weak solutions belonging to $X^2_\Pot$.
\end{rem}
Like problem $(\Lagr)$, also problem $(\Eul)$ has a particular strong equilibrium, which takes two different forms depending on the alternative $\kappa\equiv 0$ or $\kappa\not\equiv 0$:
\renewcommand{\labelenumi}{{\roman{enumi})}}
\begin{enumerate}
\item when $\kappa\equiv 0$, it is the triplet $(0,0,\mathbbm{1}_{\Gamma_1})\in H^1_\Eul$, since  $p(t)=0$, $\vvect(t)=0$ and $v(t)=\mathbbm{1}_{\Gamma_1}$, for all $t\in\R$, define a strong stationary solution of $(\Eul)$;
\item when $\kappa\not\equiv 0$, it is the triplet $(B,0,z^*)\in H^1_\Eul$, since $p(t)=B$, $\vvect(t)=0$ and $v(t)=z^*$, for all $t\in\R$, define a strong stationary solution of $(\Eul)$.
\end{enumerate}
To treat  the two cases  above at once we  set the strong equilibrium
of problem $(\Eul)$:
\begin{equation}\label{4.9}
(q^\bullet ,0,z^\bullet)=
\begin{cases}
(0,0,\mathbbm{1}_{\Gamma_1})\qquad&\text{when $\kappa\equiv 0$,}
\\
(B,0,z^*)\qquad&\text{when $\kappa\not\equiv 0$.}
\end{cases}
\end{equation}
 Moreover, we set the corresponding strong stationary solution $(p^\bullet,0,v^\bullet)\in X^2_\Eul$ of $(\Eul)$ by setting $(p^\bullet,0,v^\bullet)(t)=(q^\bullet ,0,z^\bullet)$ for all $t\in\R$.
We also denote by $\E_\Eul^\bullet$ and $\SSS_\Eul^\bullet$ the one--dimensionale spaces of strong equilibria and stationary
solutions which are generated by them, that is
\begin{equation}\label{4.10}
\E_\Eul^\bullet=\C (q^\bullet,0,z^\bullet),\qquad\text{and}\quad \SSS_\Eul^\bullet=\C(p^\bullet,0,v^\bullet).
\end{equation}
To prove Theorem~\ref{Theorem 1.4} we are going to use the following  result.
\begin{lem}\label{Lemma 4.1} We have $L_\Eul(q^\bullet,0,z^\bullet)\not=0$, and consequently, for $k=1,2$, the following decompositions hold true:
\begin{equation}\label{4.11}
H^k_\Eul=H^k_\Eulc\oplus \E^\bullet_\Eul,\quad X^k_\Eul=X^k_\Eulc\oplus \SSS^\bullet_\Eul,\quad\text{and}\quad
\mathcal{S}^k_\Eul=\mathcal{S}^k_\Eulc\oplus\SSS^\bullet_\Eul.
\end{equation}
\end{lem}
\begin{proof}When $\kappa\equiv 0$, by \eqref{4.9}, we have $L_\Eul(q^\bullet,0,z^\bullet)=-B\mathcal{H}^2(\Gamma_1)<0$. When $\kappa\not\equiv 0$,  by \eqref{4.9}, we get $L_\Eul(q^\bullet,0,z^\bullet)=B|\Omega|-B\int_{\Gamma_1}z^*$. By using Lemma~\ref{Lemma 3.3} and formula \eqref{3.12}, we obtain
$L_\Eul(q^\bullet,0,z^\bullet)=BL(\svect^\bullet,z^\bullet)<0$. Hence, in both cases, $L_\Eul(q^\bullet,0,z^\bullet)\not=0$.

Since $H^1_\Eulc=\text{Ker }L_\Eul$ and $H^2_\Eulc=H^2_\Eul\cap\text{Ker }L_\Eul$, we get $H^k_\Eulc\cap \E_\Eul^\bullet=\{0\}$ for $k=1,2$,  i.e. $(q^\bullet,0,z^\bullet)$ is not an equilibrium of problem $(\Eulc)$.
The splitting $H^k_\Eul=H^k_\Eulc\oplus \E^\bullet_\Eul$, for $k=1,2$, then trivially follows. Combining it with the definitions of $X^k_\Eulc$, $k=1,2$, see \eqref{4.4}, one gets the splitting $X^k_\Eul=X^k_\Eulc\oplus \SSS^\bullet_\Eul$ for $k=1,2$. Since $\mathcal{S}^\bullet\subset\mathcal{S}^2_\Eul\subset\mathcal{S}^1_\Eul$, also the splitting $\mathcal{S}^k_\Eul=\mathcal{S}^k_\Eulc\oplus\SSS^\bullet_\Eul$, $k=1,2$, follows.
\end{proof}
\begin{proof}[\bf Proof of Theorem~\ref{Theorem 1.4}]By \eqref{1.3bis}, \eqref{1.9} and \eqref{1.12}, the map $(\rvect,v)\mapsto (p,\vvect, v)$ defined by \eqref{1.31} trivially sets a linear bounded operator $\Psi^1_\LE$ from $\mathcal{S}^1_\Lagr$ to $X^1_\Eul$, which restricts to a bounded operator $\Psi^2_\LE$ from $\mathcal{S}^2_\Lagr$ to $X^2_\Eul$.

Moreover, by \eqref{1.6} and \eqref{1.9}, for all $(\rvect,v)\in \mathcal{S}^1_\Lagrn 3$, we have
$$L_\Eul(p,\vvect,v)=\int_\Omega p-B\int_{\Gamma_1}v=-B\int_{\Omega}\Div\rvect-B\int_{\Gamma_1}v=-BL(\rvect,v)=0.$$
Consequently, by \eqref{1.7} and \eqref{1.9}, for $n=0,1,2,3$, the operator $\Psi^1_\LE$  restricts to a bounded operator
$\Psi^1_\LnEc n$ from $\mathcal{S}^1_\Lagrn n$ to $X^1_\Eulc$ and to a bounded operator
$\Psi^2_\LnEc n$ from $\mathcal{S}^2_\Lagrn n$ to $X^2_\Eulc$.

We now claim that, for any weak solution $(\rvect,v)$ of problem $(\Lagr)$, the triplet $(p,\vvect,v)$ given by \eqref{1.31} is a weak solution of problem $(\Eul)$. To prove our claim, at first we point out that, by \eqref{2.4}, we immediately get
that $(p,\vvect,v)$  satisfies \eqref{4.6}. Moreover,  \eqref{4.5}  reduces to
$$\int_{\R\times\Omega}-\Div\rvect \,\varphi_t+\int_{\R\times\Omega}\rvect_t\cdot\nabla\varphi+\int_{\R\times\Gamma_1}v_t\,\varphi=0$$
for all $\varphi\in C^2_C(\R\times\R^3)$.
Integrating by parts in time the last two addenda in it we get
$$\int_{\R\times\Omega}\Div(\varphi_t\rvect) +\int_{\R\times\Gamma_1}v\,\varphi_t=0,$$
which, by the Divergence Theorem, reads as
$$\int_{\R\times\Gamma_1}(\rvect\cdot\nuvect+v)\varphi_t+\int_{\R\times\Gamma_0}\rvect\cdot\nuvect\,\varphi_t=0.$$
Now, the last identity holds true thanks to \eqref{2.2}--\eqref{2.3}, so proving our claim.

Moreover, when $(\rvect,v)$ is a strong solution, $(p,\vvect,v)$ is a strong solution as well. Furthermore,
when $(\rvect,v)$ solves $(\Lagrn n)$ for $n=0,1,2,3$, $(p,\vvect,v)$ solves $(\Eulc)$. We thus obtained the operators
\begin{equation}\label{4.12}
  \Psi^k_\LE\!\! \in \mathcal{L}(\mathcal{S}^k_\Lagr;\mathcal{S}^k_\Eul),\quad \Psi^k_\LnE 4\!\!\in \mathcal{L}(\mathcal{S}^k_\Lagrn 4;\mathcal{S}^k_\Eul),\quad\text{and}\quad \Psi^k_\LnEc n\in \mathcal{L}(\mathcal{S}^k_\Lagrn n;\mathcal{S}^k_\Eulc),
\end{equation}
for $n=0,1,2,3$, in the statement.

The operators $\Psi^k_\LnEc 0$, $k=1,2$, have been already introduced in \cite[Corollaries~1.9~and~5.13]{tremodelli}, although a different notation was used. Indeed, $\Psi^1_\LnEc 0$ was denoted as $\Psi_\LEc $ in Corollary~1.9 in the quoted paper, while $\Psi^2_\LnEc 0$ corresponds to the operator \begin{footnote}{by a missprint, in formula (5,39) the notation $\dot{\Psi}^2_\EcL$ is used for this operator and its inverse.}\end{footnote} $\dot{\Psi}^2_\LEc$  in Corollary~5.13.
Using these results, we know that, for $k=1,2$, $\Psi^k_\LnEc 0$ is a bijective isomorphism. From now on we fix $k$.

The inverse of $\Psi^k_\LnEc 0$ is the operator
$\Psi^k_\EcLn 0\in \mathcal{L}(\mathcal{S}^k_\Eulc,\mathcal{S}^k_\Lagrn 0)$. It is defined, for all $(p,\vvect,v)\in \mathcal{S}^k_\Eulc$, by $\Psi^k_\EcLn 0(p,\vvect,v)=(\dot{\rvect},v)$, where $\dot{\rvect}(t)\in H^k_{\curl 0}(\Omega)$ is the unique solution of problem \eqref{1.34}, for all $t\in\R$. Clearly, $(\dot{\rvect},v)$ coincides with the unique solution $(\rvect,v)$ of problem $(\Lagrn 0)$ satisfying the equations \eqref{1.31}.

We now turn to the operator $\Psi^k_\LnEc n$, for $n=1,2,3$. Since $X^k_\Lagrn 0\subset X^k_\Lagrn n$, the operator
$\Psi^k_\LnEc n$ restricts to the surjective operator $\Psi^k_\LnEc 0$, hence it is surjective as well. Differently from $\Psi^k_\LnEc 0$, the operator $\Psi^k_\LnEc n$, for $n=1,2,3$, is not injective. Indeed, by \eqref{1.31} and Theorem~\ref{Theorem 1.2}, one easily gets that $\text{Ker }\Psi^k_\LnEc n=\SSS^k_n$.

Consequently, it subordinates the bijective isomorphism $\dot{\Psi}^k_\LnEc n$ given by \eqref{1.32}. The inverse of this isomorphism is the operator $\dot{\Psi}^k_\EcLn n\in \mathcal{L}(\mathcal{S}^k_\Eulc,\dot{\mathcal{S}}^k_\Lagrn n)$, trivially given by \eqref{1.33}. Moreover, since  $\dot{\mathcal{S}}^k_\Lagrn n=\mathcal{S}^n_\Lagrn n /\SSS_n^k$ and $\SSS_n^k=\text{Ker }\Psi^k_\LnEc n$, we also have $\dot{\Psi}^k_\EcLn n(p.\vvect,v)= (\dot{\rvect},v)+\SSS^k_n$.

We now turn to the operators $\Psi^k_\LnE 4$ and $\Psi^k_\LE$. To deal with them, it is useful to point out suitable decompositions of their domains and codomains. By \eqref{1.24}, we get $\SSS^k_4=\SSS^k_2\oplus \SSS^\bullet$ and
$\SSS^k=\SSS^k_3\oplus \SSS^\bullet$. Moreover, by Theorem~\ref{Theorem 1.3}, we have
$\mathcal{S}^k_\Lagrn n=\mathcal{S}^k_\Lagrn 0\oplus\SSS^k_n$, for $n=2,3,4$, and
$\mathcal{S}^k_\Lagr=\mathcal{S}^k_\Lagrn 0\oplus\SSS^k$. Combining these decompositions with those given above, and with \eqref{4.11}, we get
\begin{equation}\label{4.13}
  \mathcal{S}^k_\Lagrn 4=\mathcal{S}^k_\Lagrn 2\oplus\SSS^\bullet,\quad
  \mathcal{S}^k_\Lagr=\mathcal{S}^k_\Lagrn 3\oplus\SSS^\bullet,\quad\text{and}\quad
  \mathcal{S}^k_\Eul=\mathcal{S}^k_\Eulc\oplus\SSS^\bullet_\Eul.
\end{equation}
We shall now give the matrix representations of $\Psi^k_\LnE 4$ and $\Psi^k_\LE$, respectively using the first and third  decompositions and the last two ones in \eqref{4.13}.

Since $(-B\Div \svect^\bullet,0,z^\bullet)=(q^\bullet,0,z^\bullet)$, the operators $\Psi^k_\LnE 4$ and $\Psi^k_\LE$ restrict to a bijective isomorphism $\Psi^\bullet$ from $\SSS^\bullet$ onto $\SSS_\Eul^\bullet$, which are both one--dimensional spaces. Moreover, $\Psi^k_\LnE 4$ and $\Psi^k_\LE$ respectively restrict to $\Psi^k_\LnEc 2$  on $\mathcal{S}^k_\Lagrn 2$, and to $\Psi^k_\LnEc 3$ on $\mathcal{S}^k_\Lagrn 3$. Consequently, the matrix representations described above are
\begin{equation}\label{4.14}
\Psi^k_\LnE 4=
\begin{pmatrix}\Psi^k_\LnEc 2 &  0\\
0 & \Psi^\bullet
\end{pmatrix}
\qquad\text{and}\quad
\Psi^k_\LE =
\begin{pmatrix}\Psi^k_\LnEc 3 &  0\\
0 & \Psi^\bullet
\end{pmatrix}.
\end{equation}
By using \eqref{4.14} and the (already proved) properties of the operators $\Psi^k_\LnEc 2$ and $\Psi^k_\LnEc 3$,
we get that also $\Psi^k_\LnE 4$ and $\Psi^k_\LE$ are surjective, with
$$\text{Ker }\Psi^k_\LnE 4=\text{Ker }\Psi^k_\LnEc 2=\SSS^k_2,\quad\text{and}\quad
\text{Ker }\Psi^k_\LE=\text{Ker }\Psi^k_\LnEc 3=\SSS^k_3.$$

Consequently, these two operators respectively subordinate the bijective isomorphisms
$\dot{\Psi}^k_\LnE 4\in \mathcal{L}(\dot{\mathcal{S}}^k_\Lagrn 4;\mathcal{S}^k_\Eul)$ and
$\dot{\Psi}^k_\LE \in \mathcal{L}(\dot{\mathcal{S}}^k_\Lagr;\mathcal{S}^k_\Eul)$, both given in \eqref{1.32}.
Their inverses are the operators $\dot{\Psi}^k_\ELn 4$ and $\dot{\Psi}^k_\EL$, trivially given in \eqref{1.33}.
\end{proof}
\begin{proof}[\bf Proof of Theorem~\ref{Theorem 1.5}]The proof essentially consists in combining Theorem~\ref{Theorem 1.4} and  \cite[Theorems~1.7~and~5.11]{tremodelli}. These results  were partially recalled in the unnumbered Theorem in Section~\ref{Section 1.1}. While the operators $\Psi^k_\PE$,
$\Psi^k_\PcEc$, $\dot{\Psi}^k_\PE$, $\dot{\Psi}^k_\PcEc$, $k=1,2$, are essentially given by formula \eqref{1.3ter}, we need to recall the expressions  of  $\dot{\Psi}^k_\EP=(\dot{\Psi}^k_\PE)^{-1}\in\mathcal{L}(\mathcal{S}^k_\Eul;\dot{\mathcal{S}}^k_\Pot)$
and $\dot{\Psi}^k_\EcPc=(\dot{\Psi}^k_\PcEc)^{-1}\in\mathcal{L}(\mathcal{S}^k_\Eulc;\dot{\mathcal{S}}^k_\Potc)$. For all $(p,\vvect,v)\in\mathcal{S}^1_\Eul$, we have
\begin{equation}\label{4.15}
 \dot{\Psi}^1_\EP(p,\vvect,v)=(u,v)+\C_{X_\Pot},
\end{equation}
where, correcting a missprint, $u$ is given, up to a space--time constant, by
 \begin{equation}\label{4.16}
  u(t)=u(0)+\tfrac 1{\rho_0}\int_0^t p(\tau)\,d\tau,\quad\text{and}\quad -\nabla u(0)=\vvect(0).
\end{equation}
The operators $\dot{\Psi}^2_\EP$, $\dot{\Psi}^1_\EcPc$, and $\dot{\Psi}^2_\EcPc$, are the appropriate restrictions of
$\dot{\Psi}^1_\EP$.

We then fix $k=1,2$ and set  the operator $\Psi^k_\LP:=\dot{\Psi}^k_\EP\cdot\Psi^k_\LE\in\mathcal{L}(\mathcal{S}^k_\Lagr;\dot{\mathcal{S}}^k_\Pot)$, with
$\dot{\Psi}^k_\EP$ given by \eqref{4.15} and $\Psi^k_\LE$ given by \eqref{1.31}. Hence, for all $(\rvect,v)\in \mathcal{S}^k_\Lagr$, we have
$$\Psi^k_\LP(\rvect,v)=\dot{\Psi}^k_\EP(-B\Div \rvect,\rvect_t,v)=(u,v)+\C_{X_\Pot},$$
where $u$ is given, up a space--time constant, by \eqref{1.37}.  Since $\dot{\Psi}^k_\EP$ is bijective , by using Theorem~\ref{Theorem 1.4}, we get that $\Psi^k_\LP$ is surjective, with $\text{Ker }\Psi^k_\LP=\text{Ker }\Psi^k_\LE=\SSS^k_3$.  Hence $\Psi^k_\LP$ subordinates the bijective isomorphism
$\dot{\Psi}^k_\LP=\dot{\Psi}^k_\EP\cdot\dot{\Psi}^k_\LE\in\mathcal{L}(\dot{\mathcal{S}}^k_\Lagr;\dot{\mathcal{S}}^k_\Pot)$,
given by \eqref{1.39}. Its inverse $\dot{\Psi}^k_\PL$ is trivially given by \eqref{1.40}.

The operator $\Psi^k_\LP$ restricts to $\Psi^k_\LnP 4\in\mathcal{L}(\mathcal{S}^k_\Lagrn 4;\dot{\mathcal{S}}^k_\Pot)$,
which is also given by $\Psi^k_\LnP 4=\dot{\Psi}^k_\EP\cdot \Psi^k_\LnE 4$. Hence, by Theorem~\ref{Theorem 1.4}, since $\dot{\Psi}^k_\EP$ is bijective, the operator $\Psi^k_\LnP 4$ is surjective as $\Psi^k_\LnE 4$, and we have $\text{Ker }\Psi^k_\LnP 4=\text{Ker }\Psi^k_\LnE 4=\SSS^k_2$.

Consequently, it subordinates the bijective isomorphism
$\dot{\Psi}^k_\LnP 4=\dot{\Psi}^k_\EP\cdot\dot{\Psi}^k_\LnE 4\in\mathcal{L}(\dot{\mathcal{S}}^k_\Lagrn 4;\dot{\mathcal{S}}^k_\Pot)$, given by \eqref{1.39}. Its inverse $\dot{\Psi}^k_\PLn 4$ is trivially given by \eqref{1.40}.

Fixing for $n=0,1,2,3$, the operator $\Psi^k_\LP$ further restricts to $\Psi^k_\LnPc n\in\mathcal{L}(\mathcal{S}^k_\Lagrn n;\dot{\mathcal{S}}^k_\Potc)$, which is also given by $\Psi^k_\LnPc n=\dot{\Psi}^k_\EcPc \cdot \Psi^k_\LnEc n$. Hence, by Theorem~\ref{Theorem 1.4}, since $\dot{\Psi}^k_\EcPc$ is bijective, the operator $\Psi^k_\LnPc n$ is surjective as $\Psi^k_\LnE n$, and we have $\text{Ker }\Psi^k_\LnPc n=\text{Ker }\Psi^k_\LnEc n=\SSS^k_n$.

Consequently, it subordinates the bijective isomorphism
$\dot{\Psi}^k_\LnPc n=\dot{\Psi}^k_\EcPc\cdot\dot{\Psi}^k_\LnEc n\in\mathcal{L}(\dot{\mathcal{S}}^k_\Lagrn n;\dot{\mathcal{S}}^k_\Potc)$, given by \eqref{1.39}. Its inverse $\dot{\Psi}^k_\PcLn n$ is trivially given by \eqref{1.40}.

Furthermore, when $n=0$, by Theorem~\ref{Theorem 1.4}, we also get that
\begin{equation}\label{4.17}
\dot{\Psi}^k_\PcLn 0[(u,v)+\C_{X_\Pot}]=(\dot{\rvect},v),
\end{equation}
where $\dot{\rvect}$ is defined by \eqref{1.41}. Since $\text{Ker }\Psi^k_\LnPc n=\SSS^k_n$ and  $\dot{\mathcal{S}}^k_\Lagrn n=\mathcal{S}^n_\Lagrn n /\SSS_n^k$, by \eqref{4.17} we get $\dot{\Psi}^k_\PcLn n[(u,v)+\C_{X_\Pot}]=(\dot{\rvect},v)+\SSS^k_n$,  completing the proof.
\end{proof}
We can finally give the proof of our last main result.
\begin{proof}[\bf Proof of Theorem~\ref{Theorem 1.6}] By using Theorem~\ref{Theorem 1.3},
for any $U_0=(\rvect_0,v_0,\rvect_1,v_1)\in \mathcal{H}^1_\Lagr=H^1_\Lagr\times H^0$, we can decompose the weak solution $(\rvect,v)$ of problem $(\Lagr_0)$ as in \eqref{1.28}. Here we denote
$(\rvect_\infty,v_\infty):=\Pi_{H^1_\Lagr;\E^1}(\rvect_0,v_0)$ and
$(\rvect_0^0,v_0^0):=\Pi_{H^1_\Lagr;H^1_\Lagrn 0}(\rvect_0,v_0)$. Moreover, we denote the solution of problem
$(\Lagr_0^0)$,  corresponding to initial data $U_0^0:=(\rvect_0^0,v_0^0,\rvect_1,v_1)\in \mathcal{H}^1_\Lagrn 0$, by $(\rvect^0,v^0)\in X^1_\Lagrn 0$.

We apply \cite[Theorem~1.5]{tremodelli}. Consequently, there is a weak solution $(u^0,v^0)$ of $(\Potc)$, unique up to an element of $\C_{X_\Pot}$, such that $-\nabla u^0(0)=\rvect_1$ and $\rho_0 u^0_t(0)=-B\Div \rvect_0^0$. Furthermore, one has $\rvect^0_t=-\nabla u^0$ and, for all $t\in\R$,
$\rvect^0(t)\in H^1(\Omega)^3$ is the unique solution  of the problem
\begin{equation}\label{4.18}
  \begin{cases}
  -B\Div \rvect^0(t)=\rho_0 u^0_t(t),\quad &\text{in $\Omega$,}\\
  \quad \curl\rvect^0(t)=0\quad &\text{in $\Omega$,}\\
  \quad \rvect^0(t)\cdot\boldsymbol{\nu}=0\quad &\text{on $\Gamma_0$,}\\
 \quad \rvect^0(t)\cdot\boldsymbol{\nu}=-v^0(t)\quad &\text{on $\Gamma_1$.}
  \end{cases}
\end{equation}

We now point out that, setting $c^2:=B/\rho_0$, problem $(\Pot)$ reduces to  problem (1.1) in \cite{mugnvit}.
Moreover, the space $\mathcal{H}^1_\Potc$ in \cite[formula (1.6)]{tremodelli}, to which $(u^0(t),v^0(t),u^0_t(t),v^0_t(t))$ belong for $t\in\R$, coincides with the space $\mathcal{H}_1$ in \cite[Chapter 1, formula (1.15)]{mugnvit}.

We can then apply \cite[Theorem~1.3.1]{mugnvit}. Consequently, we deduce that
$$(\nabla u^0(t),u^0_t(t),v^0(t),v^0_t(t))\to 0\quad\text{in $L^2(\Omega)^3\times L^2(\Omega)\times H^1(\Gamma_1)\times L^2(\Gamma_1)$}$$
as $t \to\infty$. Hence, using \cite[Lemma~4.2]{tremodelli},  by \eqref{4.18} we have $\rvect^0(t)\to 0$ in $H^1_{\curl 0}(\Omega)$ as $t\to\infty$. Hence, since $\rvect^0_t=-\nabla u^0$, recalling \eqref{1.6},   we can conclude that
$$(\rvect^0(t),v^0(t),\rvect^0_t(t),v^0_t(t))\to 0\quad\text{in $\mathcal{H}^1_\Lagrn 0$, as $t\to\infty$.}$$
Since $(\rvect(t),v(t))=(\rvect^0(t),v^0(t))+(\rvect_\infty, v_\infty)$ for all $t\in\R$,  \eqref{1.43}  follows.
\end{proof}
\appendix
\section{On the projectors in Theorem~\ref{Theorem 1.3}}\label{appendice A}
This section is devoted to give a more explicit form of the projectors $\Pi_{H^k_\Lagr;H^k_\Lagrn 0}$ and
$\Pi_{H^k_\Lagr;\E^k}$, given in Theorem~\ref{Theorem 1.3} and used in Theorem~\ref{Theorem 1.6}, in the two cases $\kappa\equiv 0$ and $\kappa\not\equiv 0$.

When $\kappa\equiv 0$, since $\svect^\bullet=0$ and $z^\bullet=\mathbbm{1}_{\Gamma_1}$, using \eqref{3.1}, for any $(\svect,z)\in H^k_\Lagr$, we have
$$
 \Pi_{H^k_\Lagr;\E^k}(\svect,z)=\left(\mvect,\dfrac{\int_\Gamma \svect\cdot\nuvect+\int_{\Gamma_1}z}{\mathcal{H}^2(\Gamma_1)}\right),\quad
 \Pi_{H^k_\Lagr;H^k_\Lagrn 0}(\svect,z)=\left(\nvect,z-\dfrac{\int_\Gamma \svect\cdot\nuvect+\int_{\Gamma_1}z}{\mathcal{H}^2(\Gamma_1)}\right),
$$
where $\mvect, \nvect\in H^k(\Omega)^3$ respectively are the unique solutions of the  problems
\begin{align*}
(M)\qquad
&\begin{cases}
\Div \mvect=0&\text{in $\Omega$,}\\
\curl \mvect=\curl\svect&\text{in $\Omega$,}\\
\mvect\cdot\nuvect=\svect\cdot\nuvect&\text{on $\Gamma_0$,}\\
\mvect\cdot\nuvect=\svect\cdot\nuvect+z-\dfrac{\int_\Gamma \svect\cdot\nuvect+\int_{\Gamma_1}z}{\mathcal{H}^2(\Gamma_1)}&\text{on $\Gamma_1$,}
\end{cases}
\\\intertext{and}
(N)\qquad
&\begin{cases}
\Div \nvect=\Div\svect&\text{in $\Omega$,}\\
\curl \nvect=0&\text{in $\Omega$,}\\
\nvect\cdot\nuvect=0&\text{on $\Gamma_0$,}\\
\nvect\cdot\nuvect=-z+\dfrac{\int_\Gamma \svect\cdot\nuvect+\int_{\Gamma_1}z}{\mathcal{H}^2(\Gamma_1)}&\text{on $\Gamma_1$.}
\end{cases}
\end{align*}
When $\kappa\not\equiv 0$, since $\svect^\bullet=\svect^*$, $z^\bullet=z^*$,
using \eqref{1.18} and \eqref{3.1}, for any $(\svect,z)\in H^k_\Lagr$, we have
$$
 \Pi_{H^k_\Lagr;\E^k}(\svect,z)=\left(\pvect,\dfrac{\int_\Gamma \svect\cdot\nuvect+\int_{\Gamma_1}z}{-|\Omega|+\int_{\Gamma_1}z^*}z^*\right),\quad
 \Pi_{H^k_\Lagr;H^k_\Lagrn 0}(\svect,z)=\left(\qvect,z-\dfrac{\int_\Gamma \svect\cdot\nuvect+\int_{\Gamma_1}z}{-|\Omega|+\int_{\Gamma_1}z^*}z^*\right),
$$
where $\pvect, \qvect\in H^k(\Omega)^3$ respectively are the unique solutions of the  problems
\begin{align*}
(P)\qquad
&\begin{cases}
\Div \pvect=\dfrac{\int_\Gamma \svect\cdot\nuvect+\int_{\Gamma_1}z}{|\Omega|-\int_{\Gamma_1}z^*}&\text{in $\Omega$,}\\
\curl \pvect=\curl\svect&\text{in $\Omega$,}\\
\pvect\cdot\nuvect=\svect\cdot\nuvect&\text{on $\Gamma_0$,}\\
\pvect\cdot\nuvect=\svect\cdot\nuvect+z+\dfrac{\int_\Gamma \svect\cdot\nuvect+\int_{\Gamma_1}z}{|\Omega|-\int_{\Gamma_1}z^*}z^*&\text{on $\Gamma_1$,}
\end{cases}
\\\intertext{and}
(Q)\qquad
&\begin{cases}
\Div \qvect=\Div\svect+\dfrac{\int_\Gamma \svect\cdot\nuvect+\int_{\Gamma_1}z}{-|\Omega|+\int_{\Gamma_1}z^*}&\text{in $\Omega$,}\\
\curl \qvect=0&\text{in $\Omega$,}\\
\qvect\cdot\nuvect=0&\text{on $\Gamma_0$,}\\
\qvect\cdot\nuvect=-z+\dfrac{\int_\Gamma \svect\cdot\nuvect+\int_{\Gamma_1}z}{-|\Omega|+\int_{\Gamma_1}z^*}z^*&\text{on $\Gamma_1$.}
\end{cases}
\end{align*}


\def\cprime{$'$}
\providecommand{\bysame}{\leavevmode\hbox to3em{\hrulefill}\thinspace}
\providecommand{\MR}{\relax\ifhmode\unskip\space\fi MR }
\providecommand{\MRhref}[2]{%
  \href{http://www.ams.org/mathscinet-getitem?mr=#1}{#2}
}
\providecommand{\href}[2]{#2}

\end{document}